\documentclass[11pt]{amsart}
\usepackage[T1]{fontenc}
\usepackage{lmodern}
\usepackage{microtype}
\usepackage[letterpaper,textwidth=5.9in,textheight=8.8in,centering,heightrounded]{geometry}
\usepackage{amsmath,amssymb}
\usepackage{needspace}
\usepackage[hidelinks]{hyperref}
\hypersetup{pdftitle={Approximation of sections over real curves},pdfauthor={Peihan Wu},
  pdfsubject={Algebraic geometry},
  pdfkeywords={Rational connectedness, tight approximation, Campana sections, integral Hodge classes},
  bookmarksnumbered=true}
\allowdisplaybreaks[2]
\numberwithin{equation}{section}

\newtheorem{theorem}{Theorem}
\newtheorem{corollary}[theorem]{Corollary}
\newtheorem{lemma}[theorem]{Lemma}
\newtheorem{proposition}[theorem]{Proposition}
\newcommand{\R}{\mathbf{R}}
\newcommand{\C}{\mathbf{C}}
\newcommand{\Pone}{\mathbf{P}^{1}}
\newcommand{\cX}{\mathcal{X}}

\newcommand{\cO}{\mathcal{O}}
\newcommand{\dbar}{\bar{\partial}}
\newcommand{\rank}{\operatorname{rk}}
\newcommand{\Supp}{\operatorname{Supp}}
\newcommand{\Gr}{\operatorname{GL}}
\newcommand{\area}{\operatorname{area}}
\newcommand{\Spec}{\operatorname{Spec}}
\newcommand{\Id}{\operatorname{id}}

\title{Approximation of sections over real curves}
\author{Peihan Wu}
\address{Sichuan University, Chengdu, China}
\email{3119567113@qq.com}
\date{December 3, 2026}

\begin{document}
\begin{abstract}
We study approximation of sections over smooth projective real curves for
proper flat regular models with smooth, geometrically rationally connected
generic fiber. Conjugation-equivariant holomorphic sections defined near
the real locus can be approximated on compact subsets of their domains by
algebraic sections matching prescribed finite jets. The construction
uses holomorphic disk replacement with rational spheres inserted in
conjugate pairs. On smooth projective models it gives polynomial degree
bounds in the inverse approximation error for fixed approximation and
interpolation data. For smooth projective, geometrically rationally
connected varieties over function fields of real curves without real points,
we obtain a degree isomorphism for zero-cycles. We also give a transfer
result for real integral Hodge classes of one-cycles along fibrations
with geometrically rationally connected generic fiber.
\end{abstract}
\maketitle

\section{Introduction}

Benoist and Wittenberg introduced tight approximation over real and complex
curve function fields \cite[Definition~2.1]{BW}. Over \(\R\) it requires
conjugation-equivariant approximation near the whole real locus, together
with interpolation of finite jets. They proved the property for rational
varieties, positive-dimensional quadrics and homogeneous spaces of connected
linear algebraic groups, among other cases
\cite[Theorems~3.1, 4.2, 5.1, 6.1 and 7.4 and Example~4.5]{BW}.
The following theorem treats the real setting of their Question~3.6.

\begin{theorem}\label{thm:main}
Let \(B\) be a smooth projective connected curve over \(\R\), \(F=\R(B)\), and \(X\) a smooth proper geometrically rationally connected variety over \(F\). Let \(p\colon\cX\to B\) be a proper flat model of \(X\) with \(\cX\) regular. Let \(\Omega\subset B(\C)\) be a conjugation-invariant open neighborhood of \(B(\R)\), and let \(u\colon\Omega\to\cX(\C)\) be an equivariant holomorphic section of \(p(\C)\). For every finite set \(b_1,\ldots,b_m\in\Omega\) with \(m\geq 0\), and every integer \(r\geq 0\), there are algebraic sections \(s_j\colon B\to\cX\) over \(\R\) such that
\[
j^r_{b_i}(s_j(\C))=j^r_{b_i}(u)\qquad(1\leq i\leq m,\ j\geq1)
\]
and \(s_j(\C)|_\Omega\) tends to \(u\) in the compact-open \(C^\infty\) topology.
\end{theorem}

Here \(j^r\) denotes restriction to the infinitesimal neighborhood of length
\(r+1\). The domain \(\Omega\) can be disconnected or, when \(B(\R)\) is empty,
empty. The two complex components of a connected real curve, when there are
two, are exchanged by conjugation. The theorem also gives the real-locus
formulation in \cite[Question~3.8]{WittenbergICM}.

Our construction first extends the input to a smooth equivariant section
in the relative smooth locus, preserving holomorphic germs at bad values.
On two conjugate good disks, inserted spheres match the relative
determinant degree of a positive twist and yield an acyclic holomorphic
parameter bundle (Lemma~\ref{lem:positive}).
Disk replacement in the other regions and a small-area contraction
then give the corrected section without assuming an initial
real algebraic section. The relative first-jet families and the analytic
replacement are due to Du, Guo and Xie
\cite[Sections~4--5 and Theorem~7.1]{DGX}, building on Gournay's construction
\cite{Gournay}. We give the correction argument in Section~\ref{sec:analytic};
the paired insertions used to prepare the bundle preserve the prescribed
holomorphic germs. The relative homology calculation adapts the
disk-class calculation of \cite[\S2.1]{HWshort}.

Over complex function fields, Benoist and Wittenberg proved tight
approximation under rational simple connectedness, including at bad
fibers \cite[Theorem~1.4]{BWcomplex}. Du, Guo and Xie treat smooth
projective morphisms over affine complex curves with a continuous lifting
that is holomorphic near the chosen compact set and interpolation points
\cite[Theorem~1.7]{DGX}. For a fixed smooth projective rationally connected
complex target, simple connectivity permits the required continuous
extension from a compact bordered neighborhood and properness extends
the resulting algebraic maps across deleted points. The fixed-target
approximation theorem was previously proved for projective manifolds
birational to algebraically elliptic ones \cite[Corollary~1.8]{FL}.
The general fixed rationally connected target case of
\cite[Question~1.1]{BWcomplex} also follows
from Theorem~\ref{thm:main} by viewing the complex base over \(\R\).
Theorem~\ref{thm:main} appeared in the first version of this preprint
\cite[Theorem~1]{WuRealVone}. The present manuscript revises and expands
the public second version \cite{WuRealVtwo}. In a preprint posted on
27 October 2026, Liu claims qualitative weak approximation at all
places over complex curve function fields \cite{Liu}. Cao and Tian
report a complex-analytic approach to that problem and announce work
of Starr and Tian on weak approximation for separably rationally
connected varieties \cite[Remark~2.3]{CT}. The real theorem here gives
equivariant approximation near the whole real locus, together with
finite jets.

For fixed approximation data, Theorem~\ref{thm:effective} bounds degrees by
counting replacement disks, and Theorem~\ref{thm:campana} extends the estimate
to weak log Fano pairs by logarithmic contact repair. The companion paper \cite{HW} applies the fixed-jet construction
and a residue interpolation argument to obtain a linear effective bound
over complex curve function fields.
Over \(\R\), Theorem~\ref{thm:unbounded-real} gives unbounded minimum
degree for admissible zero-jets at one real bad place.

Section~\ref{sec:apps} derives Lang's \(C_1\) conclusion for real curves
without real points, the Manin--Koll\'ar existence statement and the
local--global principle of P\'al and Szab\'o.

Theorem~\ref{thm:real-hodge-transfer} extends the transfer of the real
integral Hodge conjecture for one-cycles from conic bundles to fibrations
with geometrically rationally connected generic fiber. The proof combines
Theorem~\ref{thm:main} with the integral degree-two comparison of
Lemma~\ref{lem:local-complex-kato-two}.
Corollary~\ref{cor:real-mrc-hodge} applies this to maximal rationally
connected quotients of dimension at most two.
Corollary~\ref{cor:real-fivefold-level} uses the companion's fourth
relative Kato comparison to give the sixteen-square bound for point-free
real fivefolds with an odd-index rationally connected fibration of a
complex birational model. This includes the case of separable rational
connectedness in codimension one.
Corollary~\ref{cor:real-sevenfold-level} gives a bound of sixty-four
in dimension seven when the complexification has no \(2\)-torsion in
third integral homology.
The transfer theorem also gives the uniruled-threefold and rationally
connected cases of \cite[Question~2.16, published version]{BWhodgeI}.
In the case of geometrically rational surface fibrations, we give a direct
construction of horizontal cycles on relative Hilbert schemes, including
along multiple components of singular fibers. The all-dimensional
rationally connected case combines the real approximation theorem with integral
generation by rational curves on the complexification. We give the
required sphere argument within the paper. Section~\ref{sec:real-zero-cycles} proves the degree isomorphism for zero-cycles over real curve fields without real points. It uses equivariant parametric disk replacement to eliminate the kernel left by restriction and norm. Corollary~\ref{cor:zero-cycle-diagonal} describes the integral
diagonal image in the product of the local Chow groups, including at
places of bad reduction. Corollary~\ref{cor:real-chow-one} proves generation of the real one-cycle Chow group by curves with rational geometric components. Section~\ref{sec:further-apps} gives a uniform local degree bound
for even-genus curves near isolated singular limits of real quartic
threefolds.

Section~\ref{sec:real-surface-approx} uses the curve approximation
theorem on a pencil of curves on a real surface without real points.
Theorem~\ref{thm:real-divisorial-approx} gives approximation at a
specified transverse divisor. We place this divisor in a nodal
pencil member and keep the prescribed jet during the fixed Abel
deformation. Corollary~\ref{cor:real-surface-quartics} treats smooth
quartic generic fibers with a model smooth away from finitely many
base points.

\section{Smooth extension and the parameter bundle}\label{sec:prep}

We will work with a projective regular model. Write \(\Sigma=B(\C)\), \(M=\cX(\C)\), and use \(p\) also for the analytic map. Write \(\sigma\) for base conjugation and \(\tau\) for conjugation on \(M\). Let \(M^\circ\) be the submersion locus of \(p\), with holomorphic vertical tangent bundle \(T_p\). These are complex manifolds on each component because the regular schemes are smooth over the perfect field \(\R\). We call a complex vector bundle with an antilinear involution over \(\sigma\) equivariant. A holomorphic equivariant bundle is understood to have an antiholomorphic involution on its total space.

There is a finite invariant set \(Z\subset\Sigma\) such that \(p\) is smooth with connected rationally connected fibers over \(\Sigma\setminus Z\). Indeed, smoothness and geometric connectedness hold over a nonempty algebraic open by properness and the hypotheses on the generic fiber. A very free curve in the geometric generic fiber is defined over a finite extension of the function field. Spread it out over an open of the normalization in that extension. After shrinking the open, upper semicontinuity for its pulled-back tangent bundle twisted by \(\cO_{\Pone}(-2)\) gives very free curves in the fibers. The image of the discarded points in the base is finite. This argument applies to both conjugate components if \(\Sigma\) is disconnected. Good fibers are topologically simply connected \cite[Theorem~5.3.2]{CP}.
Equivalently, one can finish from finiteness of the fundamental group:
the finite universal cover is projective and contains a lift of a very
free rational curve, so both varieties are rationally connected. Their
holomorphic Euler characteristics are one. By Riemann--Roch this
characteristic is multiplied by the degree of a finite \'etale cover,
so the cover has degree one.

\begin{lemma}\label{lem:extension}
Assume \(\Omega\ne\Sigma\), and let \(K\subset\Omega\) be an invariant compact set containing \(\Sigma^\sigma\). There exist a smooth equivariant section \(s_0\colon\Sigma\to M^\circ\), holomorphic near \(K\cup Z\) and equal to \(u\) near \(K\), and a pair of disjoint closed coordinate disks \(Q,\sigma Q\) in \(\Sigma\setminus(K\cup Z\cup\Sigma^\sigma)\). The disks can be taken disjoint from the neighborhoods on which the holomorphic data are prescribed.
\end{lemma}

\begin{proof}
The open complement of \(K\) contains a nonreal point outside \(Z\). Take a sufficiently small disk \(Q\) around it and take its conjugate. Choose an invariant smooth function supported in \(\Omega\) away from these disks that equals one near \(K\). A positive regular superlevel set \(N\) at a level below one is a compact invariant subsurface with boundary disjoint from \(Z\), provided the level avoids the finitely many values of the function on \(Z\). If \(K\) is empty, take \(N=\varnothing\). Prescribe \(u\) on a neighborhood of \(N\). For the elements of \(Z\setminus N\), take disjoint small closed disks away from \(N,Q,\sigma Q\), in conjugate pairs. For one disk in each pair prescribe the restriction of a complex algebraic section on the corresponding component of \(\Sigma\), given by the Graber--Harris--Starr theorem \cite[Theorem~1.1]{GHS}. Use its conjugate on the mate. This gives holomorphic prescriptions on neighborhoods of all these disks. Each prescribed section takes values in \(M^\circ\), since differentiating the section identity gives a right inverse of \(dp\).

Remove the interiors of \(N\) and the bad-value disks. The resulting compact surface with boundary lies over the good values and the action on it is free. By Ehresmann's theorem the restriction is a smooth fiber bundle. Its quotient by conjugation is a smooth fiber bundle over the quotient surface. In an evenly covered base neighborhood, a choice of sheet identifies the quotient bundle with the bundle over that sheet. The prescribed sections define a section on a boundary collar of the quotient. Retain a smaller collar and use a relative triangulation of the complement to extend. Over a simplex the pulled-back bundle is trivial. The extension over vertices uses nonemptiness, over edges path connectedness, and over triangles simple connectivity of the fibers. There are no higher cells. The argument applies to closed components as well.

Smooth the continuous section relative to a smaller collar. For example, an embedding \(j\) of the total space in a Euclidean space gives a fiberwise embedding \(x\mapsto(p(x),j(x))\). A fiberwise tubular neighborhood, obtained from the normal bundles to the fiber embeddings, gives a smooth fiberwise retraction near the image. Relatively smooth the Euclidean coordinate of the continuous section and compose with this retraction. The quotient section lifts uniquely to an equivariant section upstairs because each quotient point has exactly one representative over each base lift. Glue to the prescribed sections along the retained collars.
\end{proof}

For a smooth section \(t\colon\Sigma\to M^\circ\), put \(E_t=t^*T_p\). Endow it with the holomorphic structure given by the \((0,1)\) part of the pullback of a Chern connection on \(T_p\). We choose the connection from an invariant Hermitian metric, so it commutes with conjugation. Integrability on \(\Sigma\) follows from dimension one. Where \(t\) is holomorphic the structure is the natural one, and the structures for two sections coincide wherever the sections coincide. An invariant effective divisor \(D\) defines the natural equivariant bundle \(E_t(-D)\) and the divisor map \(E_t(-D)\to E_t\).

\begin{lemma}\label{lem:positive}
Let \(s_0,Q,\sigma Q\) be as in Lemma~\ref{lem:extension}, and let \(D\) be an invariant effective divisor supported on \(K\). Assume \(n=\rank T_p>0\). There exist a smooth equivariant section \(s\colon\Sigma\to M^\circ\), equal to \(s_0\) off the interiors of the two disks, an equivariant holomorphic bundle \(V\) with \(H^1(\Sigma,V)=0\), and a smooth equivariant complex-linear isomorphism
\[
J\colon V\longrightarrow E_s(-D)
\]
that is holomorphic on a neighborhood of \(K\cup Z\).
\end{lemma}

\begin{proof}
Let \(P\) be the center of \(Q\) and \(Y=p^{-1}(P)\). Trivialize the smooth family over \(Q\) by radial transport for a horizontal distribution agreeing with \(ds_0(T\Sigma)\) along the graph of \(s_0\). Such a distribution is obtained by extending this prescribed right inverse of \(dp\) and using a smooth partition of unity. In the resulting product \(Q\times Y\), the fiber coordinate of \(s_0\) is constant at \(x=s_0(P)\). Under this product identification the vertical tangent remains a complex bundle. A complex-linear connection along the radial segments at fixed fiber coordinate identifies it complex-linearly with the pullback of \(TY\).

Choose a very free \(g\colon\Pone\to Y\) with \(g(\infty)=x\), using \cite[Definition--Theorem~2.1]{Kollar}, and set \(d=\deg(g^*TY)>0\). Write \(L=\cO_\Sigma(P+\sigma P)\). It has positive degree on every complex component of \(\Sigma\). For a positive integer \(k\), put
\[
V=E_{s_0}(-D)\otimes L^{kd}.
\]
By projective GAGA for coherent sheaves and Serre vanishing, \(H^1(\Sigma,V)=0\) for all sufficiently large \(k\). Fix one such \(k\).

Inside \(Q\), insert \(nk\) based spheres in disjoint subdisks, replacing the constant fiber coordinate by \(g\) composed with a smooth degree-one collapse to \(\Pone(\C)\) which is constant at \(\infty\) on a boundary collar. Such a collapse represents the positively oriented quotient of a disk by its boundary and can be made smooth relative to the collar. Make the conjugate modification inside \(\sigma Q\), producing \(s\). The relative change of the determinant degree on \(Q\) is \(nkd\). To compute it, glue the new pullback to the old one using the common collar with the orientation of the old disk reversed. The complex-linear product identification reduces the degree on the resulting sphere to evaluation of \(c_1(TY)\) on the sum of the \(nk\) inserted spheres. Each contributes \(d\). On \(\sigma Q\) the increment is again \(nkd\), because conjugation reverses the base orientation and changes the sign of the first Chern class of the complex bundle.

Off the interiors of the disks, use the canonical nonvanishing section of \(L^{kd}\) to identify \(V\) with \(E_s(-D)\). Its extension over \(Q\) is a relative clutching problem. In smooth interior frames of the two bundles the prescribed collar isomorphism is a loop in \(\Gr_n(\C)\) up to homotopy. Its determinant winding is the difference between the relative determinant increments of the two bundles compared with \(E_{s_0}(-D)\). The twist has increment \(nkd\), since \(L\) has a simple canonical zero at \(P\). The modified pullback has the same increment by the sphere calculation. The winding is therefore zero. The determinant induces an isomorphism \(\pi_1(\Gr_n(\C))\simeq\pi_1(\C^\times)\), so the collar isomorphism extends smoothly over \(Q\) while retaining a smaller collar. Use its conjugate over the other disk. Near \(K\cup Z\), the identification is the canonical one between the unmodified pullbacks after trivializing the twist and is holomorphic there.
\end{proof}

\section{Equivariant compact approximation}\label{sec:analytic}

We first recall the local geometric input. If \(p\) is smooth projective over an open set with rationally connected fibers, then around any triple \((z_0,x_0,\beta_0)\), with \(x_0\in p^{-1}(z_0)\) and \(\beta_0\in T_{p,x_0}\), there is a jointly holomorphic family
\[
G(w;z,x,\beta),\qquad w\in\Pone(\C),
\]
taking values in \(p^{-1}(z)\), with value \(x\) and derivative \(\beta\) at \(w=0\). This includes \(\beta_0=0\) \cite[Proposition~3.5]{DGX}. The finite-jet theorem of Koll\'ar \cite[Definition--Theorem~2.1(5)]{Kollar} gives a curve \(g\) with the specified length-two jet and \(H^1(\Pone,g^*T_{p^{-1}(z_0)}(-2))=0\). This cohomology group is the obstruction space for the relative restriction from maps on \(\Pone\) to their prescribed restriction to \(2[0]\). Indeed local infinitesimal lifts exist by smoothness and, after matching the restriction to \(2[0]\), their differences are sections of the tangent bundle twisted by \(-2[0]\) with coefficients in the square-zero ideal. Vanishing removes the \v Cech obstruction. Thus relative jet evaluation on the local Hom space is smooth at \(g\). A local analytic section of that evaluation and the universal map give \(G\), with the base \(z\) included among the parameters \cite[Theorem~3.2 and Proposition~3.4]{DGX}.

The following is the conjugation-equivariant form of the compact-source
approximation theorem of Du, Guo and Xie \cite[Theorem~7.1]{DGX}.
The proof records the conjugate choices and the estimates needed for the
degree bound and the parametric applications below.

\begin{proposition}\label{prop:analytic}
Let \(s\colon\Sigma\to M^\circ\) be a smooth equivariant section. Let \(K\subset\Sigma\) be invariant compact and contain \(\Sigma^\sigma\), and let \(D\) be an invariant effective divisor supported on \(K\). Suppose \(s\) is holomorphic on an invariant open \(U_h\) containing \(K\cup Z\). Suppose there are an equivariant holomorphic vector bundle \(V\) with \(H^1(\Sigma,V)=0\) and a smooth equivariant complex-linear isomorphism \(J\colon V\to E_s(-D)\) holomorphic on \(U_h\). Then \(s|_K\) can be approximated arbitrarily closely and uniformly by global equivariant holomorphic sections agreeing with \(s\) on \(D\).
\end{proposition}

\begin{proof}
\emph{The normalized equation.}
Near the compact graph of \(s\) we use an equivariant smooth local addition \(\Psi(x,v)\) along the fibers of \(p\). It is holomorphic in \(v\in T_{p,x}\), has value \(x\) at \(v=0\), and has vertical derivative \(\Id\) there. To construct it, use relative holomorphic coordinates near the graph. Near the relative diagonal, the difference of two fiber coordinates at \(x'\) and \(x\) can be identified with an element of \(T_{p,x}\) by the inverse coordinate derivative at \(x\). Take a partition-of-unity average with weights depending only on \(x\). For fixed \(x\) this map is holomorphic in \(x'\) and its derivative on the diagonal is the identity. On a common invariant neighborhood average it with its conjugate, using the antiholomorphic involutions on the target fibers and \(T_p\). This preserves holomorphy in \(x'\) and the derivative identity. Invert in that variable. Compactness gives a common neighborhood for the addition along the graph \cite[Lemma~4.1]{DGX}.

Write \(\iota\colon E_s(-D)\to E_s\) for the divisor map. On a small metric ball bundle of \(V\) put
\begin{equation}\label{eq:phi}
\Phi(z,v)=\Psi\bigl(s(z),\iota_zJ_zv\bigr).
\end{equation}
In a holomorphic frame of \(V\) and relative holomorphic target coordinates let \(\widehat\Phi(z,q)\) denote its fiber coordinate. Off \(D\), normalize the section equation to
\begin{equation}\label{eq:norm}
\dbar_V v+\mathcal A(z,v)=0,\qquad
\mathcal A(z,q)=
(\partial_q\widehat\Phi)^{-1}\partial_{\bar z}\widehat\Phi\,d\bar z .
\end{equation}
Here the base derivative is taken at fixed \(q\). The fiber map in \(\Psi\) has uniformly invertible derivative on a small ball, so the derivative in \eqref{eq:norm} is invertible off \(D\) for one fixed small ball of \(V\). Holomorphic target coordinate changes cancel from the normalized coefficient and holomorphic frame transitions give the transformation law of a \(V\)-valued \((0,1)\)-form.

In a holomorphic frame of \(E_s\) over \(U_h\), the coordinate expression for \(\Psi(s(z),w)\) and its fiber derivative at zero are holomorphic in \(z\) at \(w=0\). The coordinate expression at zero is that of \(s\) and its derivative is the coordinate matrix of the holomorphic frame. Hence the corresponding normalized coefficient \(C(z,w)\) has zero constant and fiber-linear terms. Near a point of \(D\) of multiplicity \(m\), choose holomorphic frames in which the argument of \(\Psi\) in \eqref{eq:phi} is \(z^m L(z)q\), with \(L\) holomorphic invertible. In those frames \(\mathcal A\) is
\[
(z^m)^{-1} L(z)^{-1}C(z,z^m L(z)q).
\]
Holomorphic dependence on \(w\) and vanishing to order two give a smooth extension through \(z=0\), holomorphic in \(q\). For instance Taylor's integral formula factors out \(z^{2m}\) from the numerator. Thus \(\mathcal A\) extends over \(\Sigma\), is equivariant, and has zero constant and linear terms on \(U_h\).

Let
\[
U=\{z\in\Sigma:\ |\mathcal A(z,0)|+|d_q\mathcal A(z,0)|>0\},
\]
using invariant smooth metrics throughout. Its closure is compact in the good locus and misses \(U_h\). In particular it misses the real locus and the divisor. When \(U\) is empty \(s\) is already holomorphic by \eqref{eq:norm} at \(v=0\) off the divisor and by its prescribed holomorphy near the divisor.

\medskip
\emph{Local replacement.}
Cover \(\overline U\) by finitely many coordinate neighborhoods with relative target charts and holomorphic frames of \(V\), using chart pairs from evenly covered neighborhoods in the free quotient. Write \(\zeta=z-c\) in the chart at a center \(c\). The dependence of \(\widehat\Phi\) on \(q\) is holomorphic, so its expansion at \(c\) has the form
\[
\widehat\Phi(z,q)=m(q)+a(q)\zeta+b(q)\bar\zeta+O(|\zeta|^2)
\]
with coefficients holomorphic in \(q\). Take the chart and the parameter ball small enough to use one family \(G\) for the nearby jets. For radii \(0<\rho\leq\theta<1/4\), put \(h=\theta\rho\), and take the annulus \(\rho-h\leq|\zeta|\leq\rho+h\). The prepared family on the disk is
\[
H(z,q)=G\bigl(\rho^2/\zeta;z,m(q)+a(q)\zeta,b(q)\bigr).
\]
Here point and tangent coordinates refer to the fiber over \(z\). The expression is holomorphic at \(\zeta=0\) by the sphere coordinate \(1/w=\zeta/\rho^2\). In the annulus, \(|\rho^2/\zeta|\leq 2\rho\) and \(|\bar\zeta-\rho^2/\zeta|\leq 3h\). Taylor expansion of \(G\) at \(w=0\) gives, in the common target chart,
\begin{equation}\label{eq:ann}
\bigl\|\widehat H-\widehat\Phi\bigr\|_{C_q^3}
\leq C(h+\rho^2)\leq 2Ch .
\end{equation}
For the derivative estimates, first bound the difference on a slightly larger parameter ball and use the Cauchy estimates in \(q\). All bounds may be made uniform over the centers. In each of finitely many neighborhoods the jets vary in a fixed compact subset of the domain of the chosen family \(G\). The target fiber coordinate values used on the annuli lie in a convex ball with positive margin inside the chart. Such margins follow from \(G(0;z,x,\beta)=x\) and the choice of small charts and radii.

In the annulus use the fiber coordinate
\[
\Xi=(1-\chi)\widehat\Phi+\chi\widehat H
\]
where the radial cutoff \(\chi\) is one on an inner collar, zero on an outer collar, and \(|d\chi|\leq C/h\). By \eqref{eq:ann} and the invertibility of \(\partial_q\widehat\Phi\) on the good charts, \(\partial_q\Xi\) is uniformly invertible on a smaller fixed ball. Since \(\partial_{\bar z}\widehat H=0\), the factor \(h^{-1}\) from the cutoff is multiplied by an \(O(h)\) mismatch. Consequently
\[
\mathcal A_{\rm ann}=(\partial_q\Xi)^{-1}\partial_{\bar z}\Xi\,d\bar z
\]
is uniformly bounded on the parameter ball. It is holomorphic in \(q\), so its first two fiber derivatives are uniformly bounded after shrinking the ball once. In the two collars it agrees with \(\mathcal A\) and zero, respectively, in the outer and inner collars. In the core we use \(H\) and the equation \(\dbar_V v=0\). Joint holomorphy of \(H\) makes that equation sufficient even if the fiber derivative of \(H\) is singular.

\medskip
\emph{Packing and correction.}
Choose finitely many disjoint outer disks of this kind inside \(U\), in conjugate pairs, leaving uncovered area \(O(\theta)\) in \(U\). For this, triangulate a compact neighborhood of the image of \(\overline U\) in the free quotient, subordinate to the evenly covered coordinate neighborhoods. Inside each triangle use a chosen lift and the Vitali covering by arbitrarily small coordinate disks of closure in \(U\) and in the lifted triangle interior. A finite truncation covers up to arbitrarily small area. Shrink the chosen finite collection slightly if needed to separate their closures, retaining the area bound. Lift to both sheets. The disk packing has \(\sum_i\rho_i^2\leq C\) by disjointness and the finite-chart area comparisons. Thus the union of all annuli has area at most \(C\sum_i\rho_i h_i=O(\theta)\). The charts, frames, and replacements in paired disks are obtained by conjugating the source and target coordinates, the frame, and the family \(G\). The resulting charts and frames are holomorphic. The radial cutoff and the expression \(\rho_i^2/\zeta\) respect the paired coordinates.

Use the new fiber coordinate family on the annuli, \(H\) on the cores, and \(\Phi\) off the disks. They agree on open collars and give one smooth equivariant family \(\Phi_\theta\) on a fixed smaller ball bundle of \(V\). The corresponding coefficients \(\mathcal A_{\rm ann}\), zero, and \(\mathcal A\) glue on the same collars to \(\mathcal A_\theta\), which is uniformly bounded through its second fiber derivative. The constant and fiber-linear terms of \(\mathcal A_\theta\) vanish pointwise outside the union \(Z_\theta\) of the annuli and the uncovered part of \(U\), with \(\area(Z_\theta)=O(\theta)\). The replacement leaves \(\Phi\) unchanged on \(U_h\).

By vanishing of \(H^1(\Sigma,V)\) and the elliptic estimate for \(\dbar_V\), choose a bounded right inverse
\[
Q\colon L^4(\Sigma,\Lambda^{0,1}\otimes V)\longrightarrow W^{1,4}(\Sigma,V).
\]
Sobolev embedding bounds it into \(C^0\). These norms and the inverse are fixed independently of the packing. Conjugation acts on sections and on \(V\)-valued \((0,1)\)-forms, and the two actions intertwine \(\dbar_V\). In conjugate charts this is conjugation of the coefficients evaluated in the conjugate chart. Average \(Q\) with its conjugate, obtaining an equivariant right inverse.

On the ball of radius \(\delta\) in the space of continuous equivariant sections of \(V\) consider
\[
\mathcal T_\theta(v)=-Q\mathcal A_\theta(\cdot,v(\cdot)).
\]
The uniform bound for the constant term and the area bound for its nonzero set give
\[
\|\mathcal T_\theta(0)\|_\infty\leq C\theta^{1/4}.
\]
On \(Z_\theta\) the fiber derivative is uniformly bounded. On the complement its value at the zero parameter vanishes, so its bound on the radius-\(\delta\) ball is \(C\delta\), by the uniform second derivative bound. Integrating the derivative along fiber segments and using the \(L^4\) bound for \(Q\) therefore gives
\[
\|\mathcal T_\theta(v)-\mathcal T_\theta(v')\|_\infty
\leq C(\theta^{1/4}+\delta)\|v-v'\|_\infty .
\]
Choose \(\delta\) first and then \(\theta\) so that the Lipschitz constant is at most \(1/2\) and \(\|\mathcal T_\theta(0)\|_\infty\leq\delta/2\). This gives a contraction into the ball. Its equivariant fixed point is in \(W^{1,4}\), satisfies \(\|v_\theta\|_\infty=O(\theta^{1/4})\), and solves \(\dbar_Vv_\theta+\mathcal A_\theta(\cdot,v_\theta)=0\).

\medskip
\emph{Reconstruction.}
Set \(F_\theta(z)=\Phi_\theta(z,v_\theta(z))\). It is a continuous equivariant section. Away from the divisor and the cores the normalized equation and the Sobolev chain rule show that its coordinate expressions satisfy \(\dbar F_\theta=0\) weakly. On each core the parameter \(q\) satisfies \(\dbar q=0\), so the same conclusion follows from joint holomorphy of \(H\). The formulas agree on the collars, and local source derivatives need only be bounded for a fixed finite packing. The Weyl lemma in target coordinates gives a holomorphic section off \(D\), and continuity gives removable singularities at its support. On \(K\) the section is given by \eqref{eq:phi}, which tends uniformly to \(s\). In a local relative chart at a point of \(D\) of multiplicity \(m\), \eqref{eq:phi} has coordinate displacement \(O(|z|^m)\) from \(s\), by the divisor factor and boundedness of \(v_\theta\). Holomorphy gives equality on \(D\). This proves the proposition.
\end{proof}

\section{Algebraization and the proof of Theorem~\ref{thm:main}}

The model in the theorem is integral. Its regular irreducible components are disjoint and each dominates \(B\) by flatness and openness, while the generic fiber is integral. By Chow's lemma over \(B\) and a projective resolution in characteristic zero, there is an integral regular projective model birational to it over \(B\). The projective model is flat because its local rings are torsion-free over the local discrete valuation rings of \(B\). Its generic fiber is smooth proper and remains geometrically rationally connected by birational invariance. The tight approximation property transfers between these models by \cite[Theorem~3.1]{BW}, so it suffices to prove the result for all data on projective regular models.

If the relative dimension is zero, the geometrically integral generic fiber is \(\Spec F\). The generic section extends over \(B\) by properness. Over a dense open the morphism is an isomorphism, so the extended section agrees with \(u\) on a dense subset of each component of \(\Omega\), hence everywhere. If \(\Omega=\Sigma\), algebraize \(u\) directly by projective GAGA and descend it by equivariance.

In the remaining case, enlarge the finite interpolation set to its conjugation orbit \(A\), and put \(D=(r+1)\sum_{a\in A}[a]\). Fix an invariant compact \(K\subset\Omega\) containing \(B(\R)\cup A\). Lemmas~\ref{lem:extension} and~\ref{lem:positive} give a section \(s\) and an isomorphism \(J\) as in Proposition~\ref{prop:analytic}, with \(s=u\) near \(K\). Choose an invariant open \(U_h\) containing \(K\cup Z\) on which both \(s\) and \(J\) are holomorphic. The proposition gives global equivariant holomorphic sections arbitrarily close to \(u\) on \(K\) and agreeing on \(D\). By projective GAGA \cite{Serre}, each is an algebraic section over \(\C\). Full faithfulness and conjugation equivariance give descent of this morphism to \(\R\). This argument also applies to a pair of conjugate complex components, as in \cite[Section~2.2.1]{BW}.

Choose an invariant compact exhaustion \(K_j\subset\operatorname{int}_\Omega K_{j+1}\) of \(\Omega\) with every \(K_j\) containing \(B(\R)\cup A\), and choose approximation error \(1/j\) on \(K_j\), using the same divisor \(D\). The resulting algebraic sections converge uniformly on compact subsets. On a compact subset take finitely many source disks with larger concentric closed disks in \(\Omega\) whose images under \(u\) lie in target charts. The approximating maps eventually lie in the same charts on the larger disks. Cauchy estimates on the smaller disks give convergence of all derivatives. If \(\Omega\) is empty, take \(K=A=\varnothing\) in the construction and repeat the resulting algebraic section. This proves Theorem~\ref{thm:main}.

\Needspace{17\baselineskip}
\section{A polynomial degree bound for tight approximation}\label{sec:effective-tight}

The approximation theorems of Benoist--Wittenberg and Du--Guo--Xie are qualitative \cite[Theorem~1.4]{BWcomplex}\cite[Theorem~1.7]{DGX}. Here we count the disks in the construction to obtain a degree bound for fixed holomorphic and interpolation data, including on models with bad fibers.

\begin{theorem}\label{thm:effective}
Let \(k=\R\) or \(\C\), let \(B/k\) be a smooth projective geometrically integral curve, and let \(p\colon\cX\to B\) be flat and projective with \(\cX\) smooth projective and with smooth geometrically rationally connected generic fiber. Fix a very ample line bundle \(H\) on \(\cX\). Fix a holomorphic section \(u\) over an open subset \(\Omega\subset B(\C)\), a compact set \(K\subset\Omega\), an integer \(m\geq0\), and finitely many jets of \(u\) on \(\Omega\). Over \(\R\), require equivariance and that \(\Omega\) be invariant and contain \(B(\R)\).

There exist \(C,A>0\) depending on the fixed data such that, for every \(0<\varepsilon<1\), there is a \(k\)-algebraic section \(\sigma_\varepsilon\) realizing the prescribed jets and satisfying
\begin{equation}\label{eq:effective-tight}
 d_{C^m(K)}(\sigma_\varepsilon,u)<\varepsilon,\qquad
 \deg_B(\sigma_\varepsilon^*H)\leq C\varepsilon^{-A}.
\end{equation}
The prescribed jets may lie over singular fibers.
\end{theorem}

\begin{lemma}\label{lem:quant-packing}
On a finite union of coordinate triangles in a compact Riemann surface, with disjoint interiors, one can find \(O(\theta^{-b})\) disjoint closed coordinate disks of radii at most \(\theta\) leaving area \(O(\theta)\), for some \(b>0\) and \(0<\theta<1\).
\end{lemma}

\begin{proof}
Work in one triangle and use Euclidean area in its coordinates. Use grids of side \(h_j=\theta q^j\) for a fixed \(0<q<1\). At stage \(j\) insert the concentric disk of radius \(h_j/3\) in each square whose closure lies in the open area left over. If \(a_j\) is the available Euclidean area at the start of this stage and \(N_j\) is the number of chosen squares, then \(N_jh_j^2\le a_j\). Up to a set of area zero, the part of the remaining region in omitted squares lies within \(\sqrt2h_j\) of the boundary of the triangle or of a previous disk. The former strip has area \(O(h_j)\), while the strip about a disk from stage \(i<j\) has area \(O(h_ih_j)\) since \(h_j<h_i\). As the disks occupy a fraction \(\kappa=\pi/9\) of the chosen squares, there is a constant \(C_0\) such that
\[
a_{j+1}\le (1-\kappa)a_j+C_0\theta q^j+
   C_0\sum_{i<j}a_i q^{j-i}.
\]
Fix \(1-\kappa<\tau<1\), and then choose \(q<\tau\) with \(C_0q/(\tau-q)\le(\tau-1+\kappa)/2\). Induction gives \(a_j\le M\tau^j\) by using \(M\ge a_0\) and \(C_0/M\le(\tau-1+\kappa)/2\). After stage \(J=\lceil\log(1/\theta)/\log(1/\tau)\rceil\), the remaining area is \(O(\theta)\) for the coordinate area and hence for any fixed smooth area form. The number of disks is bounded by \(a_0\sum_{j=0}^J h_j^{-2}=O(\theta^{-b})\) with \(b=2(1+\log(1/q)/\log(1/\tau))\). Apply this to each triangle. On an evenly covered triangle of a free conjugation quotient, choose the disks in one lift and take their conjugates in the other.
\end{proof}

\begin{proof}[Proof of Theorem~\ref{thm:effective}]
If \(\Omega=B(\C)\), apply GAGA to \(u\). If the relative dimension is zero, use the section of fixed degree obtained by properness as in the proof of Theorem~\ref{thm:main}. Otherwise enlarge \(K\) within \(\Omega\) to a compact \(K^+\) whose interior contains \(K\) and the interpolation points. In the real case include \(B(\R)\) and take the invariant enlargement and divisor of prescribed jets. Choose \(K^+\) to contain the larger disks of a fixed finite system of Cauchy estimates at \(K\). Apply Lemmas~\ref{lem:extension} and~\ref{lem:positive} to obtain the fixed smooth section \(s\), holomorphic on a neighborhood of \(K^+\cup Z\), and the holomorphic parameter bundle \(V\) with \(H^1(\Sigma,V)=0\). Over \(\C\) the same preparation omits the involutions: extension on the unprotected surface uses simply connected fibers, and the compensating spheres and positive twist are placed in a single good disk.

In Proposition~\ref{prop:analytic}, the unmodified constant and linear coefficients vanish near \(K^+\cup Z\). Place the closure of the set where either is nonzero in the interior of a finite union of coordinate triangles with closures disjoint from \(K^+\cup Z\). Thus a smaller neighborhood of the protected data is left unmodified. Subdivide so that each triangle carries the fixed chart, frame and jet-family data of that proposition. In the real case triangulate in the free quotient. Use Lemma~\ref{lem:quant-packing} in place of the qualitative packing. If the outer radius of disk \(i\) is \(\rho_i(1+\theta)\), use annular half-width \(\theta\rho_i\). The number of disks is \(O(\theta^{-b})\), the uncovered area is \(O(\theta)\), and \(\sum_i\rho_i^2=O(1)\) gives total annular area \(O(\theta)\). If both coefficients vanish identically, take the packing to be empty.

The \(O(\theta)\) defect area and the unchanged annular estimates give
the hypotheses of the contraction in Proposition~\ref{prop:analytic}.
Applying it with the new packing gives
\(\|v_\theta\|_\infty=O(\theta^{1/4})\) and a holomorphic reconstruction
\(F_\theta\) preserving the jets. On \(K^+\) the local addition is
unmodified. The Cauchy estimates used in the proof of
Theorem~\ref{thm:main}, now on the fixed pairs of disks in \(K^+\),
give \(d_{C^m(K)}(F_\theta,u)=O(\theta^{1/4})\).

To bound the degree, the homotopy \(\Phi_\theta(z,tv_\theta(z))\), \(0\le t\le1\), identifies the homology class of \(F_\theta\) with that of the grafted section at parameter zero. The original and grafted maps agree on the exterior collars. In each doubled disk, the original map and the annular interpolation lie in a contractible relative target chart \(W_i\). In the core expression, homotope \(z=c_i+\zeta\) to \(c_i\) in the base parameter of the rational curve family, and homotope the term \(a(0)\zeta\) to zero. The core boundary stays in \(W_i\), since the sphere coordinate on it is \(w=\rho_i^2/\zeta\) with \(|w|\le2\rho_i\). The resulting core is the vertical sphere \(g_i\) over \(c_i\), minus a cap lying in \(W_i\). The coordinate change preserves complex orientation, and \(H_2(M,\mathbb Z)\to H_2(M,W_i;\mathbb Z)\) is an isomorphism. Thus the doubled disk contributes \(+(g_i)_*[\Pone(\C)]\) in integral homology. Summing and evaluating \(c_1(H)\) gives
\[
 \deg_B(\sigma_\theta^*H)
 =\big\langle c_1(H),s_*[\Sigma]\big\rangle+\sum_i\deg(g_i^*H)
 =O(1+\theta^{-b}),
\]
where \(\sigma_\theta\) is obtained from \(F_\theta\) by GAGA and, over \(\R\), descent. The sphere degrees are bounded because the finitely many fixed local families have locally constant degrees. The real degree equals the degree after complexification. Take \(\theta=c\varepsilon^4\) for \(c>0\) fixed and sufficiently small. This gives~\eqref{eq:effective-tight} with \(A=4b\).
\end{proof}

\section{Uniform interpolation bounds and a real obstruction}\label{sec:effective-weak}

Keep the projective model and \(H\) of Theorem~\ref{thm:effective}. Fix a finite set \(S\) of closed points of \(B\). At \(b\in S\), write \(L_b=r_b+1\) for the length of the prescribed jet and set
\[
 N=\sum_{b\in S}[k(b):k]L_b.
\]
A family of jets is admissible if it is represented by local holomorphic sections at the specified points. Over \(\R\), admissibility requires these germs, with their conjugates, to be represented by a single equivariant holomorphic section on an invariant neighborhood of \(B(\R)\) and of the specified points. At the center of any such section, \(p\) is smooth because its differential admits a right inverse. In particular the jets may be based in the relative smooth locus of a singular fiber.

Hassett and Tschinkel studied effective weak approximation for quadric surface fibrations over finite fields, with height bounds for prescribed fiber values \cite[Section~10]{HTquadric}. Tanimoto and Tschinkel place this question in the study of spaces of sections matching jets, and establish homological stability for projective bundles, conic bundles and quadric surface bundles over complex curves \cite[Introduction]{TTstability}.

\begin{proposition}\label{prop:fixed-order}
Over \(\C\), fixing the entire vector of orders \((r_b)_{b\in S}\) gives a uniform finite degree bound for all admissible jets. Over \(\R\), fixing that vector and a compact set of admissible jets gives a uniform finite degree bound on that set.
\end{proposition}

\begin{proof}
Over \(\C\), let \(J\) be the finite-type scheme of jets centered in the relative smooth locus at \(S\). Every point of \(J(\C)\) is admissible, by relative holomorphic coordinates. The scheme of sections of bounded \(H\)-degree is of finite type, by the Hilbert polynomials of their graphs, and the image of its jet evaluation map is constructible. Theorem~\ref{thm:effective} makes the increasing union of these images all of \(J(\C)\). Uncountability implies that some image contains a dense open on each irreducible component, since a complex irreducible variety is not a countable union of proper closed subsets. Apply the same argument to the complement, using induction on dimension and then taking the maximum of the finitely many bounds.

Over \(\R\), choose for each admissible jet a section
\(f=\Phi_\theta(\cdot,v_\theta)\) as in Proposition~\ref{prop:analytic},
with \(\theta\) small and interpolation divisor \(D=\sum_b L_b[b]\).
For each fixed packing the coefficient is smooth, so elliptic
bootstrapping in the semilinear equation makes \(v_\theta\) smooth.
Set \(B_\theta=D_q\mathcal A_\theta(\cdot,v_\theta)\) and
\(\mathcal D=\dbar_V+B_\theta\). On the defect set of area
\(O(\theta)\), the coefficient \(B_\theta\) is bounded. On the
complement it is \(O(\theta^{1/4})\), by the second derivative bound
and \(\|v_\theta\|_\infty=O(\theta^{1/4})\). Thus multiplication by
\(B_\theta\) from \(C^0\) to \(L^4\) has norm \(O(\theta^{1/4})\).
Since \(\dbar_V Q=I\), the operator \(Q(I+B_\theta Q)^{-1}\) is a
right inverse to \(\mathcal D\) for small \(\theta\). In particular
\(H^1(\Sigma,V_{\mathcal D})=0\) for the holomorphic structure defined
by \(\mathcal D\).

In a target chart differentiate
\(\partial_{\bar z}\widehat\Phi_\theta
   =(D_q\widehat\Phi_\theta)\mathcal A_\theta\) at fixed \(q\).
Evaluating at \(q=v_\theta(z)\), the second fiber derivative term
cancels the term from \(\partial_{\bar z}v_\theta
   =-\mathcal A_\theta(z,v_\theta)\). This gives
\[
 \dbar_{f^*T_p}\bigl((D_q\Phi_\theta)w\bigr)
       =(D_q\Phi_\theta)\mathcal D w .
\]
On the cores the identity follows from joint holomorphy and
\(\partial_{\bar z}v_\theta=0\). Near \(D\), the formula
\eqref{eq:phi} has the divisor factor in its fiber derivative, so
dividing by the local equation of \(D\) gives a smooth extension of
the divided map. The identity makes it holomorphic off \(D\), hence
also across \(D\). This gives a holomorphic map
\(V_{\mathcal D}\to f^*T_p(-D)\) that is generically an isomorphism
on every component, by invertibility on the outer collars or, in the
absence of replacements, in the original relative charts off \(D\).
Its cokernel has finite support, so \(H^1(\Sigma,f^*T_p(-D))=0\).
If \(\Omega=\Sigma\), restrict \(u\) to a smaller admissible
neighborhood before the construction. In relative dimension zero the
unique section suffices.

The vanishing makes the algebraic jet evaluation map on the scheme of sections smooth at \(f\): infinitesimal local lifts are available by relative smoothness, and their gluing obstruction with the jets fixed lies in this \(H^1\). Its real-point map is thus open near \(f\), with degree locally constant on the section space. A finite subcover of the compact set gives the claimed bound.
\end{proof}

\begin{theorem}\label{thm:unbounded-real}
There exists a smooth projective geometrically integral model over
\(B=\Pone_\R\) with smooth Fano generic fiber of dimension three,
a real place \(b\) of bad reduction, and admissible zero-jets \(a_j\)
at \(b\), each realizable by a real algebraic section, such that
for every ample \(H\) the minimum \(H\)-degree of an interpolating
section tends to infinity.
\end{theorem}

\begin{proof}
Use the special fiber of the separate geometrically rational
example in \cite{WuUnbounded}. In coordinates \([x:y:z:u:v:w]\)
on \(\mathbf P^5_\R\), it is the complete intersection
\[
 Q_0=u^2+v^2+w^2-z^2=0,\qquad F_0=x^3+y^3+z^3=0 .
\]
Projection onto \([x:y:z]\) is defined at every real point of
this fiber because its center meets the fiber in a conic without
real points. In the elliptic curve \(F_0=0\) in \(\mathbf P^2_\R\),
the point \([1:-1:0]\) has exactly one real point above it in the
special fiber, namely \(a_0=[1:-1:0:0:0:0]\). For a morphism from \(\Pone_\C\)
whose image contains \(a_0\), the projection extends across the
indeterminacy points by properness and is constant by
Riemann--Hurwitz. Its image therefore contains no other real point.

The special fiber contains the path
\[
 a(s)=[1:-(1+s^3)^{1/3}:s:s:0:0],\qquad 0<s\leq 1.
\]
Form the pencils \(Q_0+tQ_1\) and \(F_0+tF_1\) with general real
forms \(Q_1,F_1\) of degrees two and three vanishing at \(a(1)\)
and with \(Q_1(a_0)\ne0\). Bertini away from \(a(1)\) and independence
of the general differentials there give smoothness of the geometric
generic complete intersection. Adjunction makes it Fano. Resolve its
projective closure over \(B\), keeping the smooth locus unchanged.
The resulting integral model dominates the regular curve \(B\) and
is therefore flat over it.
The fiber differentials of \(Q_0,F_0\) are independent along \(a(s)\).
At \(a_0\) the first has zero fiber differential and nonzero total
differential \(Q_1(a_0)\,dt\), while the second has nonzero fiber
differential. The total space is smooth at all these points.
They belong to the horizontal closure because the local total
dimension is four and the central complete intersection has dimension
three. Hence the resolution is an isomorphism there. The differential
of the resulting projection vanishes at the lift of \(a_0\).

The constant section with value \(a(1)\) lifts to the resolution
by properness. For each \(0<s\leq1\), a time-dependent vertical vector
field supported in the relative smooth locus transports this real
smooth section to one with value \(a(s)\) at \(t=0\). Theorem~\ref{thm:main}
and \cite[Proposition~2.3 and Corollary~2.4]{BW} interpolate that value
by an algebraic section. Set \(a_j=a(1/j)\) on the resolved model.

Suppose infinitely many of these values were realized by sections
of bounded degree for a very ample line bundle. After passing to
a common degree and taking a limit in the projective coarse moduli
space of one-pointed genus-zero stable maps \cite[\S1.2]{FP}, apply
the obstruction of \cite{WuUnbounded} as follows. The marked value
of the limit is the lift of \(a_0\). The degree over \(B_\C\) remains
one, giving a unique degree-one component and a section on that
component. The limit is fixed by conjugation in the coarse moduli
space. An isomorphism with the conjugate stable map fixes the
marked and degree-one vertices of the dual tree and every edge on
the path between them. The images of the nodes on this path are
therefore real. Starting at the marked value, project its vertical
components to the original special fiber. The elliptic projection
property forces the image of each successive real node to be \(a_0\).
Since the resolution is an isomorphism there, the degree-one
component passes through the lift of \(a_0\). This contradicts the
vanishing of the differential of the projection at that point.
Passing to a very ample power gives the assertion for every ample
line bundle.
\end{proof}

The obstruction disproves any bound depending only on \((p,H,S)\)
even with \(S=\{b\}\) and \(N=1\), although each jet is realized.

\Needspace{18\baselineskip}
\section{Rational points and approximation of real curves}\label{sec:apps}

\subsection{Rational points and the reciprocity obstruction}

\begin{corollary}\label{cor:points}
Suppose \(B(\R)=\varnothing\), and put \(F=\R(B)\). Every smooth proper geometrically rationally connected variety over \(F\) has an \(F\)-point and satisfies weak approximation at every closed point of \(B\). Moreover \(F\) is a \(C_1\) field.
\end{corollary}

\begin{proof}
After shrinking to a nonempty open of \(B\), spread the given variety out to a smooth integral scheme over that open. Take its schematic closure in a Nagata compactification over \(B\), and resolve it by a morphism that is an isomorphism over the regular locus. The generic fiber is unchanged, and flatness follows from torsion-freeness over the local discrete valuation rings of \(B\). Apply Theorem~\ref{thm:main} with \(\Omega=\varnothing\) and no jets. An algebraic section gives the desired \(F\)-point. For weak approximation apply \cite[Corollary~2.10]{BW}. When \(B(\R)=\varnothing\), the compatibility condition at the real points of the base is vacuous.

Thus \(F\) is geometrically \(C_1\) in the sense of \cite[Definition~1.4]{HX}. Hogadi and Xu prove that a geometrically \(C_1\) field of characteristic zero is \(C_1\) \cite[Corollary~1.5]{HX}. Their result includes singular hypersurfaces, so every homogeneous form of positive degree \(d\) in more than \(d\) variables over \(F\) has a nonzero zero.
\end{proof}

The existence assertion is the Manin--Koll\'ar conjecture over these particular real function fields \cite[Conjecture~4]{BenoistPI}. The final assertion is the curve case over \(\R\) of Lang's 1953 conjecture \cite[p.~379]{Lang}\cite[Conjecture~3]{BenoistPI}. See \cite[Section~3.6]{WittenbergICM} for this implication from tight approximation. Earlier results include the odd-degree case proved by Lang, Witt's quadratic case, and Benoist's result on degree-four del Pezzo surfaces over the function field of any real curve without real points \cite[\S0.5 and Theorem~0.13]{BenoistSurfaces}. For compactifications of homogeneous spaces of connected linear algebraic groups, existence was already known by Steinberg's theorem, as explained in \cite[Conjecture~4 and following discussion]{BenoistPI}, and the all-places weak approximation assertion follows from \cite[Theorems~7.4 and 7.7]{BW}.

\begin{corollary}[A level-two case of Lang's conjecture]
\label{cor:lang-c2-level-two}
Let \(S\) be a smooth projective geometrically integral real
surface. If \(-1\) is a sum of two squares in \(\R(S)\), then
\(\R(S)\) is \(C_2\). In particular this holds when
\(S(\R)=\varnothing\) and \(H^2(S,\cO_S)=0\).
\end{corollary}

\begin{proof}
Let \(Q=\{x_0^2+x_1^2+x_2^2=0\}\subseteq\mathbf P^2_\R\).
The level hypothesis gives an \(\R(S)\)-point of \(Q\), hence a
rational map \(S\dashrightarrow Q\). Its image is dense, since
\(\R\) is algebraically closed in \(\R(S)\) and \(Q\) has no real
point. By Corollary~\ref{cor:points}, \(\R(Q)\) is \(C_1\). The
Lang--Nagata transition theorem
\cite{LangTransition,Nagata} applies to the extension
\(\R(Q)\subseteq\R(S)\) of transcendence degree one. Benoist and
Wittenberg prove the required level bound under the final
hypotheses \cite[Theorem~1.8]{BWWu}.
\end{proof}

The analogous reduction to the anisotropic real conic in transcendence
degree one appears in \cite[Lemma~22]{BevelacquaMotley}. Here the new
input is Corollary~\ref{cor:points}, applied over the function field of
that conic before passage to transcendence degree two.

\begin{corollary}\label{cor:rec}
For smooth proper geometrically rationally connected varieties over \(F=\R(B)\), the reciprocity obstruction is the only obstruction to weak approximation in the sense of \cite[Definition~2.7]{BW}. For smooth projective such varieties, the Colliot-Th\'el\`ene obstruction is the only obstruction to the local--global principle for rational points, as formulated in \cite[Conjecture~2.6]{PS}.
\end{corollary}

\begin{proof}
For the first assertion, let \(F_b\) be the completion at the closed point \(b\) of \(B\). On a proper regular model \(p\), consider the families in \(\prod_b X(F_b)\) whose real components are induced by a single smooth section of \(p(\R)\) on \(B(\R)\), in the sense of formal jets. Theorem~\ref{thm:main} and \cite[Proposition~2.9 and Corollary~2.10]{BW} give density of the diagonal image in this locus, and consequently in its closure. This concerns actual local points at the selected places, including places of bad reduction.

For the second assertion we may assume \(B(\R)\ne\varnothing\) by Corollary~\ref{cor:points}. The nonemptiness of the cohomologically unobstructed locus gives, by Ducros's theorem as stated in \cite[Theorem~1.3]{PS}, a continuous semialgebraic section on \(B(\R)\) with image in the relative smooth locus. Approximation of sections of a submersion makes it smooth, as in \cite[Section~2.3]{BW}. Then Theorem~\ref{thm:main} and \cite[Corollary~2.4]{BW} give an algebraic section. Its generic value is an \(F\)-point.
\end{proof}

The local--global assertion answers the conjecture formulated by P\'al and Szab\'o in 2020 \cite[Conjecture~2.6]{PS}. It uses Ducros's 1998 topological description of the Colliot-Th\'el\`ene obstruction \cite[Th\'eor\`eme~4.3]{Ducros}. Earlier cases include smooth compactifications of homogeneous spaces of connected linear algebraic groups, fibrations into Severi--Brauer varieties over the projective line and the fibrations with rational base in the fibration theorem of P\'al and Szab\'o \cite[Theorems~1.2, 1.4, 1.5 and 5.3]{PS}. The density assertion answers the compatible real interpolation question, for data coming from local sections, formulated by Ducros, as recalled in \cite[Question~3.3(2) and following discussion]{WittenbergICM}. The homogeneous-space result of Benoist and Wittenberg already gives weak approximation for arbitrary stabilizers; their introduction records Scheiderer's earlier theorem for connected stabilizers \cite[Theorem~7.7]{BW}.

\subsection{Real rational curves and smooth maps}\label{sec:real-curves}

\begin{corollary}[Koll\'ar's question]\label{cor:conic}
Every positive-dimensional smooth proper geometrically rationally connected variety \(Y\) over \(\R\) with \(Y(\R)=\varnothing\) contains a geometrically integral curve of geometric genus zero.
\end{corollary}

\begin{proof}
Take the real conic \(Q=\{x_0^2+x_1^2+x_2^2=0\}\) and apply Theorem~\ref{thm:main} to \(Y\times Q\to Q\), with empty domain and interpolation set. The resulting real morphism \(Q\to Y\) is nonconstant. Otherwise it would factor through a closed point whose residue field would embed in \(H^0(Q,\cO_Q)=\R\), giving a real point of \(Y\). Its image is a geometrically integral curve. Over \(\C\), the induced nonconstant map from \(Q_\C\simeq\Pone_\C\) to the normalization of this image forces its genus to be zero by Riemann--Hurwitz.
\end{proof}

This answers Koll\'ar's question \cite[Remark~20(1)]{AK}, restated in
\cite[Conjecture~5]{BenoistPI}, by the constant-family specialization
anticipated in \cite[Remark~3.7(i)]{BW}. Earlier anisotropic-conic
constructions include the examples of Allums near ordinary-double-point
degenerations and in genus twelve \cite[Corollary~2.7.5 and \S2.8]{Allums}.
Benoist and Wittenberg obtained even geometric genus for all smooth real
Fano threefolds \cite[Corollary~7.2]{BWhodgeII}.

\begin{corollary}\label{cor:realmaps}
Let \(Y\) be smooth, proper and geometrically rationally connected over \(\R\), and let \(B\) be a smooth projective connected real curve. Every smooth map \(B(\R)\to Y(\R)\) is a compact-open \(C^\infty\) limit of real algebraic morphisms \(B\to Y\) on real loci. Moreover \(H_1(Y(\R),\mathbb Z/2)\) is generated by Borel--Haefliger classes of real rational curves.
\end{corollary}

\begin{proof}
Apply Theorem~\ref{thm:main} and \cite[Corollary~2.4]{BW} to the constant family \(Y\times B\to B\). The corollary uses real-analytic approximation and extension of real-analytic maps to holomorphic germs. For the second assertion take \(B=\Pone_\R\). The classes in \(H_1(Y(\R),\mathbb Z/2)\) are generated by smooth loops, and sufficiently close maps of the circle are homotopic. The pushforward of the fundamental class of \(\Pone(\R)\) by a real algebraic map is the Borel--Haefliger class of its pushforward one-cycle. The nonconstant images are geometrically rational curves. The loop approximation and the generation assertion also answer \cite[Questions~3.1(1) and (2)]{WittenbergICM} in the parametrized formulation described there.
\end{proof}

This answers \cite[Questions~1.1 and 8.13]{BW}. The deduction of
rational-curve generation from tight approximation is in
\cite[\S8.4]{BW}. Earlier cases include real rational targets
\cite[Section~2]{BK1999} and the varieties in \cite[Theorem~A]{BW}.
For Question~8.13, see \cite[Remark~8.12 and Proposition~8.14]{BW}
for cubics and geometrically rational surfaces.

\Needspace{11\baselineskip}
\section{Real integral Hodge classes and unramified cohomology}\label{sec:real-hodge}

\subsection{Complex degree-two comparison used in the descent}

We include the sphere argument of Hu and Wu
\cite[Theorem~1.1 and \S2.2]{HWshort} in the form needed for descent,
using the analytic construction developed in this paper.

\begin{lemma}[Hu--Wu]\label{lem:complexcycles}
If \(Y\) is smooth projective and rationally connected over \(\C\), then \(H_2(Y(\C),\mathbb Z)\) is generated, integrally, by the pushforwards of the fundamental class of \(\Pone(\C)\) under morphisms \(\Pone_\C\to Y\).
\end{lemma}

\begin{proof}
The dimension-zero case is immediate. For a class \(\beta\in H_2(Y(\C),\mathbb Z)\), choose a sign \(\varepsilon\) for which \(e=\langle c_1(TY),\varepsilon\beta\rangle\geq 0\). The integral Hurewicz theorem and simple connectivity give a smooth sphere \(s\colon\Pone(\C)\to Y(\C)\) with \(s_*[\Pone(\C)]=\varepsilon\beta\). By the classification of complex bundles on a sphere using determinant clutching, \(s^*TY\) is smoothly isomorphic to \(V=\cO(e)\oplus\cO^{\dim Y-1}\), for which \(H^1(V)=0\).

Apply Proposition~\ref{prop:analytic} over \(\C\) to the product family
on the sphere, with empty protected sets and first-jet families from
\(Y\). For the complex version, take the disjoint union of the
source with its conjugate and of the target with its conjugate, with
the involutions exchanging the two copies. Use \(s\) and \(V\) on the
first copy and their conjugates on the second. Let \(F\) be the
restriction of the equivariant correction to the first copy and \(g_i\) its vertical
rational spheres. The homotopy from the correction to the zero parameter
and the doubled-disk computation in the proof of Theorem~\ref{thm:effective}
give, before pairing with any divisor,
\[
 F_*[\Pone(\C)]=\varepsilon\beta+\sum_i(g_i)_*[\Pone(\C)]
 \quad\text{in }H_2(Y(\C),\mathbb Z).
\]
Indeed the old disk, annulus and cap lie in a contractible target chart,
and the change \(w=\rho_i^2/\zeta\) preserves complex orientation. Thus
each disk contributes the entire sphere with its integral multiplicity.
GAGA makes \(F\) and the \(g_i\) algebraic, and the identity applies
also to torsion classes.
\end{proof}

We first prove the complex input within this paper. We use singular
Borel--Moore homology with integer coefficients. On a smooth complex
variety \(Y\) of positive dimension, the integral cycle classes in
\(H_2(Y(\C),\mathbb Z)\) for proper \(Y\) form the kernel of
\(H_2(Y(\C),\mathbb Z)\to\operatorname{KH}_2(Y,\mathbb Z)\).
The same statement with Borel--Moore homology holds for a reduced
quasi-projective complex variety \(Y\), possibly singular. Indeed, in the
homological niveau spectral sequence the terms at \((p,q)\) vanish for
\(q<0\) and for \(q>p\). To see this, take the direct limit over smooth
affine dense opens of each \(p\)-dimensional subvariety. Duality identifies
the term there with degree-\((p-q)\) cohomology, and a smooth affine
complex variety of dimension \(p\) has the homotopy type of a CW complex
of dimension at most \(p\). Thus \(E^2_{2,0}=\operatorname{KH}_2(Y,\mathbb Z)\)
survives, and the remaining piece of the filtration on \(H_2^{\rm BM}\)
is generated by the fundamental classes of integral curves by localization.
This gives the natural exact sequence
\begin{equation}\label{eq:local-complex-kato-two}
 \operatorname{CH}_1(Y)\longrightarrow H_2^{\rm BM}(Y(\C),\mathbb Z)
 \longrightarrow\operatorname{KH}_2(Y,\mathbb Z)\longrightarrow0.
\end{equation}

\begin{lemma}\label{lem:local-complex-kato-two}
Let \(g:P\to T\) be a projective surjective morphism of smooth
connected quasi-projective complex varieties with geometrically
rationally connected generic fiber. Then \(g_*\) induces an
isomorphism on the quotients of \(H_2^{\rm BM}(-,\mathbb Z)\) by
integral algebraic one-cycle classes. If \(T\) is a curve and all
reduced fibers above a finite set \(D\) form a simple normal crossings
divisor in \(P\), then each reduced fiber over \(D\) has its integral
second homology generated by algebraic curves.
\end{lemma}

\begin{proof}
Work with the properties \(R(n,2),F(n,2),D(n,2),S(n,2)\) of
\cite[Theorem~6.1, preprint version]{TianLaurent}, for this singular
homology theory. A smooth projective rationally connected complex
variety is simply connected. Together with Lemma~\ref{lem:complexcycles}
and \eqref{eq:local-complex-kato-two}, this gives \(RC(n,2)\) in all
dimensions. The statement at dimension zero is immediate from the
complex for a point. Theorem~6.1(3) therefore gives \(S(n,2)\) for
all \(n\). Both \(R(1,2)\) and \(F(1,2)\) hold: klt curves are smooth,
and a fibration of a smooth curve is either over a point with rational
total curve or is birational and hence an isomorphism. The corresponding
Kato groups follow from the same niveau sequence. From \(R(n,2)\)
and \(F(n,2)\), parts (1) and (2) give \(R(n+1,2)\) and \(D(n+1,2)\).
Then part (4) gives \(F(n+1,2)\). This proves the assertion on
quotients by \eqref{eq:local-complex-kato-two}.

If \(T\) is a curve, apply the degeneration property \(D(n,2)\) to the
divisor of reduced fibers above \(D\); its second Kato homology
identifies with that of the finite set \(D\) and hence is zero. If
\(\dim P=1\), the reduced fibers are already points. The fibers are
disjoint, and their ordinary and Borel--Moore homology agree by
properness. The last assertion now follows from
\eqref{eq:local-complex-kato-two}.
\end{proof}

\subsection{Transfer along rationally connected fibrations}\label{sec:real-hodge-transfer}

\begin{theorem}\label{thm:real-hodge-transfer}
Let $f\colon X\to B$ be a surjective morphism of smooth projective geometrically
integral real varieties with geometrically rationally connected generic
fiber. If $B$ satisfies the real integral Hodge conjecture for one-cycles,
then so does $X$.
\end{theorem}

Benoist and Wittenberg proved the conic-fiber case
\cite[Theorem~6.1]{BWhodgeII}. The real-locus step below adapts
their approximation argument from \cite[proof of Theorem~8.5]{BW}.
For arbitrary relative dimension, Theorem~\ref{thm:main} supplies
sections after restriction to real curves, including curves without
real points. The free-action step uses the ordinary transfer and the
signed norm on the two quotient spaces, together with the degree-two
comparison of Lemma~\ref{lem:local-complex-kato-two}. The absolute
input to that comparison is Lemma~\ref{lem:complexcycles},
and its passage to fibrations uses Tian's induction
\cite[Theorem~6.1, preprint version]{TianLaurent}.

The case \(\dim X=0\) is immediate. Put $G=\operatorname{Gal}(\C/\R)$, $F_Y=Y(\R)$,
$Q_Y=Y(\C)/G$ and $Y_G=EG\times_G Y(\C)$. Write $\mathbb Z_-$ for the sign
local system from $BG$. For $Y$ smooth projective of dimension $d$, the
topologically admissible group
\[
 A_Y=H_G^{2d-2}(Y(\C),\mathbb Z(d-1))_0
\]
identifies with $H_2(Y_G,F_Y\times BG;\mathbb Z_-)$.
Indeed, \cite[Lemma~3.5]{BWhodgeI} identifies $A_Y$ with the corresponding
cohomology of $Y(\C)\setminus F_Y$. Poincar\'e duality on its free
quotient gives Borel--Moore homology with sign coefficients, since the
orientation character is $(-1)^d$. Relative Borel chains compute this
same homology. On an equivariant triangulation with fixed subcomplex
$F_Y$, every relative cell has trivial stabilizer. In
particular, with $\mathbb F_2$ coefficients the relative Borel group equals the
ordinary relative homology of $(Q_Y,F_Y)$.

\begin{lemma}\label{lem:transfer-conic}
Let $C$ be a smooth real conic without real points. Every real integral
Hodge class in $A_Y$ with zero Borel--Haefliger class is the pushforward
of a real integral Hodge class on $Y\times C$.
\end{lemma}

\begin{proof}
The pair connecting map
\[
 H_2(Y_G,F_Y\times BG;\mathbb Z_-)
   \longrightarrow H_1(F_Y\times BG;\mathbb Z_-)
   \xrightarrow{\simeq} H_1(F_Y;\mathbb F_2)
\]
is the Borel--Haefliger map on $A_Y$. For the second arrow, reduction
modulo two followed by projection induces an isomorphism, since
$H_0(BG;\mathbb Z_-)=\mathbb Z/2$, $H_1(BG;\mathbb Z_-)=0$ and $H_0(F_Y;\mathbb Z)$ is free.
Naturality identifies the displayed composite with the mod-two boundary
for $(Q_Y,F_Y)$ after reduction. To compare with Borel--Haefliger, this
boundary is dual to the connecting map on $H^1(F_Y;\mathbb F_2)$ in the
cohomology exact sequence of that pair. In the relative Borel
description it is $u_2$ of \cite[(1.46)]{BWhodgeI} on the $H^1(F_Y;\mathbb F_2)$
summand. The factor $\gamma$ has degree-zero term $1$, and higher terms
are discarded in that degree. Duality in \cite[Theorem~1.12 and
Section~3.2]{BWhodgeI} therefore identifies the boundary with the degree
$d-1$ component of restriction to the real locus, expressed in homology.

The vanishing of this boundary lifts the given class to
$H_2(Y_G;\mathbb Z_-)$. The Borel construction of $Y\times C$ over $Y_G$ is a
bundle with fiber $C(\C)\simeq S^2$. Its pushforward on $H_2$ with
coefficients pulled back from $\mathbb Z_-$ is onto by the homological Serre
spectral sequence. The row $q=0$ in total degree two has no incoming
differential and its only possible outgoing differential has target
in the zero row $q=1$. As $Y\times C$ has empty real locus, the lift is
admissible. To compare pushforward with the equivariant cohomological
pushforward, restrict to the inverse image of $Y(\C)\setminus F_Y$ and
use proper pushforward in Borel--Moore homology there. This is exactly
the relative homology map of pairs, and Poincar\'e duality is compatible
with this operation. Hence the lifted class pushes to the given class
in $A_Y$.

Finally, the ordinary complex degree-two homology of the product is the
direct sum of the ordinary degree-two homology of its factors. The
component of the lift in $Y(\C)$ is the original complex Hodge class, and
the component in $C(\C)$ is algebraic. Thus the lift is Hodge.
\end{proof}

\begin{lemma}\label{lem:transfer-bh}
In the situation of the theorem, if $\beta\in H_1(X(\R);\mathbb F_2)$ has
algebraic pushforward in $H_1(B(\R);\mathbb F_2)$, then $\beta$ is the
Borel--Haefliger class of a real algebraic one-cycle.
\end{lemma}

\begin{proof}
If $\beta=0$, take the zero cycle. Otherwise we use the method of
\cite[proof of Theorem~8.5]{BW} on a single homology class.
First, for any projective birational morphism $W\to V$
of smooth projective real varieties, a class in $H_1(W(\R);\mathbb F_2)$ with
algebraic image in $H_1(V(\R);\mathbb F_2)$ is algebraic. For a blowup of a smooth
center this follows from the Mayer--Vietoris description of the blowup.
The kernel is generated by real projective lines in the exceptional
fibers. Indeed, the real projective normal bundle has the same $H_0$ as
the center and surjects onto its $H_1$, with kernel generated by the
fiber classes. The Chow moving lemma followed by strict transform lifts
an algebraic class on $V$ to an algebraic class on $W$. For a general
birational morphism, resolve its inverse by a sequence of blowups of
$V$ with smooth centers,
giving a smooth projective variety dominating $W$. Its map to $W$
surjects on $H_1$ of real loci, since loops can be moved off the
codimension-at-least-two subset where the map is not an isomorphism.
Apply the blowup argument on the dominating variety and push the
algebraic cycle to $W$.

Take a dense open $B^0\subset B$ over which $f$ is smooth with
geometrically rationally connected fibers. Apply weak toroidalization
\cite[Theorem~1.1]{ADK} to $f$, with $f^{-1}(B\setminus B^0)$ as its
specified closed subset. This gives a toroidal morphism between smooth
projective modifications with strict normal crossings boundaries,
as in \cite[Section~5]{BW}. The morphism is smooth on the source
toroidal open. On a dense open of the modified base, the fibers remain
smooth and geometrically rationally connected by birational invariance
and openness. Lift $\beta$ to the modified source by the
surjectivity on $H_1$ just described. Its projection to the modified
base remains algebraic by the preceding paragraph. Choose a real
analytic representative $\iota\colon L\to X(\R)$ of the lifted class, where
$L$ is a disjoint union of circles. Perturb first so each component meets
the source toroidal open above the good base open, then use real analytic
approximation. The inverse image of the toroidal boundary is finite.
We use $X,B$ for the modified varieties.

Choose an analytic embedding $L\hookrightarrow\mathbf A^2(\R)$ in an
affine chart of $\mathbf P^2(\R)$. Replacing $f$ by $f\times
\mathrm{id}_{\mathbf P^2}$ and taking the graph of this embedding
makes both the loop and its projection embeddings. The graph adds
zero in the $\mathbf P^2$ summand of $H_1$ and thus preserves
algebraicity of the projected class by K\"unneth. At every boundary
intersection the projected analytic germ pulls back some algebraic
local coordinate to a uniformizer, since its differential is
injective. Each discrete valuation ring of real analytic germs has
residue field $\R$, and its generic point lies in the good open.
Thus \cite[Proposition~5.2]{BW} applies and gives modifications on which the
morphism is submersive along the whole lifted loop. The projected
strict transform is an embedding because its composite with the
modification is an embedding.
The projected homology class is still algebraic by the birational
argument above.

The submersion has a smooth section on a neighborhood of this
embedded projected loop, taking the prescribed values along the loop,
by \cite[Lemma~6.11]{BWhodgeII}. By the Akbulut--King theorem
\cite[Theorem~6.8]{BWhodgeII}, choose a real algebraic curve, nonsingular
along its real locus, whose real locus approximates that loop in the
neighborhood and represents its mod-two class there by an isotopy.
Each circle of the projected loop contains a point in the good open.
We can keep the approximation sufficiently close at those points so
that every component with real points meets the good open. Components
without real points have zero Borel--Haefliger class and may be
discarded. On the normalizations of the remaining components, the
smooth section on the neighborhood pulls back to a smooth real section
on the entire real locus. It lifts to a smooth projective model of
the main pullback, chosen to be an isomorphism on the relative smooth
locus through which the section passes. The geometric generic fiber of
each pullback is rationally connected. By \cite[Proposition~2.3]{BW},
first approximate the smooth section by an equivariant holomorphic
section near the real locus. Theorem~\ref{thm:main} then gives algebraic
sections arbitrarily close to the smooth section on the real locus.
Sufficiently close sections are homotopic there and induce the same
mod-two classes. Pushing their sum through the modifications gives
the desired one-cycle on the original $X$.
\end{proof}
\begin{lemma}\label{lem:transfer-free}
The theorem holds if $B(\R)=\varnothing$.
\end{lemma}

\begin{proof}
Let $\alpha$ be a real integral Hodge class on $X$. Its pushforward is
algebraic on $B$. Use the Chow moving lemma to choose a representative
one-cycle with no component contained in the complement of the good
open. The normalization of each component has no real points, since
it maps to $B$. On a smooth projective model of its pullback,
Theorem~\ref{thm:main} with empty approximation domain supplies a section.
Subtracting the corresponding algebraic one-cycle on $X$
reduces to $f_*\alpha=0$.

For $Y=X,B$ the ordinary double cover $Y(\C)\to Q_Y$ gives a chain
sequence with ordinary transfer $T_Y$ and signed pushforward $N_Y$,
\[
 H_2(Q_Y;\mathbb Z)\xrightarrow{T_Y}H_2(Y(\C);\mathbb Z)
  \xrightarrow{N_Y}H_2(Q_Y;\mathbb Z_-)
  \longrightarrow H_1(Q_Y;\mathbb Z).
\]
The map $N_Y$ is the cycle class norm under duality.

The map $f_*:H_1(Q_X;\mathbb Z)\to H_1(Q_B;\mathbb Z)$ is an isomorphism, and
$f_*:H_2(Q_X;\mathbb Z)\to H_2(Q_B;\mathbb Z)$ is surjective. To prove this, remove
the complex codimension at least two strata in the base discriminant and
shrink its divisorial strata. At every remaining point of a discriminant
divisor there is a local complex section. Indeed, stratify its inverse
image so that each remaining stratum submerses onto the divisor. A
transverse complete intersection curve through the point has a section
on the resolved dominant pullback by Graber--Harris--Starr \cite{GHS}.
At its section value the derivative gives a transverse direction, and
the stratum gives the tangent directions along the divisor. Thus $f$
is submersive there. The same holds in the free quotients.
Represent an integral degree-two class by a map from a closed
oriented surface, allowing a disjoint union. By smooth general
position the map avoids strata of real codimension at least four
and meets divisorial strata in finitely many transverse points.
Lift on disks around those points using the local sections, then
extend across the remaining two-dimensional domain in the proper
smooth bundle. Its fibers are connected and simply connected,
so a lift along the boundary extends over the successive one- and
two-cells. This proves surjectivity in degree two. Loops on the
base can be moved into the good open and lifted there, using
connectedness of the fiber to close the lift. A loop on
the source can be moved into the inverse image of the good open.
If its projection is null-homotopic, use the same lifting argument
on a bounding disk, prescribing the original lift on the boundary.
This gives an isomorphism on fundamental groups and hence on their
abelianizations.

Now the connecting image of $\alpha$ in $H_1(Q_X;\mathbb Z)$ is zero, so
$\alpha=N_X\beta$ for some $\beta\in H_2(X(\C);\mathbb Z)$. Its pushforward
under $f$ lies in the image of $T_B$. Lift a preimage under $T_B$ to
$H_2(Q_X;\mathbb Z)$ and subtract its transfer from $\beta$. This keeps the
norm and makes $f_*\beta=0$. The relative degree-two comparison
of Lemma~\ref{lem:local-complex-kato-two} makes this corrected class algebraic. Its norm
represents $\alpha$.
\end{proof}

\begin{proof}[Proof of Theorem~\ref{thm:real-hodge-transfer}]
The pushforward of the Borel--Haefliger
class of a real integral Hodge class $\alpha$ on $X$ is algebraic,
by the hypothesis on the base and functoriality of the real cycle class
condition \cite[Theorem~1.21]{BWhodgeI}. Lemma~\ref{lem:transfer-bh} supplies an
algebraic correction, leaving $\operatorname{BH}(\alpha)=0$. Take a real conic $C$
without real points. The real integral Hodge conjecture for $B\times C$
follows from \cite[Theorem~6.1]{BWhodgeII}. Lemma~\ref{lem:transfer-free} applied to
$X\times C\to B\times C$ gives the real integral Hodge conjecture on
$X\times C$. By Lemma~\ref{lem:transfer-conic}, $\alpha$ is the pushforward of a
Hodge class on this product and is therefore algebraic.
\end{proof}

\begin{corollary}\label{cor:real-mrc-hodge}
Let $X$ be a smooth projective geometrically integral real variety whose
maximal rationally connected quotient admits a smooth projective real
model of dimension at most two. Then $X$ satisfies the real integral
Hodge conjecture for one-cycles.
\end{corollary}

\begin{proof}
Resolve the rational quotient map on a smooth projective birational
model of $X$. The base satisfies the conjecture vacuously in
dimension zero and by its fundamental class in dimension one.
In dimension two it satisfies Krasnov's real Lefschetz theorem
\cite[Proposition~2.8]{BWhodgeI}. Apply
Theorem~\ref{thm:real-hodge-transfer} to the resolved morphism.
Birational invariance gives the result on $X$
\cite[Proposition~2.13]{BWhodgeI}.
\end{proof}

\subsection{Geometrically uniruled threefolds}

\begin{theorem}\label{thm:uniruled-real-hodge}
Let \(X\) be a smooth proper geometrically integral threefold over \(\R\). If \(X_\C\) is uniruled, then \(X\) satisfies the real integral Hodge conjecture for one-cycles.
\end{theorem}

For conic bundles over surfaces this assertion is \cite[Corollary~6.2]{BWhodgeII}. Smooth Fano threefolds without real points were settled in \cite[Theorem~7.1]{BWhodgeII}.

We give a separate relative Hilbert scheme argument for fibrations
into geometrically rational surfaces. Benoist and Wittenberg treated del Pezzo fibrations of degree at least five, those with empty real total space over a curve with real points, those of degree one with connected smooth real fibers, and all cubic-surface fibrations over \(\R\) \cite[Theorem~8.1 and Example~8.2]{BWhodgeII}. The degree-four case over a curve without real points follows from Benoist's real period-index theorem and the vertical cycle and Leray arguments of Benoist--Wittenberg \cite[Theorem~0.13]{BenoistSurfaces}\cite[Propositions~8.6 and 8.11]{BWhodgeII}\cite[Theorem~1.2.14]{BenoistHDR}. Their Proposition~8.3 gives the equivalence with section existence over the anisotropic conic. Here the relative Hilbert scheme construction produces horizontal cycles for arbitrary geometrically rational surface fibrations, including when the real fibers have several connected components or contain components of multiplicity greater than one. It retains the vertical cycle construction and the equivariant Leray argument of \cite[\S8]{BWhodgeII}.

\paragraph{Paths of real configurations.}\label{par:real-hilbert-paths}
In a smooth projective geometrically rationally connected real variety
of positive dimension \(n\), reduced real configurations of length
\(N\) with the same parity on each real component are joined in the
smooth locus of the main Hilbert component. For fixed real-point
counts on each component, the unordered real configurations can move
within those components. The nonreal pairs can move in the connected
quotient of the complement of the real locus by conjugation. For
\(n\geq2\) that complement is connected by real codimension at least
two; for \(n=1\) a split real conic has its two complementary
hemispheres exchanged by conjugation. For a point-free conic the
complement is also connected. The unordered configurations of
distinct nonreal conjugate pairs are therefore path connected.

To change the counts, convert surplus real pairs to conjugate pairs
before using them at components with deficits. The local family
\[
 z_1^2-v=0,\qquad z_2=\cdots=z_n=0
\]
runs through a curvilinear double point. The chart in which the
quotient has basis \(1,z_1\) is given by
\(z_1^2-az_1-b\) and \(z_j-c_jz_1-d_j\) for \(2\leq j\leq n\).
The \(2n\) coefficients are free, and the other points vary
independently. These coordinates also give relative smoothness over a
smooth base family. The chart meets the reduced locus, hence belongs
to the main Hilbert component.

\begin{proposition}\label{prop:real-hodge-fibration}
Let \(p\colon X\to B\) be a morphism from a smooth projective geometrically integral threefold over \(\R\) to a smooth projective geometrically integral curve. Suppose that the geometric generic fiber is a smooth rational surface. Then
\[
  \operatorname{cl}\colon \operatorname{CH}_1(X)\longrightarrow
  H^4_G(X(\C),\mathbf Z(2))_0
\]
is surjective, where \(G=\operatorname{Gal}(\C/\R)\) and the subscript denotes the topological condition of \cite[Definition~1.19]{BWhodgeI}.
\end{proposition}

\begin{proof}
The morphism is flat since \(X\) is integral and \(B\) is regular of dimension one. The coherent cohomology satisfies \(H^2(X,\cO_X)=0\) by \cite[Lemma~6.6]{BWhodgeII}, so the displayed group is the target for the real integral Hodge conjecture. By its birational invariance \cite[Proposition~2.13]{BWhodgeI}, make an embedded resolution of the supports of the finitely many fibers outside an open over which \(p\) is smooth with geometrically rational fibers. We continue to write \(X\) for the resulting smooth projective threefold and \(p\) for its morphism to \(B\). All fibers outside this open now have reduced supports with simple normal crossings.

\smallskip
\emph{Fiber restrictions.}
Fix \(\alpha\in H^4_G(X(\C),\mathbf Z(2))_0\). Its restriction to each good complex fiber has an integer degree \(d\), independent of the fiber by local constancy over the connected complex open curve. On a good real fiber \(X_b\), its restriction belongs to \(H^4_G(X_b(\C),\mathbf Z(2))_0\). Such a class is determined by its degree and by its Borel--Haefliger components in \(H^2(X_b(\R),\mathbf Z/2)\). Indeed, the kernel of the Borel--Haefliger component map is the image of the norm, whose composition with degree is multiplication by two in this case \cite[Lemma~2.11]{BWhodgeI}. In the component coordinates of \(H^2(X_b(\R),\mathbf Z/2)\), the sum of the entries is \(d\bmod 2\). This follows by writing the class as a sum of classes of real points and norms of complex points, as in \cite[Proposition~2.10]{BWhodgeI}.

Let \(\beta\in H_1(X(\R),\mathbf Z/2)\) be the Poincar\'e dual of the degree-two component of the restriction of \(\alpha\) to the real locus, in the notation of \cite[(1.57)]{BWhodgeI}. Represent \(\beta\) by a compact smooth embedded one-dimensional submanifold \(L\subset X(\R)\), using the usual smoothing of mod-two one-dimensional cycles and general position in a three-manifold. We take \(L=\varnothing\) if \(\beta=0\). At good real values transverse to \(L\), the parity entries on \(X_b(\R)\) are the intersection counts of \(L\) with its connected components, by naturality of the component map and mod-two Poincar\'e duality. These entries vary locally constantly along the smooth real fibration.

If \(B(\R)=\varnothing\), Corollary~\ref{cor:points} gives a generic
point, hence a section by properness. Subtract \(d\) times its class
and proceed to localization. In the rest of the horizontal construction
assume \(B(\R)\ne\varnothing\).

\smallskip
\emph{Subschemes over bad real values.}
Perturb \(L\) to avoid the intersections of the components of each bad fiber and to meet the smooth locus of its reduced support transversely. The excluded strata on the real locus have dimension at most one. At each intersection point \(x\) over a real value \(b\), the support has a single local component, of multiplicity \(m\geq1\). In coordinates compatible with conjugation at \(x\) and a real local parameter \(t\) at \(b\), write \(p\) as
\[
  t=\varepsilon_x w^m,\qquad \varepsilon_x\in\{1,-1\},
\]
with the other two coordinates denoted \(z_2,z_3\). This follows by absorbing an \(m\)-th root of the unit after extracting its real sign. The germ \(z_2=z_3=0\) is finite flat of degree \(m\) over a small disk at \(b\). It defines a local holomorphic section of the relative Hilbert scheme of length-\(m\) subschemes. For nonzero real \(t\) on either side, its real points are in the same respective connected components of \(X_t(\R)\) as those of \(L\) near \(x\). Indeed \(w\) is a parameter on \(L\) there by transversality, and at each real root one joins the two points by varying \(z_2,z_3\) with \(w\) fixed. Thus the sum \(D_b\) of the multiplicities \(m\) over the intersections of \(L\) with \(X_b\) is congruent to \(d\) modulo two. For even \(m\) there are either two or zero real roots, and for odd \(m\) there is one on either side.

Fix an integer \(N\equiv d\pmod 2\) at least as large as all \(D_b\) and as the number of nonzero parity entries on each good real interval and circle. We may take \(N\geq2\). At each bad value complete the disjoint union of the chosen germs to degree \(N\) by adjoining \((N-D_b)/2\) distinct conjugate pairs of local sections through nonreal relative-smooth points. There are such points even over a nonreduced fiber. The Graber--Harris--Starr theorem supplies a section after complexification, and its derivative shows that it meets the fiber in the relative smooth locus. This locus contains a nonempty complex open part of a surface, from which the nonreal points can be chosen distinct and away from all previous supports.

The relative Hilbert scheme is smooth over \(B\) at the resulting subschemes over \(b\). At the length-\(m\) subscheme at \(x\) this can be checked on local Artinian real algebras \(A\) mapping to \(B\) at \(b\). Write \(t_A\) for the image of \(t\). The finite flat quotient algebra of any such deformation receives a map from \(A[w]/(w^m-\varepsilon_x t_A)\) that is an isomorphism on the special fiber. Both sides are free of rank \(m\) over \(A\), so this map is an isomorphism. The two transverse coordinates can then independently take arbitrary values in this algebra reducing to zero on the central subscheme. This gives \(2m\) free relative deformation coordinates in the completed local ring. The real analytic coordinate normal form also yields the required formal normal form and the finite flat holomorphic germs above. At a disjoint union the local deformation problems are independent, since idempotents of a finite algebra lift over nilpotent ideals. This also handles the added simple nonreal points.

\Needspace{4\baselineskip}
\smallskip
\emph{Connecting the subschemes over the real base.}
On a closed subinterval of the good locus meeting both prescribed
collars, trivialize the underlying complex surface family by
conjugation-equivariant smooth transport. Such transport is obtained
by averaging a horizontal lift of the base vector field with its
conjugate. In this trivialization, reparametrize each reduced end
germ so that it becomes constant after a smaller collar, and extend
it toward an interior value \(c\). The transport takes real
components to real components and preserves distinct support points.
The two configurations obtained at \(c\) therefore have the same
component parities, since those parities are the restrictions of
\(\alpha\) in the smooth family.

Join these configurations by the real configuration path above in
\(\operatorname{Hilb}^N(X_c)(\R)\), smoothed with its endpoints fixed.
The relative Hilbert scheme is
smooth over the good locus by \cite[Theorem~2.9]{Fogarty}. In a
smooth trivialization of its real points near \(c\), the two
transported germs are close to the respective endpoints of the
path. Use short interpolations in coordinate balls at those
endpoints and traverse the path over a small interval around \(c\).
Smooth cutoffs in the endpoint charts and a path stationary near
its endpoints give a smooth section agreeing with the transported
germs outside that interval. A reduced configuration of any
required parity vector exists by the choice of \(N\), using one
point on each odd component and filling the remaining length by
distinct conjugate pairs. On a good real circle, choose a reduced
germ with the required parities at a cut point and keep the same
germ at both ends after cutting. The interval construction then
glues smoothly because the parity vector is invariant under the
component monodromy. This gives a smooth real section along all of
\(B(\R)\), with a reduced sample in each component of the good
real locus chosen outside the transition interval.

The generic Hilbert scheme of length-\(N\) points is smooth and geometrically integral by \cite{Fogarty}. It is birational to the symmetric power of the generic surface, which is geometrically rationally connected as it is dominated by the product of \(N\) geometrically rational surfaces. Take the closure of this generic Hilbert scheme with its reduced structure and a projective resolution that is an isomorphism on its regular locus. Every point of the constructed smooth real section lies in the relative smooth locus of that closure. In fact all the local charts used for the section are smooth over the base, so the unique component that dominates the base contains them locally. The resolved model is flat by integrality over the regular curve. By \cite[Proposition~2.3]{BW}, the smooth real section can be approximated by equivariant holomorphic sections near \(B(\R)\). Theorem~\ref{thm:main} then gives an algebraic real section sufficiently close to preserve the reduced configurations and their component parities at the sample values.

Pull back the universal length-\(N\) family along the resulting morphism to the Hilbert scheme. It is finite flat over \(B\) and its pushforward gives a one-cycle of degree \(N\) on \(X\), whose restriction at each sample is the prescribed parity vector. The norm of the class of a complex algebraic section of \(p_\C\) has degree two and zero real parity vector. Subtract \((N-d)/2\) times this norm from the cycle. The resulting class has exactly the fiber restrictions of \(\alpha\) at the samples by \cite[Lemma~2.11]{BWhodgeI}, and on all good fibers by local constancy.

\smallskip
\emph{Localization.}
After subtracting the horizontal cycle constructed above, the remaining class restricts to zero on each good fiber. Choose a good affine open \(B^0\subset B\) with no circle components in \(B^0(\R)\), by removing at least one point from each remaining real circle. The ordinary groups \(H^3(X_b(\C),\mathbf Z)\) vanish on this open since its fibers are geometrically rational surfaces. Proposition~8.11 of \cite{BWhodgeII} therefore shows that the class is zero on \(p^{-1}(B^0)\). The localization sequence expresses the remaining class as a sum of images of classes supported on omitted fibers. At real values these images are algebraic by \cite[Proposition~8.6]{BWhodgeII}. At a nonreal value the equivariant support group identifies with ordinary cohomology with support on one of the two conjugate fibers, with the global image given by the norm of the corresponding complex image. The proof of \cite[Proposition~8.6]{BWhodgeII} with the trivial group gives algebraicity of the ordinary image: the reduced normal crossings fiber has support classes generated by the component support groups as in Lemmas~8.7--8.8, and each smooth surface component has \(H^2(\cO)=0\), so ordinary Lefschetz \((1,1)\) applies after purity. Taking norms of those cycles gives the required real algebraic classes. This proves surjectivity.
\end{proof}

\begin{proof}[Proof of Theorem~\ref{thm:uniruled-real-hodge}]
Replace \(X\) by a smooth projective birational model by Chow's lemma,
resolution and \cite[Proposition~2.13]{BWhodgeI}. The maximal rationally
connected quotient of \(X_\C\) has dimension at most two, since
\(X_\C\) is uniruled. Its function field \(E\subset\C(X)\) is stable
under conjugation \(\sigma\) by the birational uniqueness of the
quotient \cite[\S1.1, before Corollary~1.4]{GHS}. The fixed field
\(E^\sigma\subset\R(X)\) gives a real quotient because
\(E^\sigma(\sqrt{-1})=E\). Take a smooth projective real model of
this quotient and apply Corollary~\ref{cor:real-mrc-hodge}.
\end{proof}

\begin{corollary}\label{cor:uniruled-real-unramified}
Let \(X\) be as in Theorem~\ref{thm:uniruled-real-hodge}. Then
\[
  H^3_{\mathrm{nr}}(X,\mathbf Q/\mathbf Z(2))_0=0,
\]
where the subscript denotes classes whose restriction to every real point vanishes. If \(X(\R)=\varnothing\), then \(H^3_{\mathrm{nr}}(X,\mathbf Q/\mathbf Z(2))=0\).
\end{corollary}

\begin{proof}
For uniruled threefolds, the normalized unramified group measures the real integral Hodge defect by \cite[Proposition~5.2 and Remark~5.3(iii)]{BWhodgeI}. Theorem~\ref{thm:uniruled-real-hodge} gives its vanishing. If the real locus is empty, the restriction condition is vacuous.
\end{proof}

\subsection{Geometrically rationally connected varieties}

The cycle-class assertion follows from Theorem~\ref{thm:real-hodge-transfer}
over a point. The stronger assertion about representatives uses the
Borel--Haefliger and anisotropic-conic constructions of
Section~\ref{sec:real-curves}.

\begin{corollary}\label{cor:realhodge}
Every smooth projective geometrically rationally connected variety over \(\R\) satisfies the real integral Hodge conjecture for one-cycles. Its real integral Hodge classes in this degree are represented by real algebraic one-cycles whose geometric components are rational.
\end{corollary}

\begin{proof}
The dimension-zero case is vacuous. In positive dimension rational
connectedness gives \(H^2(Y,\cO_Y)=0\) and, by simple connectivity
and the Kummer sequence, \(\operatorname{Pic}(Y_\C)[2]=0\).

For the assertion on representatives, use the norm exact sequence in \cite[Proposition~3.14]{BWhodgeI}. The quotient by the norm image is \(H_1(Y(\R),\mathbb Z/2)\) when the real locus is nonempty and \(\mathbb Z/2\) otherwise. In the first case its generators lift to real rational-curve classes. In the second case the class of the curve in Corollary~\ref{cor:conic} maps to the generator by \cite[Theorem~3.6]{BWhodgeI}. The kernel is generated by norms of complex rational-curve classes by Lemma~\ref{lem:complexcycles}. These norms are cycle classes of the pushforwards of the complex curves along \(Y_\C\to Y\). This supplies representatives with rational geometric components, allowing sums of conjugate curves and integral coefficients.
\end{proof}

For the smooth proper cycle-class assertion, pass to a projective birational
model and apply \cite[Proposition~2.13]{BWhodgeI}. Corollary~\ref{cor:realhodge}
answers the rationally connected case of \cite[Question~2.16]{BWhodgeI} in
all dimensions. Its complex input in dimension three was already known by
Voisin, as used in \cite[Corollary~3.23]{BWhodgeI}. Over non-archimedean real
closed fields there are counterexamples \cite[Example~9.7]{BWhodgeII}.

\subsection{Sum-of-squares bounds in dimensions five and seven}

\begin{corollary}\label{cor:real-fivefold-level}
Let \(X\) be a smooth projective geometrically integral real fivefold
with \(X(\R)=\varnothing\). Suppose \(X_\C\) has a smooth projective
birational model \(Y\) with a surjective morphism to a smooth projective
connected complex variety \(B\). Assume that its geometric generic
fiber is positive-dimensional and rationally connected, and that its
generic fiber has odd index over \(\C(B)\), so it admits a zero-cycle
of odd degree. Then \(-1\) is a sum of at most sixteen squares in
\(\R(X)\). In particular this holds if
\(X_\C\) is separably rationally connected in codimension one.
\end{corollary}

\begin{proof}
For a positive-dimensional base, \cite[Corollary~8.8, p.~38]{HW} gives an odd
annihilator for the kernel of
\(\operatorname{KH}_4(Y,\mathbb Z)\to\operatorname{KH}_4(B,\mathbb Z)\).
For a point base use the absolute fourth Kato vanishing in
\cite[Proposition~10.2, p.~41]{HW}. As \(\dim B\le4\), the target embeds in top field cohomology
when \(\dim B=4\) and is zero below that dimension. The top field group
is torsion-free \cite[Theorem~3.1]{CTVcomplex}. Every element of
\(\operatorname{KH}_4(Y,\mathbb Z)\) killed by two is therefore in
the relative kernel and is zero.

Let \(N^2\) denote coniveau at least two. The bottom niveau edge
identifies \(H^6(Y^{\rm an},\mathbb Z)/N^2\), by Poincar\'e duality,
with the quotient of \(H_4(Y^{\rm an},\mathbb Z)\) by homology
supported on subvarieties of dimension at most three. This quotient
is \(E^\infty_{4,0}\) and embeds into
\(E^2_{4,0}=\operatorname{KH}_4(Y,\mathbb Z)\), since the negative
rows vanish and there are no incoming higher differentials.
Consequently \(H^6(Y^{\rm an},\mathbb Z)/N^2\) has no two-torsion.
This quotient is a birational invariant of the complex fivefold
\cite[after Question~7.15]{BWWu}. Since \(X_\C\) is uniruled,
\(H^5_{\mathrm{nr}}(X(\C),\mathbb Q)=0\)
\cite[Proposition~7.4]{BWWu}. The level bound
\cite[Theorem~7.11]{BWWu} now gives \(s(\R(X))\le16\).

In the final case, Tian's characterization
\cite[Lemma~2.2]{TianSrcCodimOne} gives either a rationally connected
variety or a model fibered over a curve with rationally connected
geometric generic fiber. Over the curve the index is one by
Graber--Harris--Starr \cite{GHS}.
\end{proof}

\begin{corollary}\label{cor:real-sevenfold-level}
Let \(X\) be a smooth projective geometrically rationally connected
real sevenfold with \(X(\R)=\varnothing\).
If \(H_3(X(\C),\mathbb Z)[2]=0\), then \(-1\) is a sum of at most
sixty-four squares in \(\R(X)\).
\end{corollary}

\begin{proof}
The companion's degree-six Lawson theorem and its corollary on
codimension-two support
\cite[Theorem~11.10, p.~48, and Corollary~11.16, p.~53]{HW} give
\(H^8(X(\C),\mathbb Z)=N^2H^8(X(\C),\mathbb Z)\) under the
hypothesis on \(H_3\). Uniruledness of \(X_\C\) gives
\(H^7_{\mathrm{nr}}(X(\C),\mathbb Q)=0\)
\cite[Proposition~7.4]{BWWu}. Apply Benoist--Wittenberg's level bound
\cite[Theorem~7.11]{BWWu} with \(n=7\).
\end{proof}

\section{Zero-cycles over real curve function fields}\label{sec:real-zero-cycles}

Colliot-Th\'el\`ene's question concerns the degree-zero group over
complex curve function fields \cite[Question~11.3(a) in the preprint]{CTAlmostRational}.
His survey records vanishing for geometrically rational surfaces
and certain hypersurfaces and intersections of quadrics over
fields of cohomological dimension one. Pirutka used rational
simple connectedness for the complex function-field case
\cite[Theorem~1.4 and the discussion after Corollary~1.5 in arXiv v2]{Pirutka}.

We first record the restriction--norm argument and then prove the degree
isomorphism for point-free bases. We include proofs of the complex
curve-field input and of integral algebraicity on reduced special
fibers. The equivariant cycle comparison at disk collapse is the
additional step over the reals.

\subsection{Complex curve fields used for restriction and norm}

\begin{lemma}\label{lem:local-complex-chow-zero}
Let \(K\) be the function field of a smooth projective connected
complex curve and let \(V/K\) be smooth, projective and geometrically
rationally connected. Then \(\deg:\operatorname{CH}_0(V)\to\mathbb Z\)
is an isomorphism.
\end{lemma}

\begin{proof}
We give the curve-field argument also used in \cite{HW}, with its
inputs proved within this paper. We use the analytic construction and neck estimates from the proof of
Lemma~\ref{lem:eq-replace}, which precede its cycle comparison by norms.
For a smooth projective
rationally connected complex variety \(Y\), homologous one-cycles on
\(Y\) are algebraically equivalent. To check this, Tian--Zong's theorem
and comb smoothing write their difference, modulo algebraic equivalence,
as a signed sum of very free rational maps
\cite[Theorem~1.3 and Proposition~2.4]{TZ}. For the smoothing step,
attach very free teeth to a rational map and smooth. The smoothed rational
map and the teeth then express the original class as a difference of
very free classes. A very free rational handle through the marked
values of all signed maps exists by prescribing length-two jets in the
pointed rational-curve theorem \cite[Definition--Theorem~2.1(5)]{Kollar}.
Attach the positive and
negative maps separately to that same handle at distinct source
points. On either tree the normalization sequence gives vanishing of
the first cohomology of the pulled-back tangent bundle: each tooth can
prescribe its value at its attaching point. Blow up the attaching
points on the central fiber of \(\mathbb P^1\times\mathbb A^1\);
this supplies a fixed-source smoothing, which lifts in the relative
Hom scheme by the vanishing. The central and smoothed cycle classes
are algebraically equivalent by the graphs of this family over an
algebraic curve in the Hom scheme transverse to the central fiber.
The two smoothed maps have acyclic
pullbacks by semicontinuity and positive Euler characteristic by
their positive total anticanonical degrees. Their homology classes
agree, so by simple connectivity and Hurewicz they are homotopic.

Apply the fixed-source construction to this homotopy. To use it without
a real form of \(Y\), work on \(C\amalg\overline C\) with target
\(Y\amalg\overline Y\), where \(C=\mathbb P^1_\C\) and bars denote
the conjugate complex structures. The free involutions exchange the
two copies by the identity on the underlying manifolds. On the second
copy use the conjugates of the homotopy and endpoint operators from
the first. Equivalently, these are the complexifications of the same
varieties viewed over \(\R\) by composing their structure morphisms
with \(\operatorname{Spec}\C\to\operatorname{Spec}\R\). The analytic
part of the fixed-source construction on this disjoint union restricts
to the desired family on \(C\).
It preserves the endpoints and gives continuous families of graphs in
projective Hilbert schemes between disk switches. The endpoints of each
such family give algebraically equivalent cycles by taking chains of
algebraic curves in the connected parameter scheme. At a switch the
one-cycle of the Hilbert limit is the graph after collapse plus the
bubble, with its pushforward multiplicity, by the neck calculation
in Lemma~\ref{lem:eq-replace}. Each bubble class modulo algebraic
equivalence is constant on its active interval: its first-jet sphere
family, including the correction parameter at the disk center, gives
a continuous path of rational maps of fixed degree, and hence of
pushed-forward cycles in a projective Chow variety. Subtracting the
active bubble classes, with their multiplicities, therefore gives one
constant algebraic equivalence class through all switches. At the two
endpoints every disk is off. This proves the consequence.

Now let \(K=\C(t)\). Choose a nonempty finite subset \(D\) of
\(\mathbb P^1\) and a smooth projective model \(p:Y\to\mathbb P^1\)
with rationally connected smooth fibers over \(U=\mathbb P^1\setminus D\)
and reduced simple normal crossings boundary. The total space is
rationally connected: a section by Graber--Harris--Starr and rational
chains in good fibers connect general points, and rational chain
connectedness equals rational connectedness here. Put \(d=\dim V\).
The case \(d=0\) is immediate. Poincar\'e duality and Leray identify
\(H_2^{\rm BM}(Y_U,\mathbb Z)\) with \(\mathbb Z\) by generic degree.
Indeed, the base has the homotopy type of a graph, the fiber next-to-top
cohomology vanishes by simple connectivity, and the top cohomology
has its canonical orientation. Extend a given degree-zero zero-cycle
to a one-cycle \(Z\) on \(Y\). Localization lifts its homology class to
the reduced boundary. By Lemma~\ref{lem:local-complex-kato-two} a
vertical curve cycle \(C\) represents that lift. Thus \(Z-C\) is
algebraically trivial by the preceding consequence. Its restriction to
\(V\) vanishes: the algebraically trivial one-cycle group is divisible,
as it is generated by images from complex Jacobians, while the
degree-zero Chow group of \(V\) has bounded exponent by the diagonal
argument \cite[Proposition~11.1, preprint version]{CTAlmostRational}.
The class of \(C\) itself restricts to zero.

For general \(K\), take a finite inclusion \(\C(t)\subset K\). For a
finite extension \(L/K\) and two \(L\)-points of \(V_L\), their
corresponding points on \(\operatorname{Res}_{L/\C(t)}(V_L)\) have equal
zero-cycle classes by what we have proved. This restriction of scalars
is smooth projective and geometrically rationally connected, since
over an algebraic closure it is a product of the conjugate varieties.
Pull the equality back to \(L\) and use the proper evaluation map to
\(V_L\). Compare the resulting equality for each closed point of
\(V\), over its residue field, with a fixed \(K\)-point supplied by
Graber--Harris--Starr, then push down to \(V\). Every closed point thus
represents its degree times the fixed point class.
\end{proof}

\begin{corollary}\label{cor:real-function-field-zero-cycles}
Let \(B\) be a smooth projective connected curve over \(\R\), put \(F=\R(B)\), and let \(X\) be a smooth projective geometrically rationally connected variety over \(F\). Then the kernel \(A_0(X)\) of the degree map
\[
 \deg\colon\operatorname{CH}_0(X)\longrightarrow\mathbb Z
\]
is annihilated by two. The index of \(X\) divides two. Whenever \(X(F)\) is nonempty, the choice of \(x\in X(F)\) gives a splitting \(n\mapsto n[x]\), and
\[
 \operatorname{CH}_0(X)\simeq \mathbb Z\oplus A_0(X),
 \qquad 2A_0(X)=0.
\]
The splitting applies under the local--global hypothesis of Corollary~\ref{cor:rec}. If \(F\) contains \(\C\), the degree map is an isomorphism.
\end{corollary}

\begin{proof}
If \(F\) contains \(\C\), it is a complex curve function field, so the degree isomorphism follows from Lemma~\ref{lem:local-complex-chow-zero}. Otherwise \(E=F(\sqrt{-1})\) is a quadratic extension and is the function field of the smooth complex curve \(B_\C\). The same lemma gives \(\operatorname{CH}_0(X_E)\simeq\mathbb Z\) by degree. For the finite flat projection \(q\colon X_E\to X\) the operation \(q^*\) preserves the degree over the respective ground fields and \(q_*q^*=2\). Thus \(q^*\) sends \(A_0(X)\) to zero, which gives \(2A_0(X)=0\). Pushing forward a degree-one class on \(X_E\) yields a degree-two class on \(X\), proving the index assertion.

A rational point supplies the displayed splitting. Corollary~\ref{cor:rec} supplies it when the Colliot-Th\'el\`ene obstruction to the existence of rational points vanishes in the sense of \cite[Conjecture~2.6]{PS}.
\end{proof}

\begin{theorem}[Zero-cycles over point-free real curves]\label{thm:chow}
Let \(B\) be a smooth projective geometrically integral real curve
with \(B(\R)=\varnothing\). If \(X/\R(B)\) is smooth, projective
and geometrically rationally connected, the degree map
\(\deg:\operatorname{CH}_0(X)\to\mathbb Z\) is an isomorphism.
\end{theorem}

\subsection{Real families and the anisotropic conic}

The complex degree isomorphism and restriction followed by norm
show that the degree-zero group is killed by \(2\).
The real argument below uses integral twisted homology to choose
vertical corrections, and uses conjugation-equivariant disk replacement
to compare the corrected sections in the Chow group.

Let \(Q\) be the real conic \(x_0^2+x_1^2+x_2^2=0\).
For the proof of Theorem~\ref{thm:chow} it is enough to establish
equality in \(\operatorname{CH}_0(V)\) of the classes of any two rational points
on any smooth projective geometrically rationally connected
\(V/\R(Q)\). Indeed, Corollary~\ref{cor:points} applied to
the constant conic over \(B\) gives a nonconstant morphism
\(B\to Q\) by properness and geometric integrality. In particular
\(F=\R(B)\) is finite over \(F_0=\R(Q)\).
For a finite extension \(L/F\) and points \(a,b\in X(L)\), the
Weil restriction \(\operatorname{Res}_{L/F_0}(X_L)\) is smooth
projective and geometrically rationally connected,
since its geometric fiber is a product of conjugates of \(X_L\).
The two points correspond to \(F_0\)-points on this restriction.
Pull an equality of their zero-cycle classes back to \(L\) and
push along the proper evaluation map to \(X_L\). This gives
\([a]=[b]\) in \(\operatorname{CH}_0(X_L)\). A fixed \(F\)-point of \(X\)
exists by Corollary~\ref{cor:points}. Taking \(L\) to be the residue field of
each closed point, and then taking the norm, shows that its class
is its degree times the fixed rational-point class.

Fix \(V/F_0\) as above and a smooth projective model
\(p:W\to Q\), with a reduced simple normal crossings divisor over
the bad fibers. Such a model is obtained by resolution of a
projective closure and resolution of the vertical boundary.
The relative dimension-zero case is immediate, so assume
\(\dim V>0\).
The total space \(W_\C\) is rationally connected by the
Graber--Harris--Starr theorem \cite{GHS} and rational connectedness
of the general fibers. Also \(W(\R)=\varnothing\).
Write \(G=\operatorname{Gal}(\C/\R)\) and \(Y=W(\C)\).
Let \(N:\operatorname{CH}_1(W_\C)\to\operatorname{CH}_1(W)\) be proper pushforward from
the complexification, and write \(z_\eta\) for the restriction
of a one-cycle on \(W\) to \(V\).
Corollary~\ref{cor:real-function-field-zero-cycles} gives
\(2A_0(V)=0\), where \(A_0(V)=\ker(\deg:\operatorname{CH}_0(V)\to\mathbb Z)\).
Over \(F_0(\sqrt{-1})=\C(Q_\C)\), apply the degree isomorphism of
Lemma~\ref{lem:local-complex-chow-zero} and take norms to obtain
\begin{equation}\label{eq:normdegree}
      (N A)_\eta=(N B)_\eta
      \quad\text{if}\quad
      \deg(A_{\eta_\C})=\deg(B_{\eta_\C}).
\end{equation}

\begin{lemma}\label{lem:real-families}
Let \(M\) be a projective real scheme and let \(Z\) be a flat
family of one-dimensional subschemes of \(P\times W\) over \(M\),
where \(P\) is projective over \(\R\).
For \(m\in M(\R)\), push the top-dimensional cycle of \(Z_m\)
forward to \(W\) and restrict to \(V\). The resulting class in
\(\operatorname{CH}_0(V)\) is constant on connected components of \(M(\R)\).
\end{lemma}

\begin{proof}
For a smooth projective real curve \(T\) mapping to \(M\), the
pulled-back family gives a correspondence on divisor classes of
\(T\), by refined fiber intersection, proper pushforward to \(W\)
and restriction to \(V\). The class \([a]-[b]\) for
\(a,b\in T(\R)\) in the same component lies in
\(\operatorname{Pic}^0(T)(\R)^0\) by the Abel map based at a real point.
The real point ensures that the real Picard points are represented
by actual line bundles. This connected real Lie group is divisible
by two.
The image of its elements in \(\operatorname{CH}_0(V)\) has degree zero
as the degree homomorphism to \(\mathbb Z\) vanishes on this group.
Their image is therefore zero since \(2A_0(V)=0\).

Stratify \(M(\R)\) into a finite collection of connected smooth
semialgebraic pieces in smooth real algebraic loci, by successive
smooth and singular loci and real components. On each piece,
nearby points can be joined by an arc in an algebraic curve,
using inverse images of short line segments in an \'etale chart.
The result just proved shows local constancy on the piece.
A point in the closure of a piece is approached by a Nash arc
\cite[Proposition~8.1.13]{BCR}. The Zariski closure of the
nonconstant arc is a real algebraic curve. Normalization and
projective completion lift its branch to a real arc with limit
in the same component of a smooth projective real curve.
Properness of \(M\) extends the map at the endpoints.
Thus the value on a piece agrees with the value on its closure.
The finite stratification gives the assertion.
\end{proof}

In a flat family of this kind over a smooth real parameter curve,
the associated twisted homology classes on \(Y/G\) are likewise
constant along real paths. The equivariant class of the family,
after pushing forward from the projective factor, restricts to
the class of each fiber by compatibility with refined pullback
\cite[\S1]{BWhodgeI}. Along the real parameter path the quotient
is the product with \(Y/G\) and the local coefficients come
from \(Y/G\).

\begin{lemma}\label{lem:conic-obstruction}
Suppose \(Y\) is a simply connected manifold with a free
involution \(g\) generating \(G\).
Two equivariant maps \(Q(\C)\to Y\) are equivariantly homotopic
if their pushforwards of the twisted fundamental class agree in
\(H_2(Y/G,\mathbb Z(1))\). Here \(\mathbb Z(1)\) is the sign local system.
\end{lemma}

\begin{proof}
Use the free \(G\)-cellular structure lifted from the standard
one on \(Q(\C)/G=\mathbf{RP}^2\). Construct a homotopy over the
zero- and one-cell orbits, using simple connectedness.
The obstruction over the two-cell orbit lies in
\[
 H^2(\mathbf{RP}^2;\pi_2(Y))=
 H_2(Y,\mathbb Z)/\langle m+gm\rangle .
\]
The left side uses the action on \(\pi_2(Y)\) and the
identification on the right uses Hurewicz and the cellular
boundary \(1+g\). Changing the homotopy on the one-skeleton
changes the obstruction by this relation.
Gluing the two maps of the lifted cell along the homotopy
identifies the image of the obstruction under the natural map
\[
 H_2(Y,\mathbb Z(1))_G \longrightarrow H_2(Y/G,\mathbb Z(1))
\]
with the difference of twisted pushforwards. The lifts specify
the coefficient identifications. The map displayed is injective
by the Cartan--Leray spectral sequence. Its only possible incoming
differentials are from \(H_2(G,H_1(Y,\mathbb Z(1)))=0\) and
\(H_3(G,\mathbb Z(1))=0\). The latter vanishing follows from the
periodic resolution for the sign action. The obstruction is thus
zero. Equivariant smooth approximation relative to the endpoints
gives a smooth homotopy stationary near its endpoints.
\end{proof}

\subsection{Equivariant parametric replacement}

\begin{lemma}\label{lem:acyclic}
Let \(E\) be a smooth complex vector bundle on a compact Riemann
surface with a real structure covering a free antiholomorphic
involution. If \(\chi(E)\ge1\), two real Dolbeault operators with
\(H^1=0\) on \(E\) can be joined by a smooth family of such
operators, stationary near the endpoints.
\end{lemma}

\begin{proof}
Join the operators in the affine space of real Dolbeault
operators and use a stationary collar at each endpoint.
Serre duality identifies failure of \(H^1=0\) with a nonzero
holomorphic section of \(K_{\Sigma}\otimes E^\vee\).
A nonzero such kernel has a nonzero real section.
Work on the \(L^2\)-unit sphere of the real \(W^{1,p}\) sections,
for \(p>2\), and write \(D^\vee\psi=0\) for the normalized
dual equation, with values in real \(L^p\) forms.

The universal linearization along the normalized kernel is
surjective when one can vary \(D\) by real zero-order forms.
Indeed, restricting the operator on \(\psi\) to the tangent
sphere leaves its image unchanged, as the omitted direction
\(\psi\) itself lies in the kernel. A vector annihilating the
image solves the adjoint equation and is smooth by elliptic
regularity. On a small disk where \(\psi\ne0\), any variation in
the equation can be made by a zero-order variation of \(D\)
with support there and on the conjugate disk. The two real
pairings agree. Thus the annihilator vanishes on this open set
and then everywhere by unique continuation for the adjoint
Dolbeault equation.

The normalized kernels along the original path form a compact
set by uniform elliptic estimates and Rellich compactness.
Surjectivity of the augmented linearization can therefore be
arranged by finitely many smooth perturbation directions with
support away from the endpoint collars. For sufficiently small
coefficients the universal equation remains transverse along
all its zeros, since any sequence of normalized zeros for
coefficients tending to zero has a convergent subsequence.
If the number of perturbation directions is \(m\), its zero
manifold has real dimension
\[
 m+1+\chi(K_\Sigma\otimes E^\vee)-1=m-\chi(E)<m.
\]
The \(1\) accounts for the path parameter and the \(-1\) for
normalization. The indices here are real indices on the fixed
spaces of conjugation and hence equal the complex indices.
By Sard's theorem the projection of this manifold to
\(\R^m\) has measure zero. A sufficiently small choice of
coefficients outside its image gives the required path.
\end{proof}

\begin{lemma}\label{lem:eq-replace}
In the situation of \(W\to Q\) above, suppose two real morphisms
\(u_0,u_1:Q\to W\) are equivariantly homotopic and have
\(H^1(Q_\C,u_i^*T_{W_\C})=0\) with
\(\chi(u_i^*T_{W_\C})\ge1\).
Then
\((u_{0*}[Q])_\eta=(u_{1*}[Q])_\eta\) in \(\operatorname{CH}_0(V)\).
\end{lemma}

\begin{proof}
Smooth the equivariant homotopy relative to stationary collars.
Its pullback tangent bundle admits an equivariant smooth
trivialization in the parameter direction by a real connection.
Apply Lemma~\ref{lem:acyclic} to join the endpoint Dolbeault
operators, keeping the natural operators on the collars.

Use the normalized equation of Proposition~\ref{prop:analytic}
with the operators \(D_t\) of Lemma~\ref{lem:acyclic}. After smooth
identification of the bundles over the interval, a right inverse
\(Q_{t_0}\) at \(t_0\) gives the right inverses
\(Q_{t_0}(D_tQ_{t_0})^{-1}\) near \(t_0\) by a Neumann series.
Use a finite partition of unity in \(t\) to obtain uniformly bounded
right inverses continuous in the operator norm \(L^4\to W^{1,4}\).
Average them under conjugation.

For clarity, the following fixed-source version of the disk construction
in the proof of Lemma~\ref{lem:real-chow-vary} is the one used here.
The source is \(Q_\C\simeq\mathbf P^1_\C\) with fixed complex structure,
and its real involution is free. Choose a finite disjoint packing in
the quotient with arbitrarily small radii and lift it in conjugate
pairs. On a disk of radius \(r\), with fixed holomorphic coordinate
\(\zeta\), the local addition in a moving holomorphic bundle frame has
Taylor terms \(m(q)+a(q)\zeta+b(q)\bar\zeta\). Use the core formula
\[
    G\bigl(\lambda r^2/\zeta; m(q)+a(q)\zeta,b(q)\bigr).
\]
Here \(G\) is a holomorphic first-jet family from
\cite[Proposition~3.5]{DGX}, chosen on successive parameter intervals.
When the factor \(\lambda\) is zero the core is the holomorphic
Taylor expression \(m(q)+a(q)\zeta\). At full scale \(\lambda=1\),
on the annulus of width \(\theta r\) the mismatch is
\(O(\theta r+r^2)\), including first and second fiber derivatives on
a smaller ball. The bound follows from
\(\bar\zeta-r^2/\zeta=O(\theta r)\) and Taylor expansion. Taking
\(r\le\theta\), normalization as in Proposition~\ref{prop:analytic}
gives uniform estimates for the coefficients and their fiber derivatives
on these annuli. The standard annuli and unpacked set have total area
\(O(\theta)\).

To change the local family, widen the annulus to inner radius
\(r/2\), take \(\lambda\) to zero, then turn the new family on and
narrow the annulus. At zero scale the two core expressions agree
if the same target chart is used for the switch. If a chart change is
also needed, express both Taylor approximants in one common chart and
interpolate there with \(\lambda=0\). The two are holomorphic on the
core and differ by \(O(r^2)\), since their value and linear terms agree.
The widened annulus then has a uniform bound depending on \(\theta\)
for the coefficients and their fiber derivatives. Take sufficiently
small radii and separate the switch intervals so that any one widened
pair contributes \(O(\theta^{1/4})\) to the \(L^4\) estimate. On a
stationary endpoint collar switch the pairs off successively. After
setting \(\lambda=0\) on a widened pair, multiply its gluing cutoff by
a scalar decreasing to zero. This interpolates between the holomorphic
Taylor expression and the original local addition in one target chart.
The difference is \(O(r)\), also after two fiber derivatives, and the
cutoff derivative is \(O(1/r)\). The inverse fiber derivative remains
bounded, so the same small-support estimate applies on the whole pair
while the cutoff is being removed. Reversing this procedure turns a
pair on. Where no modification is active on the stationary collars,
the natural endpoint operator makes the constant and linear terms vanish.

On the holomorphic cores the normalized coefficient is zero, including
at \(\lambda=0\) before cutoff removal. On the annuli and during cutoff
removal its local expression varies continuously. The small-support contraction of
Proposition~\ref{prop:analytic} therefore gives a continuous family of
corrections in \(C^0\) and \(W^{1,4}\), unique on the fixed small
parameter ball. The paired construction and invariant right inverses
make the solutions equivariant. Reconstruction is holomorphic by
the chain rule on the annuli and by joint holomorphy on the cores.
When all disks are off on the endpoint collars, the zero parameter
is the unique solution, so the endpoint maps are exactly \(u_0,u_1\).

Here is the cycle comparison at the switches. On an interval
with a fixed choice of active branches the corrected maps
\(v_t:Q\to W\) are real and give a continuous path of graphs
in a projective Hilbert scheme of \(Q\times W\). Indeed, the
holomorphic maps vary uniformly on such an interval, so their
degrees are locally constant and their graphs have a fixed
Hilbert polynomial. A Hilbert limit there has the limiting smooth
graph as its one-cycle. The polynomial rules out an additional
zero-dimensional part, so the limit is that graph.
By Lemma~\ref{lem:real-families} their generic cycle class is
constant on that interval. At a collapsing disk use the local
coordinate \(\zeta\) centered at the disk and the rescaling
\(w=\mu/\zeta\) in the replacement formula. Off the center
the maps tend to the main map, and in the rescaled chart
they tend to one holomorphic rational bubble. The correction
section is holomorphic in a fixed local frame on the widened core,
so Cauchy estimates bound its derivative near the collapsing center.
In the transition annulus
\(\mu/\delta\le|\zeta|\le\delta\), the sphere coordinate satisfies
\(|\mu/\zeta|\le\delta\). The chain rule and the derivative bound on
the correction give an upper bound \(C(1+\mu/|\zeta|^2)\) on the map
derivative. Integrating its square in polar coordinates gives image
area \(O(\delta^2)\). Integrating the derivative along radial segments
and an arc of the outer circle gives image diameter \(O(\delta)\),
since \(\mu\le\delta^2\). Thus the two ends have the same limit, and the one-cycle of a real Hilbert limit consists of the
graph of the limiting main map and the bubbles over the source
points, with their pushforward multiplicities. This also
follows by applying the positive-current limit for the
replacement to the graph in \(Q\times W\), then using the
Hilbert--Chow morphism from the reduced seminormalization of the
Hilbert scheme \cite[Definition--Theorem~3.1(6)]{Ross}. The finite
seminormalization map induces homeomorphisms on real and complex
points, so this applies to the Hilbert limits used here. On the graph
the bubble projection to the source is constant. In \(W\) its projection can be horizontal
for \(p\).

Each of these bubbles occurs with its conjugate, so its
contribution on the generic fiber is the norm of a complex
zero-cycle. Along a branch of a disk family the degree of its
bubble over the target \(Q_\C\) is constant, by homotopy and
proper pushforward, including at the branch limit.
Equation~\eqref{eq:normdegree} therefore gives a constant
generic class for each bubble pair along that branch.
Subtract the active bubble contributions from the generic
cycle class of the main graph. Lemma~\ref{lem:real-families},
applied also to its one-sided Hilbert limits, gives the same
corrected class on both sides of every switch. At both
endpoints no bubbles are subtracted. This proves the assertion.
\end{proof}

\subsection{Vertical correction and proof of Theorem~\ref{thm:chow}}

Extend the two \(F_0\)-points of \(V\) to sections \(s_i:Q\to W\).
Attach to each section conjugate pairs of general very free
vertical teeth in distinct good fibers and smooth, using
\cite[Lemma~25 and Proposition~24 in the authors' revised version]{HT}. Here the open
conditions on the attachment points and tangent directions can
be satisfied by conjugate choices. In the product of a complex
parameter chart and its conjugate the paired locus is Zariski dense.
One may work in \(W\times\mathbf P^N\) over \(Q\) to use embedded
combs. The normal-sheaf vanishing makes the Hilbert scheme
smooth, and one can choose a real tangent vector smoothing
all nodes because this is a nonempty open condition. A real smooth
curve in the Hilbert scheme with that tangent gives nearby smooth
fibers. Degree one over \(Q\) is preserved, so the resulting maps
are sections. The normal-sheaf vanishing persists by semicontinuity.
For a smoothed section in the auxiliary product, its normal sheaf
is the sum of the pullbacks of \(T_{W/Q}\) and \(T_{\mathbf P^N}\).
Projection to \(W\) therefore preserves the needed vanishing.
By the tangent exact sequence of a section and
\(H^1(Q_\C,T_{Q_\C})=0\), the new sections \(v_i\) satisfy
\(H^1(Q_\C,v_i^*T_{W_\C})=0\).
Lemma~\ref{lem:real-families} identifies their generic Chow classes
with those of the original sections, since the added curves are
vertical.

Choose a good affine open \(U\subset Q\) with rationally
connected complex fibers and let \(n=\dim V\).
Poincar\'e duality with local coefficients on the free quotient
gives
\[
 H_2^{\mathrm{BM}}(W_U(\C)/G,\mathbb Z(1))
     = H^{2n}(W_U(\C)/G,\mathbb Z(n)).
\]
The orientation system of the total space is \(\mathbb Z(n+1)\).
The base \(U(\C)/G\) has the homotopy type of a graph. In the
Leray spectral sequence \(R^{2n-1}\) vanishes by simple
connectedness and duality on the fibers. The degree \(2n\)
cohomology with its \(\mathbb Z(n)\) twist is constant with value \(\mathbb Z\).
Thus restriction to any good fiber identifies the displayed
group with \(\mathbb Z\), and takes the class of a one-cycle to its
generic degree. The difference of the section classes has degree
zero. Borel--Moore localization lifts it to the homology of the
boundary \((W\setminus W_U)(\C)/G\). Each quotient fiber there
identifies with one of two conjugate complex reduced fibers,
and the sign local system trivializes by choosing that lift.
Lemma~\ref{lem:local-complex-kato-two} applies to this model over
the algebraic curve \(Q_\C\), since the reduced fibers above the
boundary form a simple normal crossings divisor. It supplies a complex vertical
one-cycle \(C\) whose norm gives the required twisted class.
The norm here corresponds to pushforward from a single
complex fiber in each pair, with its sign trivialization.

Write the class of \(C\) as a difference of sums of classes of
very free rational maps into \(W_\C\).
This follows from Tian--Zong's generation by rational curves and
their very free comb deformation \cite[Theorem~1.3 and
Proposition~2.4]{TZ}. Attach the negative summands to \(v_0\)
and the positive ones to \(v_1\), if the class of \(N C\)
is \(\operatorname{cl}(v_{0*}[Q]-v_{1*}[Q])\).
To arrange the attaching points, for each summand choose
one very free auxiliary map through the two proposed node
values and a point on the summand. Very free interpolation
allows these incidences at distinct source points.
On the side requiring the summand, smooth its union with
the auxiliary curve while fixing the relevant node value
on the auxiliary map. For this smoothing \(H^1\) of the
pullback tangent bundle with a double twist at that mark
vanishes by normalization and very freeness. Hence it
can be smoothed fixing the value and the smoothed map is
still very free. On the other side use the auxiliary map
directly. The added auxiliary classes therefore cancel.
Use the conjugate of each tooth at the conjugate node.

These two absolute combs in \(W\) have acyclic tangent
pullbacks by normalization, since the handles are acyclic
and the teeth supply all node values. They can now be smoothed
with fixed smooth source \(Q\): blow up \(Q\times\mathbf A^1\)
at the attaching closed points on the central fiber and
use the relative Hom scheme into \(W\). It is smooth over
the parameter at the comb by the \(H^1\)-vanishing.
Taking a real algebraic curve through that point transverse
to the central fiber gives real smoothings \(w_i:Q\to W\).
Their tangent pullbacks are acyclic by semicontinuity.
Their degrees over \(Q\) are positive since each handle
has degree one and all teeth have nonnegative degree.
In particular the differential gives a nonzero map
\(\mathcal O_{\mathbf P^1}(2)\to w_i^*T_{W_\C}\) over \(Q_\C\), so
\(H^0(Q_\C,w_i^*T_{W_\C})\ne0\). Acyclicity in degree one
then gives \(\chi\ge1\).

The central and smoothed fibers are in the same real
parameter component. Their cycle classes specialize by
flatness, and their generic Chow classes agree by
Lemma~\ref{lem:real-families}, applied to graphs in a
relative projective embedding of the source family.
Thus \(w_0,w_1\) have equal twisted homology classes.
Moreover, the difference of the attached norm cycles
has generic Chow class zero. After cancellation of the auxiliary
summands, the complex teeth added on the second side minus those
added on the first side have the homology class of \(C\).
Their relative degree is the pairing with
\(p_\C^*c_1(\cO_{Q_\C}(1))\), which vanishes on the vertical
cycle \(C\). Now \eqref{eq:normdegree} applies to the two
sums of complex teeth.
The difference of the generic classes of the new maps
therefore equals \([s_0(\eta)]-[s_1(\eta)]\).

Lemma~\ref{lem:conic-obstruction} gives an equivariant
homotopy between \(w_0\) and \(w_1\). Lemma~\ref{lem:eq-replace}
makes their generic Chow classes equal. This proves equality of
the two original point classes. The Weil restriction argument
completes the proof of Theorem~\ref{thm:chow}.

\subsection{Local classes and closed points of fixed degree}
\label{sec:zero-cycle-local-global}

Let \(B\) be a smooth projective geometrically integral curve over
\(\R\), and put \(F=\R(B)\). Let \(X/F\) be smooth, projective and
geometrically rationally connected of dimension \(n>0\). Choose a regular
projective flat
model \(f\colon\cX\to B\). A place will be called good if \(f\) is smooth
there and the special fiber is geometrically rationally connected. Write
\(F_b\) for the completion at a closed point \(b\) of \(B\).

Let \(T\subset B(\R)\) be finite, contain every bad real place, and meet
every real component of \(B\). For each interval \(I\) of
\(B(\R)\setminus T\), the real fiber components are identified along \(I\)
by the smooth proper fibration. Each end of \(I\) at \(b\in T\) specifies
one of the two orderings of \(F_b=\R((t))\), according to the sign of a
local parameter \(t\). After extension to the corresponding real closure,
a zero-cycle on \(X_{F_b}\) gives a vector over \(\mathbb Z/2\) by counting its
real points with multiplicity on each semialgebraic component. These
components are identified with the fiber components along \(I\), by
semialgebraic triviality near the end. The two vectors associated with a
marked point can differ at a bad fiber. We call the data compatible if
the two end vectors agree on every interval, counting both ends if the
same marked point occurs twice.

\begin{corollary}[Interpolation by a single closed point]
\label{cor:closed-point-interpolation}
Let \(N\geq1\). Prescribe a reduced effective zero-cycle of degree \(N\) at
each \(b\in T\) and at any finite set of non-real closed points of \(B\).
There exists a closed point \(P\in X\) of degree \(N\) inducing the
prescribed local Chow classes and the prescribed local finite
\(\acute{\mathrm e}\)tale support algebras if and only if the real data
are compatible. One may also prescribe any finite jet of the
corresponding formal Hilbert-scheme sections on a regular projective
model as in the proof. For \(N\geq2\), the normal closure of \(F(P)/F\)
may be required to have group \(S_N\) and to be linearly disjoint
from any given finite extension of \(F\).
\end{corollary}

For a relative unramified class \(h\in
H^q_{\mathrm{nr}}(F(X)/F,\mathbb Z/2)\), pairing with a local
zero-cycle \(z_b\) uses evaluation on the residue fields of its points
followed by corestriction to \(F_b\). This descends to
\(\operatorname{CH}_0(X_{F_b})\) by the pairing for cycle modules
\cite[Section~2.4]{MerkurjevUnram}. For a connected component \(S\) of
\(B(\R)\) put
\[
 \rho_{S,h}((z_b))=\sum_{b\in S}\partial_b\langle z_b,h\rangle ,
\]
where \(\partial_b\) is the Laurent-series residue. For a fixed \(h\)
only finitely many terms are nonzero. The Gersten resolution identifies
the residue-free classes on a smooth model with sections of the Zariski
sheaf \(\mathcal H^q\) associated to \(H^q_{\acute{\mathrm e}t}(-,\mathbb Z/2)\)
\cite[Theorem~4.2(2)]{BlochOgus}. All horizontal residues of \(h\) vanish,
and only finitely many vertical residues can be nonzero. Deleting the
corresponding fibers and the bad fibers therefore extends \(h\) as a
section of \(\mathcal H^q\) on the remaining model. At any remaining
place, properness extends a point over a finite extension of \(F_b\)
to its valuation ring. A Zariski-local \'etale representative of the
section near the center pulls back to a class on that ring. Its
evaluation is unramified, as is its corestriction by compatibility
of corestriction with residues.

\begin{corollary}[The integral diagonal image]
\label{cor:zero-cycle-diagonal}
A tuple in \(\prod_b\operatorname{CH}_0(X_{F_b})\) is the image of an
element of \(\operatorname{CH}_0(X)\) if and only if its degrees are all
one integer \(d\) and \(\rho_{S,h}=0\) for every \(S\) and every
\(h\in H^{n+2}_{\mathrm{nr}}(F(X)/F,\mathbb Z/2)\). The same criterion
holds with \(n+1\) in place of \(n+2\). Equivalently, require the
common degree and compatibility at the real ends for every finite
refinement \(T'\supseteq T\) of one such \(T\).
Fix the degree-two effective zero-cycle \(\nu\) obtained by
corestriction of an \(F(\sqrt{-1})\)-point. A global preimage can be
chosen of the form \([P]-a\nu\) for every sufficiently large integer
\(a\), where \(P\) is a single closed point of degree \(d+2a\) with
the Galois and disjointness conditions of
Corollary~\ref{cor:closed-point-interpolation}.
\end{corollary}

\paragraph{Compatibility and reciprocity.}
For \(q\geq2\), the two real-closure evaluations of a class in
\(H^q(\R((t)),\mathbb Z/2)\) differ by its residue. This follows from
the decomposition as the sum of a constant class and \( (t)\) cupped
with a constant class of degree \(q-1\). Corestriction is compatible
with passage to the two real closures. Contributions from non-real
factors vanish, while real factors contribute once per real embedding.
At a good place a real residue point of ramification index \(e\)
contributes parity \(e\) to the same fiber component on each side.
Indeed, in a normalized real branch the base parameter has the form
\(\pm s^e\) after absorbing a real unit. This has one real branch on
both sides for odd \(e\), and two on one side for even \(e\).
Thus the two vectors of a good real place agree with one another.
Compatibility after inserting that place also requires its vector
to equal the vector transported from the old interval ends.

For a section of \(\mathcal H^q\) on a smooth open model, the evaluation
on its real points is locally constant and transfers to real closed
extensions. Indeed, by the Gersten theorem for Milnor \(K\)-theory
and the mod-two norm residue isomorphism, the sheaf \(\mathcal H^q\)
is locally generated by Kummer symbols of regular units
\cite[proof of Proposition~2.2]{BenoistHilbert17}. On a finite Zariski
cover these expressions evaluate over any real closed extension by
the same finite sign tests, using the generator \((-1)^{\cup q}\). The
assertion that those tests have a fixed value on a given semialgebraic
component therefore transfers to every real closed extension.
Semialgebraic triviality over a sufficiently short interval at each
end identifies the component labels by such formulas. In an ordered
real closure of \(F_b\), the base point is infinitesimally inside the
corresponding interval, so the value on a formal branch is the value on
its labeled component. This applies in particular to the sections
induced by global \(\acute{\mathrm e}\)tale classes.

To separate the end vectors for a finite refinement \(T'\), take
the smooth open model over
\(B\setminus(T'\cup\{\text{bad places}\})\). For a real component of
this open model the mod-two equivariant Thom class of the fixed
locus yields a global class of degree \(n+1\) taking the
characteristic function of that component as its real evaluation.
Indeed, the real normal bundle has \(n+1\) sign directions, so its
restriction over a real point has equivariant Euler class
\(u^{n+1}\), where \(u\) generates
\(H^1(\operatorname{Gal}(\C/\R),\mathbb Z/2)\).
Use the Thom class on the specified component and zero on the
others. The equivariant comparison identifies this with an
\(\acute{\mathrm e}\)tale class on the open model and hence produces
a relative unramified class. Cupping with the real constant
\((-1)\) gives the same evaluation in degree \(n+2\).
See \cite[Remark~1.6(ii) and (1.36)--(1.39)]{BWhodgeI}.
Each such component lies over an interval and represents one
continuously identified fiber component. The resulting
\(\rho_{S,h}\) is the sum of the two end parities on that component.
This also shows, using the evaluation pairing, that the end
vectors factor through the local Chow group.

Conversely, add the real points of the deleted fibers for \(h\) to
\(T\) to obtain a finite refinement. Apply the preceding evaluation
formula to the section of \(\mathcal H^q\) on the remaining model.
The two-ordering formula writes \(\rho_{S,h}\) as a sum of end parities
weighted by the constant evaluations on the components over these
intervals. Compatibility for this refinement makes the sum zero.
Necessity for every global zero-cycle now follows from the
componentwise reciprocity law
\cite[Proposition~2.3]{PS} applied to its global evaluation class.

\begin{proof}[Proof of Corollary~\ref{cor:closed-point-interpolation}]
Take a smooth projective resolution \(\mathcal H_N\) of the main
component of \(\operatorname{Hilb}^N(\cX/B)\), preserving its
smooth locus over the good fibers. Its generic fiber is
geometrically rationally connected, since it is birational to the
symmetric power of \(X\), dominated by \(X^N\). The generic
reduced-configuration locus is also preserved. Every prescribed
configuration lifts to a formal section of \(\mathcal H_N\) by
properness. Its center lies in the relative smooth locus. In the
regular local ring at the center, the pullback of the base
parameter cannot belong to the square of the maximal ideal,
since the section sends it to a uniformizer. Relative
\(\acute{\mathrm e}\)tale coordinates allow replacement of each
formal section by a convergent holomorphic section germ
with an arbitrarily high prescribed jet. Over a real point the
germ can be chosen equivariant, and over a non-real point one
chooses the conjugate germ on the other disk.

These germs can be chosen in the required local Chow classes.
In the reduced locus of the Hilbert scheme, its universal support
is finite \(\acute{\mathrm e}\)tale. Its isomorphism scheme with the
given local support algebra is finite \(\acute{\mathrm e}\)tale and
has a point over the prescribed configuration. Consequently nearby
configurations in the valuation topology have the same support
algebra and their corresponding points are nearby over each fixed
residue field. The \(R\)-equivalence classes are open over these
complete discrete valuation fields by
\cite[Remark~1.10]{KollarLocal}. Corestriction preserves the
resulting equality of Chow classes. Sufficiently high jet interpolation
preserves the class, the support algebra and both end vectors. For the
end vectors, view the real embeddings of each support field in the two
ordered real closures of \(F_b\). Each semialgebraic fiber component is
open in its fiber over the real closure. In an affine chart about a
support point, choose a positive radius whose ball in the fiber stays
in that component. Each radius has finite Puiseux valuation. Under the
support-algebra isomorphism, a sufficiently high jet makes every
coordinate difference have valuation greater than the valuations of
these finitely many radii, preserving all the component labels. For a convergent germ,
the normalization of its finite analytic support gives convergent local
parametrizations with \(t=\pm s^e\) on real branches. Membership in each
component over a sufficiently short one-sided interval is determined by
finitely many semialgebraic sign conditions. The signs after substitution
are the signs of the formal Puiseux series in the corresponding ordering.
The analytic germs therefore have the prescribed end vectors on actual
real fibers near the marked point.

On each good real interval, use the equivariant transport on the
underlying fibers from the proof of
Proposition~\ref{prop:real-hodge-fibration} to bring the reduced
ends to a common fiber. Join them there by
\hyperref[par:real-hilbert-paths]{the real configuration path}.
Its reduced and curvilinear charts lie in the relative smooth
locus of the main component, hence lift unchanged to \(\mathcal H_N\).
A smooth horizontal lift on a neighborhood of the compact path
gives the local product used in that interpolation. Near the two
reduced endpoints, choose it to agree with the configuration
transport, extending by smooth cutoffs in the relative smooth
locus. Keeping smaller collars of the analytic section germs gives
a smooth section over the whole real locus in the submersion locus
of the model. By \cite[Proposition~2.3]{BW} we may approximate it by
an equivariant holomorphic section on a neighborhood of \(B(\R)\),
preserving the real jets. Adjoin disjoint disks at the nonreal marked
places, as in \cite[proof of Proposition~2.9]{BW}, and apply
Theorem~\ref{thm:main} to realize all
these jets by an algebraic section.

At an additional good non-real place prescribe an irreducible
local support algebra of degree \(N\). Such places can be chosen
outside any finite set. Over \(\C((t))\) all extensions are cyclic,
and each factor of degree \(e\) can be obtained in a smooth
relative coordinate chart by taking one coordinate equal to
\(t^{1/e}\) over its totally ramified extension. Distinct factors
may be taken at distinct complex residual points. This test forces
the generic support to be one closed point. For \(N\geq3\) also
prescribe at two further good non-real places the decomposition
types \((N-1,1)\) and \((2,1,\ldots,1)\). The inertia in the normal
closure contains an \((N-1)\)-cycle fixing one letter and a
transposition. A transitive subgroup with such an \((N-1)\)-cycle
is doubly transitive, and the conjugates of the transposition then
give every transposition. For \(N=2\) the irreducible quadratic test
itself supplies the transposition. Choose the transposition place to
be unramified
in the normal closure of the given finite extension of \(F\).
The corresponding transposition is killed on the intersection
of the two normal closures. Since its normal closure in \(S_N\)
is all of \(S_N\), the intersection is \(F\). This gives the
asserted disjointness.
\end{proof}

\begin{proof}[Proof of Corollary~\ref{cor:zero-cycle-diagonal}]
We first describe how finitely many interpolated places detect a
tuple of common degree that is compatible for every finite
refinement. At a good place, Koll\'ar's
specialization isomorphism identifies the local Chow group with
the Chow group of the special fiber
\cite[Theorem~2]{KollarSpecialization}. The completed local ring
is Henselian, its residue field is perfect and the special fiber
is smooth proper and separably rationally connected. For a complex
fiber the degree suffices. For a real fiber use the degree and the
parity on each real component. In fact two real points on the
same component are \(R\)-equivalent
\cite[Corollary~1.7]{KollarLocal}, and a non-real point contributes
the norm of a degree-one complex class. If there is a real point
this norm is twice that point. If there is no real point every
point is a norm and the degree distinguishes all classes.
Specialization sends a local real residual branch of ramification
index \(e\) to \(e\) times its residual point, so its parity equals
the common end parity. For any good real point \(b\notin T\), apply
compatibility to \(T\cup\{b\}\). It identifies that parity with the
vector transported from the ends of the interval containing \(b\).
Hence two tuples of the same common degree, both compatible for
every refinement, agree everywhere as soon as they agree at all
bad places and all marks of \(T\). If this set is empty, one may
add a place to prescribe the degree.

The Graber--Harris--Starr theorem \cite{GHS} gives a point over
\(F(\sqrt{-1})\) and hence \(\nu\). Its end parities vanish.
There is an integer \(M>0\), depending only on \(X\) and \(\nu\),
such that
\[
  2M[z]=M\deg(z)\nu_L\qquad
  \text{for every zero-cycle \(z\) on \(X_L\) and every \(L/F\)}.
\]
Here and below brackets denote Chow classes. Over an algebraic
closure of \(F(X)\), rational connectedness gives
\(2[\eta]=\nu\) for the generic point \(\eta\). The relation descends
to a finite extension, so its norm gives \(2M[\eta]=M\nu\) over
\(F(X)\). By Chow localization, \(2M\Delta_X-M[X]\times\nu\)
is supported over a proper closed subset in the first factor.
The moving lemma applied to zero-cycles over \(L\) gives the
displayed identity.

For each place in the chosen finite set, start from a representative
of \(z_b\) and add sufficiently large multiples of the negative
points, in multiples of \(2M\). The identity expresses the result
as an effective representative of \(z_b+a_b\nu_{F_b}\).
Add further copies of \(\nu_{F_b}\) to make \(a_b=a\) for a common
arbitrarily large \(a\). The local effective representatives can
be chosen reduced in their Chow classes. Over each residue field
move repeated points separately in their open \(R\)-equivalence
class, keeping distinct conjugates and the original residue
field. The latter is open near a point of exact residue field,
as a function whose value generates the finite extension can be
kept away from its finitely many proper subfields. The smooth
local coordinate chart has positive dimension.

By compatibility and the vanishing end parities of \(\nu\),
Corollary~\ref{cor:closed-point-interpolation} applies to these
effective data. Its point \(P\) gives \([P]-a\nu\) in the prescribed
classes at the test places. Global reciprocity makes this new
tuple compatible for every finite refinement. The finite detection
just established proves
equality of the tuples at all places.
\end{proof}

\begin{corollary}[Degrees of closed points with good real reduction]
\label{cor:real-closed-point-degrees}
Assume all real places are good. For a component \(S\) of \(B(\R)\)
let \(q_S\) be the least odd length of a monodromy orbit on the
connected components of a real fiber. Put \(q_S=\infty\) if there
is no odd orbit. The positive degrees of closed points on \(X\)
are precisely
\[
  \{2,4,6,\ldots\}\ \cup\
  \{N>0:\ N\text{ is odd and } N\geq q_S
                             \text{ for each } S\}.
\]
For every allowed degree \(N\geq2\) one can impose the Galois and
linear-disjointness conditions in
Corollary~\ref{cor:closed-point-interpolation}.
\end{corollary}

\begin{proof}
A global effective zero-cycle gives a monodromy-invariant parity
vector on each real circle. If \(N\) is odd, a set of components
with parity one must contain an odd orbit. Each represented
component requires at least one of the \(N\) support points.
Conversely choose a good marked point on each circle. For even
\(N\) take conjugate pairs only. For odd \(N\) choose one real
point on each component of a shortest odd orbit and fill the
remaining degree with conjugate pairs. The chosen points and
pairs lift at these good places to reduced local configurations.
The transported parity is invariant under the monodromy, so
Corollary~\ref{cor:closed-point-interpolation} gives the required
closed point. If the real locus is empty only the non-real
tests of that corollary are needed.
\end{proof}

\begin{corollary}[Degrees of closed points with arbitrary real reduction]
\label{cor:real-closed-point-graph}
The spectrum of degrees in Corollary~\ref{cor:real-closed-point-degrees}
has the same form without the good real reduction assumption, with the
least odd orbit length replaced by the following finite graph invariant.
Resolve the bad real fibers to normal crossings, separating conjugate
geometric components at real points. Take \(T\) as above. On each
circle \(S\) the vertices of the graph are the connected components of
the real fibers identified over each interval of \(S\setminus T\).
For each \(b\in S\cap T\) and each connected component of the real
smooth locus of the reduced fiber, take an edge. Its two ends are
the fiber components reached by a real transverse arc on either side
of the given smooth locus. Attach the weight given by that fiber
component's multiplicity. Keep track of which sides of \(b\) its two
ends occupy. In particular both ends may be on the same side.

For a set of edges with even incidence at every vertex, its mod-two
winding around \(S\) is the parity of the number of edges crossing
from one side to the other at any marked point. This parity is
independent of the marked point. Let \(q_S\) be the minimum, over
such edge sets of odd winding, of the maximum over marked points
of the sum of the selected weights there. If there is no such edge
set, put \(q_S=\infty\). Then all positive even degrees occur, and
an odd degree \(N\) occurs exactly when \(N\geq q_S\) for every
\(S\). In the minimum it suffices to use a single simple graph
cycle of odd winding for each \(S\). The Galois and disjointness
conditions of Corollary~\ref{cor:closed-point-interpolation} are
available for all permitted \(N\geq2\).
\end{corollary}

\begin{proof}
The graph may have loops and parallel edges. A loop counts as a simple
cycle and contributes twice to the incidence at its vertex. The indicated
model can be obtained by log resolution and blowing up the intersections
of any geometrically distinct conjugate
components which meet at real points. In normal crossings these are
smooth centers and the two strict transforms are disjoint. The
exceptional components through real points are defined over \(\R\).
These modifications leave good fibers unchanged. At a real point
of the resulting fiber, complete the local algebraic equations
\(z_i\) of its real components to an \'etale coordinate system.
Here \(f^*t=U\prod_i z_i^{m_i}\). On a small real coordinate box
the local inverse is Nash and \(U\) has constant sign. Thus the
chamber descriptions apply in the real closures for both
orderings of the formal branch field. Every chamber with prescribed nonzero
coordinate signs is connected and maps into one neighboring real
interval, so lies in one of its fiber components. The two resulting
labels of a transverse arc are locally constant, up to exchange,
on a real smooth stratum. The graph is finite by semialgebraic
finiteness.

A closed point of \(X_{F_b}\) with real residual field has a
normalized branch \(z_i(s)=s^{e_i}u_i(s)\) for the coordinates
vanishing at its center. Its degree is \(\sum_i m_i e_i\).
Its two endpoint sign chambers differ exactly in the coordinates
with \(e_i\) odd. Change those signs one by one inside the same
box. Each step crosses the real smooth locus of one divisor and
gives one graph edge of weight \(m_i\). The intermediate chambers
cancel in pairs in the end vector. The weight sum has therefore
decreased by an even number while preserving these parities.
A local point with non-real residual field contributes only an
even degree and no end parity. For a sum of local points, repeated
edges may be deleted in pairs, decreasing the weight by another
even number. Thus an effective local configuration gives a set
of edges of weight no more than its degree and of the same parity.

Conversely each edge is realized by a formal arc transverse to
the corresponding fiber component at a real smooth point. Set
\(z_i=s\) and keep the other coordinates constant in the completed
\'etale chart. The extension \(\R((s))/\R((t))\) has degree \(m_i\)
since \(t=s^{m_i}U(s)\). The image point has this exact degree
because its residue field contains the algebraic coordinate
\(z_i=s\), which has valuation \(1\) in \(\R((s))\).
Choose distinct supports and fill any remaining even degree by
closed points with residue field \(\C((t))\). To obtain these,
start with the restriction of the point defining \(\nu\). In
the smooth positive-dimensional generic fiber, \'etale
coordinates and Hensel's lemma give nearby \(\C((t))\)-points
avoiding any finite set of supports and their conjugates. If the
starting point is real, take the chart over \(\R((t))\) and vary
one coordinate by \(it^M\) for large \(M\), giving a point that
does not descend to \(\R((t))\). If it is nonreal, sufficiently
small perturbations remain nonreal. Each such conjugate pair
thus gives one reduced closed point of degree two.
It follows that edge sets describe precisely the end parities
of degree \(N\) effective local configurations when their weight
at each mark is at most \(N\) and congruent to \(N\) modulo two.

Compatibility means that the edge set has zero boundary modulo
two at every vertex. The parity of the weights over a mark
equals the parity of the crossing edges because the sign of the
base parameter changes on a transverse arc exactly for odd
multiplicity. Zero boundary makes this parity the same at
successive marks. Even degree allows the empty edge sets. For
odd degree, apply Corollary~\ref{cor:closed-point-interpolation}
to the selected and filled local configurations. Finally an
even-incidence edge set decomposes into edge-disjoint simple
cycles. If its winding is odd, at least one cycle has odd
winding. Deleting the others can only reduce the local weight
sums, so does not worsen the minimum.
\end{proof}

\subsection{Generation in the Chow group of one-cycles}\label{sec:real-chow-one}

Let \(G_1(X)\subset\operatorname{CH}_1(X)\) be the subgroup
generated by integral real curves whose complex irreducible
components have rational normalization. A conjugate pair of
complex rational curves is allowed.

\begin{corollary}\label{cor:real-chow-one}
For every smooth projective geometrically rationally connected
variety \(X/\R\), one has \(G_1(X)=\operatorname{CH}_1(X)\).
\end{corollary}

This gives a Chow-group version of the real cycle-class generation
problem highlighted by Wittenberg \cite[p.~2351, after Questions~3.2]{WittenbergICM}.

\begin{proof}
By Tian--Zong \cite[Theorem~1.3]{TZ}, the norm from
\(\operatorname{CH}_1(X_\C)\) takes values in \(G_1(X)\).
Restriction followed by the norm is multiplication by two,
so \(2M_X=0\) for \(M_X=\operatorname{CH}_1(X)/G_1(X)\).
A geometrically reducible integral real curve is itself the
norm of one of its complex components.

Apply the correspondence argument of Lemma~\ref{lem:real-families}
with values in \(M_X\), using pushforward to \(X\). Divisibility by
two on real Jacobian identity components and the Nash curve selection
step require only \(2M_X=0\). Thus a pushed one-cycle has constant
class in \(M_X\) along each connected component of the real points of
its projective Hilbert scheme.

Let \(v:B\to X\) normalize a geometrically integral curve. If
\(B(\R)=\varnothing\), Corollary~\ref{cor:points} gives
real maps \(B\to Q\) and \(Q\to X\) for the anisotropic conic
\(Q\). Their composite \(w\) pushes \([B]\) into \(G_1(X)\).
By Theorem~\ref{thm:chow}, the generic points of the two
graphs have equal classes in \(\operatorname{CH}_0(X_{\R(B)})\).
By Chow localization their difference is supported on finitely
many closed fibers of \(B\times X\to B\). All of these
have complex residue field, so their pushforwards are norms
and belong to \(G_1(X)\).

If \(B(\R)\ne\varnothing\), the construction below gives
a map \(w:B'\to X\) with \([w_*B']=0\) in \(M_X\).
After adding conjugate very free teeth, both maps have
acyclic tangent pullbacks of positive Euler characteristic
and are equivariantly homotopic on an identified oriented
differentiable source. The teeth have classes in \(G_1(X)\).
Lemma~\ref{lem:real-chow-vary} then gives \([v_*B]=[w_*B']=0\)
in \(M_X\).
\end{proof}

The degeneration into rational components uses the source-degeneration
approach of Tian and Zong \cite{TZ}. For a source with real points,
the choice of degeneration tracks each real circle. The comparison of
the resulting maps uses equivariant disk replacement with varying source
complex structures. Families of graphs then give the required relations
in \(M_X\).

\paragraph{Rational components of a real nodal source.}

We now assume \(v:B\to X\) is a real map from a smooth
projective geometrically integral real curve with real points.
Set \(\Sigma=B(\C)\), let \(\sigma\) be conjugation, and let
\(Y=X(\C)\) with its conjugation. We first construct a smooth
real curve \(B'\) with a map \(w:B'\to X\), as follows.

The quotient \(\Sigma/\sigma\) is a compact bordered surface
whose boundary circles correspond to \(B(\R)\). Over its interior
the two-fold cover is its orientation cover. Cut off an annular
collar of each boundary circle by a parallel interior curve.
In the residual orientable surface, cut disjoint meridians
that remove its handles. If the residual surface is nonorientable,
cut off disjoint M\"obius bands, one for each crosscap, so that
the residual part is planar. These choices follow from the
classification of compact surfaces with boundary. The chosen
cut circles are all two-sided and lift to pairs of disjoint
conjugate circles. The components obtained upstairs are planar.
After capping the cut ends they are spheres. The two capped
ends of each cut are marked for identification to a node.
The spheres above the original annular collars have one fixed
circle each. The spheres above the M\"obius bands carry the
fixed-point-free real structure. Other spheres occur in conjugate
pairs. For genus zero the source is already a rational generator,
so collars that cut off disks cause no exceptional case.

Realize the marked spheres algebraically by real split projective
lines, anisotropic conics and conjugate pairs of complex projective
lines, according to these types. The marked points belonging
to the cut circles are nonreal and distinct. Identify the
marked pairs to form a projective real nodal curve \(D\).
The nodes occur in conjugate pairs. Its complex normalization
has only rational components. Smoothing all nodes gives back
the equivariant differentiable surface \(\Sigma\): complex
plumbing inserts the annuli removed by the cuts, and there is no
selection of real smoothing signs because the nodes are nonreal.
The identification can track and label the original real circles.
One may take a projective real one-parameter smoothing.
For instance, a sufficiently ample embedding of \(D\) satisfies
\(H^1(D,T_{\mathbf P^m}|_D)=0\) by the Euler sequence. Embedded
deformations then lift the unobstructed abstract nodal
deformations, and a real tangent vector nonzero in each local
node-smoothing direction gives the required curve of parameters.

Choose a point \(x_*\in X(\R)\). On every split component of \(D\),
approximate the corresponding smooth loop from \(v|_{B(\R)}\)
by a real algebraic map, prescribing value \(x_*\) at all marked
node branches. Apply Theorem~\ref{thm:main}
to the constant family. Real-analytic approximation of the loop
gives an equivariant holomorphic germ along the real circle.
Prescribe constant germs on small disjoint disks about the nonreal
marked points, with the conjugate germs on the other disks.
On the remaining components one can start with the constant map
to \(x_*\). Each set of paired node values is therefore compatible.

We use the usual comb deformation for sections at this point.
For a map \(a:S\to X\) from a fixed real smooth curve, view its
graph as a section of \(S\times X\). Attach sufficiently many
very free vertical teeth in general nonreal fibers, along with
the conjugates. With a prescribed divisor \(A\) of marks disjoint
from the nodes of the comb, the graph deformation makes
\(H^1(S_\C,a_{\mathrm{new}}^*T_{X_\C}(-A))=0\),
and can preserve the given values on \(A\).
This is the comb construction of Hassett--Tschinkel
\cite[Lemma~20, Proposition~24 and Lemma~25, revised version]{HT}.
One may embed the graph in an auxiliary projective bundle to
separate teeth. Any required open general-position conditions
can be met with conjugate pairs: their real locus in the product
of a smooth complex parameter chart and its conjugate is
Zariski dense. The deformation scheme with the prescribed
mark conditions is smooth at the comb, and a real tangent
direction can smooth all nodes. Projection to \(S\) is of
degree one on the smoothed curve, and hence an isomorphism.
Projection from an auxiliary projective space preserves the
needed vanishing, since the relative normal bundle of a
smoothed graph in the product splits off the auxiliary factor.
Over a collar of the real locus the teeth are absent, so a
small real smoothing preserves the homotopy classes of the real
loops. The same construction on a pair of conjugate components
is made on one component and conjugated on the other.

Apply this with \(A\) containing each node branch on each
component. The normalization sequence for the nodal map
\(w_D:D\to X\) now gives \(H^1(D_\C,w_D^*T_{X_\C})=0\):
vanishing with the marks removed ensures surjectivity of every
node evaluation in that sequence. Relative Hom into \(X\)
over the chosen source smoothing is smooth at \(w_D\), its
obstruction being this \(H^1\). Take a real algebraic parameter
curve transverse to the central fiber in the relative Hom scheme.
This produces a nearby \(w:B'\to X\). Its source has the
equivariant differentiable type of \(B\), and its real loops
are homotopic to those of \(v\) under the tracked identifications.
Using a relative projective embedding of the source and graph
cycles gives a real Hilbert family through the central and
smoothed maps. The family comparison shows that
\[
                      [w_*B']=0\quad\text{in }M_X ,
\]
since the central curve has only rational geometric components.

\paragraph{Correction of the equivariant homotopy.}

Use an orientation-preserving equivariant diffeomorphism
\(B'(\C)\simeq\Sigma\) as above. First make the maps \(v,w\)
acyclic and arrange \(\chi(v^*T_Y),\chi(w^*T_Y)>0\) by
the same conjugate comb deformation on the fixed curves.
The teeth raise the Chern degree without changing the genus.
The classes in \(M_X\) and the homotopies of the real loops
are unchanged by the family comparison.

There is an equivariant homotopy on a collar of the real locus,
by retracting onto that locus and using the homotopies of the
real loops. Triangulate the complement of a smaller invariant
collar by free equivariant cells. Extend the homotopy over
zero- and one-cell orbits, using connectedness and simple
connectedness of \(Y\). For each remaining two-cell orbit the
obstruction on a chosen lift is a sphere class in \(\pi_2(Y)\).
Changing one endpoint map inside that disk by a based sphere
changes the obstruction by that class. Since \(Y\) is simply
connected, basepoint transport is unique on \(\pi_2(Y)\), and
Hurewicz identifies it with \(H_2(Y,\mathbb Z)\).

By Lemma~\ref{lem:complexcycles} each such class is a signed sum of classes of
rational curves. They can be taken very free using the comb
construction of \cite[Proposition~2.4]{TZ}. Choose distinct
attaching points inside the selected disk, and attach the
positive and negative parts to the appropriate endpoint maps,
and attach the conjugate teeth over the conjugate disk.
To realize one very free summand as a difference of teeth with
specified attaching images, choose a common very free auxiliary
map through both attaching images and one point on the summand,
at distinct marks. On the chosen side, join the auxiliary map
to the summand and smooth while fixing the attaching image.
The normalization sequence gives \(H^1=0\) after a double twist
at that mark, so the smoothing can still be taken very free.
On the other side use the auxiliary map alone. This realizes the
signed correction with the same auxiliary class on both sides.

Smooth the two combs with each handle's source fixed. More
explicitly, blow up the product of that source with \(\mathbf A^1\)
at the attaching points in the central fiber. The central
fiber is the comb. Its map has acyclic tangent pullback by
normalization, since the handle is acyclic and the teeth
evaluate surjectively at their nodes. Relative Hom to \(X\)
is smooth there, so take a real transverse parameter curve
and a nearby smooth fiber. Under the plumbing identification,
the smoothed map is homotopic to the map with the indicated
spheres inserted. Indeed it converges to the handle on the
complement of the disks and to each tooth on its part of the
plumbing. In a sufficiently small neck the image lies in a
contractible chart about the node image, by continuity of the map
on the total source. This identifies the local pinching with
based sphere insertion relative to the outer disk boundary
up to homotopy. The construction is conjugation equivariant.
All cell obstructions have thus been canceled. The two smoothed
maps are equivariantly homotopic on \(\Sigma\), remain acyclic
with positive Euler characteristic, and their classes in \(M_X\)
are those of \(v,w\), respectively. The last assertion follows
from the family comparison and the rationality of the added
teeth.

\paragraph{Disk replacement with varying source.}

\begin{lemma}[Varying source]\label{lem:real-chow-vary}
Let two real maps from smooth projective real curves to \(X\)
have equivariantly identified oriented differentiable sources and
be equivariantly homotopic on that source. Suppose the
pullback tangent bundles for the endpoint maps have \(H^1=0\)
and positive Euler characteristic. Then their pushed-forward
one-cycles have the same class in \(M_X\).
\end{lemma}

\begin{proof}
Use a smooth equivariant homotopy \(h_t\) stationary on endpoint
collars. Equivariant smooth approximation relative to these
collars follows by embedding \(Y\) equivariantly in a real
representation and projecting an averaged approximation from
an invariant tubular neighborhood. Join the source complex
structures \(J_0,J_1\) by \(J_t\) stationary on the same collars.
For example take a convex interpolation of invariant metrics
compatible with the endpoint structures and use the fixed
orientation. Every such almost complex structure on a
surface is integrable, and \(\sigma\) is antiholomorphic.

\emph{Acyclic operators.}
Use Lemma~\ref{lem:acyclic} with the principal symbol replaced
by \(\bar\partial_{J_t}\). The bundles in the dual equation can
be identified smoothly along the interval, and uniform elliptic
estimates apply. A nonzero dual section is nonzero on a nonreal
disk, where arbitrary zero-order variations and their conjugates
give the same surjective universal differential. The index stays
\(-\chi(h_t^*T_Y)\), so the finite-parameter argument gives real
acyclic operators \(D_t\) natural on the endpoint collars. Their
right inverses are uniformly bounded \(L^4\to W^{1,4}\), as in
Lemma~\ref{lem:eq-replace}.

\emph{Uniform local disks.}
Use the equivariant fiber-holomorphic local addition of
Proposition~\ref{prop:analytic}, with its uniform small \(q\)-ball.
For each disk choose a smoothly varying \(D_t\)-holomorphic
frame. In this frame the local addition at its center \(c\)
is a holomorphic coordinate parametrization on a smaller ball.
Use its inverse as the moving target chart throughout the
parameter interval, independently of the sphere families.
The local expression has the expansion
\[
  m(q)+a(q)\zeta+b(q)\overline\zeta+O(|\zeta|^2).
\]
Here \(\zeta\) is centered at \(c\) in a \(J_t\)-holomorphic
chart of uniformly positive size. The bounds also hold after
one \(q\)-derivative on a smaller parameter ball. Such charts
and holomorphic frames are available on finitely many
parameter and source neighborhoods, uniformly by compactness
and smooth elliptic dependence \cite{AB}. Frames can also be
obtained on small charts by solving the matrix
\(\bar\partial\)-equation with a Cauchy--Green operator.
For a fixed disk one can choose charts and frames smoothly
throughout the interval on a common smaller neighborhood.

Use fixed differentiable disks in reference coordinate charts,
in conjugate pairs off the real locus. On a disk of reference
radius \(r\), choose its \(J_t\)-conformal coordinate \(w\)
onto \(\{|w|<r\}\), centered at \(c\). We retain \(\zeta\)
for the larger chart used in the Taylor expansion. In the
normalized coordinate \(z=w/r\), write
\[
           \zeta=r F_{t,r}(z),\qquad
           F_{t,r}(z)=\sum_{\nu\ge1}c_{\nu}(t,r)z^\nu .
\]
The change of variable and its inverse have uniformly bounded
derivatives in scaled coordinates on the closed disks.
In fact the image in the scaled \(\zeta\)-coordinate is a
smooth family of planar domains that extends, as \(r\to0\),
to ellipses in a compact family of nonsingular ellipses.
The smooth parameter dependence and boundary regularity of the
normalized Riemann maps apply \cite[Corollary~9.4]{BG}.
Alternatively, solve the Dirichlet problem with boundary
value \(\log|\zeta/r|\) in the scaled domain. If \(u\) is its
solution, then \((\zeta/r)\exp(-u-i\widetilde u)\) is the
normalized map up to rotation, where \(\widetilde u\) is a
harmonic conjugate. The maximum principle and the argument
principle give a degree-one map, and boundary regularity and
the Hopf lemma give nondegeneracy. Parametric Dirichlet
estimates on the fixed reference disk then give the claimed
uniform bounds, including as \(r\to0\).

\emph{Laurent replacement.}
Fix a small width parameter \(\theta>0\) and put
\(M=\lceil\theta^{-1/2}\rceil\). Define on the punctured disk
\[
 R_{t,r}(w)=
 r\left(\sum_{\nu=1}^M\overline{c_\nu(t,r)}
                     (w/r)^{-\nu}
       +\theta^2(w/r)^{-M-1}\right).
\]
The last term fixes the pole order \(J=M+1\). On the boundary
the Laurent series of \(\bar\zeta\) is exactly the infinite
reflected series. The uniform smooth bounds give
\(|c_\nu|\le C_k(1+\nu)^{-k}\) for each fixed \(k\).
Thus the truncation error on the circle is \(O(\theta r)\),
and the radial derivatives of the polynomial and of
\(\bar\zeta\) are bounded uniformly on
\(1-O(\theta)\le|w|/r\le1\). Indeed \((1-O(\theta))^{-M-2}\)
is bounded, and \(\sum_\nu\nu|c_\nu|\) is bounded. In
particular \(R_{t,r}-\bar\zeta=O(\theta r)\) on this
collar. The same estimates hold in conjugated coordinates on the
paired disk.

Choose local holomorphic first-jet sphere families \(G\) by
\cite[Proposition~3.5]{DGX}, with
\(G(0;m,b)=m\) and derivative \(b\) in the sphere coordinate
at zero. On a disk core the replacement is
\[
 G\bigl(\lambda R_{t,r}(w);m(q)+a(q)\zeta,b(q)\bigr),
                    \qquad 0<\lambda\le1 .
\]
It is holomorphic in source and fiber variables on the core,
including at the pole using the chart at infinity. At
\(\lambda=0\) use the holomorphic Taylor expression
\(m(q)+a(q)\zeta\) in the moving target chart. The standard
gluing annulus has width comparable to \(\theta r\), with
\(\lambda=1\), and the difference between the two expressions
there is \(O(\theta r+r^2)\), also after one \(q\)-derivative.
Take \(r\le\theta\). The cutoff derivative then gives a
uniform bound on the normalized \((0,1)\)-coefficient and its
first fiber derivative on all standard annuli. The second fiber
derivative is uniformly bounded on a smaller ball by Cauchy's
estimates in \(q\). The fiber
derivative of the interpolated map is invertible on the annuli
for small \(r\), because both expressions there differ by
\(O(r)\) in \(C^1\) in \(q\) from \(m(q)\), which is a
coordinate diffeomorphism on the parameter ball.

For a change of local sphere family, widen just that annulus
to inner radius \(r/2\), collapse \(\lambda\) to zero, and
then increase it with the new family and narrow the annulus.
At zero the two Taylor expressions agree in the moving chart
chosen independently of the sphere family. With \(\theta\) fixed,
the Laurent polynomial and its
derivatives on the widened region have bounds
\(|R|\le C_\theta r\), \(|dR|\le C_\theta\). At intermediate
widths the cutoff derivative is at most \(C/(\theta r)\).
Consequently the normalized coefficient and its fiber
derivative there are bounded by a constant depending on
\(\theta\), uniformly in small \(r\). One can likewise turn
a disk off at \(\lambda=0\) by interpolating the cutoff to
zero on the whole disk. The same bounds follow from
\(\Phi-(m+a\zeta)=O(r)\) there. Choose the radii small enough
that a wide switch on at most one conjugate pair contributes
at most \(\theta^{1/4}\) to the \(L^4\) bounds for these
coefficients. This is possible since its area is \(O(r^2)\)
and, for fixed \(\theta\), the coefficient bound is independent
of \(r\). It also ensures all annular images stay in the
chosen target coordinate ball.

The disks may be selected after setting all these bounds:
take a finite disjoint Vitali packing on the free quotient
covering all but \(O(\theta)\) of the reference area, with
arbitrarily small radii, then lift in pairs. The fixed locus
has area zero. The standard annuli also have total area
\(O(\theta)\) since \(\sum r^2\) is bounded. Compactness of
the jet paths gives finitely many branch intervals on each
disk. Choose distinct switching times in their overlaps
and disjoint intervals so that only one pair is widened
at any time. Turn on and off pairs one at a time within the
stationary endpoint collars. All constants before the choice
of \(\theta\) come from compact chart, jet and operator data.
In particular letting the maximum disk size depend on
\(\theta\) does not change the contraction constants for the
standard annuli.

\emph{Correction and the cycle limits.}
Away from cores let \(A_t(z,q)\) be the \((0,1)\)-coefficient
obtained by applying the inverse fiber derivative of the local
addition to its source \(\bar\partial\)-derivative in a
\(D_t\)-holomorphic frame. On the holomorphic cores it is zero.
For the contraction of Proposition~\ref{prop:analytic}, the coefficient
is continuous in \(t\) in the fixed point norm. The annular expressions
vary continuously up to zero scale and the core coefficient remains
zero. On unmodified parts of stationary collars the natural endpoint
operator makes \(A_t(z,0)=D_qA_t(z,0)=0\). The preceding area and
derivative bounds give size \(C(\theta^{1/4}+\rho^2)\) and Lipschitz
bound \(C(\theta^{1/4}+\rho)\) on the \(C^0\)-ball of radius \(\rho\),
as in Lemma~\ref{lem:eq-replace}. The equivariant contraction and
holomorphic reconstruction therefore give algebraic maps with
\(C^0\)- and \(W^{1,4}\)-continuous corrections. On the endpoint
collars, where the disks are off, the zero correction is the unique
solution.

Use embeddings of the varying sources into a fixed real
\(\mathbf P^l\). For these take a fixed smooth real line bundle of
degree at least \(2g+1\), from a real divisor on one endpoint
and transported over \(\Sigma\). Choose varying real
Dolbeault operators on it. Its \(H^1\) vanishes for each
\(J_t\), and its \(H^0\) has a smoothly varying real basis
over the interval by elliptic projection. The resulting
embeddings vary smoothly. The graphs of the corrected maps
are therefore continuous in a real Hilbert scheme of
\(\mathbf P^l\times X\) between switches. Indeed their degrees are
locally constant and their genus is fixed. A Hilbert limit
on such an interval has the limiting smooth graph as its
top cycle. The fixed polynomial precludes any extra
zero-dimensional part, so it is that graph.

At a collapse with pole order \(J\), the rescaling
\(w=\lambda^{1/J}z\) gives one bubble parametrized by
\(G(Az^{-J};m(q(c)),b(q(c)))\), with \(A\ne0\) the top
Laurent coefficient at the limit. This is the rational
bubble with pushforward multiplicity \(J\). To check the
absence of other top-cycle contributions, the correction
is holomorphic on the core in a uniformly varying frame.
Cauchy estimates give uniform derivative bounds near its
center. The map derivative on the neck is bounded by
\(C(1+\lambda |w|^{-J-1})\) in the region
\(\lambda^{1/J}/\delta\le|w|\le\delta\) for sufficiently
small \(\delta\), with constants for the fixed packing.
Integration gives neck area \(O(\delta^2+\delta^{2J})\).
The source embedding tends to a point in the rescaled
chart. Thus the Hilbert limits have top cycle the graph
of the limiting map plus the rational bubbles over source
points, by current convergence and the Hilbert--Chow map
\cite[\S3.1(6)]{Ross}. The same holds at a turn-on.
The conjugate bubbles push into \(G_1(X)\). Apply
the family comparison on each branch and to its limiting
real Hilbert points. The classes agree across every switch
in \(M_X\). This proves the lemma.
\end{proof}

\section{Effective Campana tight approximation}\label{sec:campana}

Let \(k=\R\) or \(\C\), and let \(p\colon\cX\to B\) be a flat projective morphism with \(\cX\) smooth and \(B/k\) a smooth projective geometrically integral curve. Let \(\mathcal D_i\) be horizontal effective Cartier divisors, let \(m_i\ge2\) be integers, and put \(\Delta=\sum_i(1-1/m_i)\mathcal D_i\). A local section is Campana if no component of its image is contained in any \(\mathcal D_i\), and every positive order of its pullback of \(\mathcal D_{i,\C}\) is at least \(m_i\). These conditions are imposed at all complex points of the domain.

\begin{theorem}\label{thm:campana}
Suppose that the geometric generic fiber is smooth and integral, that its boundary sum is reduced with simple normal crossings, and that \(-\bigl(K_{\cX_{\bar\eta}}+\Delta_{\bar\eta}\bigr)\) is nef and big. Let \(u\) be a Campana holomorphic section on \(\Omega\), with the real conditions of Theorem~\ref{thm:effective} if \(k=\R\). Fix \(K\subset\Omega\) compact, \(m\ge0\), a very ample \(H\), and finitely many jets of \(u\). There exist \(C,A>0\) such that for every \(0<\varepsilon<1\) there is a \(k\)-algebraic section inducing these jets, satisfying~\eqref{eq:effective-tight}, and satisfying the Campana condition at every point of \(B(\C)\).

In particular, there are Campana algebraic sections with the prescribed jets converging to \(u\) in the compact-open \(C^\infty\) topology on \(\Omega\).
\end{theorem}

\begin{proof}
Write \(\mathcal E=\bigcup_i\mathcal D_i\) for the reduced support. The geometric generic variety is rationally connected by the theorem on klt log Fano pairs \cite{Zhang}. Indeed, write the nef and big class \(L=-(K+\Delta)\) as \(A_0+E_0\), with \(A_0\) ample and \(E_0\) effective. For small \(e>0\), the pair with boundary \(\Delta+eE_0\) is klt and its anticanonical class \((1-e)L+eA_0\) is ample. Relative dimension zero and \(\Omega=B(\C)\) are handled by properness and GAGA as before.

\smallskip
\emph{Vertical curves.}
By Chen--Tseng \cite[Theorem~1.1 and \S5.4]{ChenTseng}, the geometric generic pair has a free rational log curve of interior nef curve class, with contacts of exact order \(m_i\) and avoiding the deeper boundary strata. It meets every geometric boundary component because an interior nef curve class has positive intersection with each nonzero effective divisor. The existence theorem applies over algebraically closed fields of characteristic zero by the Lefschetz principle and finite-type descent. Spread the marked curves over a finite cover of an open of \(B_\C\) where the pair is relatively log smooth. Here freeness gives \(H^1(N_{\mathrm{rel}}(-q))=0\) for the relative logarithmic normal complex and any one contact marking \(q\). In the finite-type relative spaces of deformations of these marked curves the contact evaluations are smooth there \cite[Proposition~2.9 and Definition~4.2]{CLT}. Consequently in each horizontal component there is a dense open of allowed contact values, in good fibers and on just one smooth boundary branch, covered by such vertical curves. The degrees of curves required over these opens are bounded by an integer \(M_H\). In fact, a fixed finite collection of numerical classes suffices, by taking one family of the marked curves for each of the finitely many generic geometric components. The complement of these allowed opens in the boundary is contained in a closed subset of the total space of codimension at least two. In the real case take the opens stable under conjugation and use the curves in conjugate pairs.

\smallskip
\emph{Protecting the local data.}
Choose a compact neighborhood \(K^+\) of \(K\) in \(\Omega\), with
the interpolation points in its interior. Over \(\R\), take \(K^+\)
invariant under conjugation and include the real locus and the
conjugate interpolation points in its interior. Take its boundary
disjoint from the contacts of \(u\). Prescribe jets of lengths greater
than all contact orders at the finitely many contacts in \(K^+\).
Add the original interpolation points and the finitely many base
values outside the good open used in the preceding step. Use the
given germ at the latter values in \(K^+\). At such a value outside
\(K^+\), a complex algebraic section exists by
\cite[Theorem~1.1]{GHS}, so the fiber has a relative smooth point.
No horizontal divisor contains a component of a fiber. Thus one
may choose holomorphic germs through relative smooth points away
from the boundary. Over \(\R\) these additional values are nonreal
and the germs are chosen in conjugate pairs.

If there is no prescribed center off the boundary, choose one where
the given local section avoids the boundary, or choose a new
boundary-avoiding germ on a disk over a good value. In particular,
the latter choice supplies data when \(\Omega\) is empty and there
are no bad values. Over \(\R\) take its conjugate as well. Choose
the neighborhood of \(K^+\) and the auxiliary disks small enough
to leave a nonempty open subset of \(B(\C)\) unused. Denote the
interpolation divisor by \(D=\sum_b L_b[b]\), taking each \(L_b\)
greater than the contact orders at its center and large enough for
the required jets.

Choose a fixed effective divisor \(P\) on \(B\) of degree at least
\(g(B)+2\), supported at a closed point where these holomorphic data
are available. A conjugate pair gives such a closed point over
\(\R\). Use Theorem~\ref{thm:effective} with the interpolation
divisor \(D+P\). This gives an ordinary algebraic approximant \(s\)
of degree \(O(\varepsilon^{-A})\) with
\begin{equation}\label{eq:campana-positive}
 H^1(B_\C,s^*T_p(-D-P))=0.
\end{equation}
For the vanishing, use the linearization in the proof of
Proposition~\ref{prop:fixed-order}, omitting conjugation over \(\C\).
On the quantitative packing of Theorem~\ref{thm:effective}, the defect
area is still \(O(\theta)\) and the parameter solution is
\(O(\theta^{1/4})\) in the sup norm. Thus the same Neumann estimate
applies to the very sections with the claimed degree bound.
Apply the approximation with tolerance \(c\varepsilon\) on \(K^+\)
for a sufficiently small fixed \(c>0\). By Rouch\'e's theorem on fixed
small disks around the protected contacts and uniform closeness away
from them, all contacts on \(K^+\) have the same centers and orders as
those of \(u\).

Make an arbitrarily small algebraic deformation of \(s\) fixing \(D\), so every unprotected boundary contact is transverse and in the allowed open. To justify genericity, \eqref{eq:campana-positive} implies \(H^1(s^*T_p(-D-Q))=0\) for every effective \(Q\) of degree at most two, by multiplying by a nonzero section of \(\cO(P-Q)\). These vanishings persist on a neighborhood by upper semicontinuity over the proper parameter spaces for \(Q\). The fixed-\(D\) section space is smooth there, and its varying value and first-jet evaluations away from \(\Supp(D)\) are smooth by the same infinitesimal evaluation argument. The excluded values have codimension at least two. Non-transverse relative section jets meeting a relatively smooth boundary branch also have codimension two, given by the boundary value and an independent derivative equation. Each corresponding incidence over the one-dimensional base projects into a proper closed subset after taking closure. In the real case the real points in a smooth local chart are Zariski dense in every analytic neighborhood of \(s\), so the perturbation can be chosen real. The degree is locally constant and closeness on \(K^+\) persists. In particular all contacts requiring repair are outside \(K^+\).

\smallskip
\emph{Repairing the short contacts.}
At each \(b\in\Supp(D)\), blow up \(L_b-1\) successive points of the transforms of \(s\), as in the jet blowups of \cite[\S6.2]{CLT}. The calculation here is local at the section. Choose a base parameter \(t\) and relative coordinates \(x\) in which \(s\) is \(x=0\). After \(j\) blowups the chart along the section has \(x=t^j y\). A section through the final point has \(y(0)=0\), so its blowdown gives the required jet. The induced map on the relative normal sheaves in a single blowup is multiplication by \(t\), giving a twist by \(-[b]\). If a Cartier boundary equation \(F(t,x)\) vanishes to order \(\mu\) at the center, its strict transform in the first chart has equation \(t^{-\mu}F(t,ty)\). While the center lies on the strict transform, \(\mu\geq 1\), so the contact order along the section drops by at least one. The original orders are less than \(L_b\). Hence the final marked points are off the strict boundary, also when the original boundary is singular there. Remove the singular locus of the reduced strict boundary from the blown-up model. The remaining log target has smooth, disjoint boundary branches and contains the transformed handle and the teeth used below. We use the local deformation theory of log maps on this quasi-projective target \cite{GrossSiebert}.

Mark the final specified points as ordinary incidence constraints. At each simple boundary contact, the logarithmic normal bundle of the handle equals its ordinary normal bundle: the boundary tangent condition removes one direction both in the pullback of the target tangent bundle and in the source tangent bundle, and the transverse derivative identifies those quotient directions. Thus the normal bundle of the handle with incidence constraints is \(s^*T_p(-D)\), with vanishing \(H^1\). At each remaining contact attach a vertical free log curve of degree at most \(M_H\), as in the first step. Glue its attaching marking of order \(m_i\) and the handle contact of order one by a contracted three-pointed bridge \cite[\S2.6.2]{CLT}. The new marking on the bridge has order \(1+m_i\). Other marked tooth contacts already have the required orders.

For the simultaneous gluing, use the normalization computation of \cite[Proposition~2.10]{CLT}. Each bridge has normal complex \(\cO^{\dim\cX-1}\) plus a skyscraper sheaf away from both nodes, with the constant summand evaluating to the tangent of the boundary at either node. Each tooth has vanishing \(H^1\) of its normal complex after subtracting the attaching point, by freeness. This uses that over a good fiber the quotient of the absolute logarithmic tangent by the relative one is the trivial base tangent on the tooth. In the normalization sequence the constant summand of each bridge supplies the discrepancy at its handle-side node. The independent tooth evaluation supplies the discrepancy at the other node. Thus the node evaluation is surjective even with all the protected handle constraints, and the normal complex of the comb with these constraints has vanishing \(H^1\).

\smallskip
\emph{Smoothing and degree.}
This vanishing gives smoothness over the stack of log structures with the incidence constraints fixed \cite[\S2.4 and Proposition~2.10]{CLT}. The basic characteristic monoid of the comb has one independent generator \(e_j\) for each bridge. In the gluing chart of \cite[\S2.6.2]{CLT} the two node lengths at that bridge are \(m_i e_j\) and \(e_j\), since both noncontracted components have zero degeneracy and meet the same boundary component. There are no relations between distinct bridges. Use the affine line with the chart sending every \(e_j\) to its coordinate \(t\), and pull the deformation stack back to this log line. At \(t\ne0\) the base log structure is trivial and all the nodes smooth.

Over \(\R\), the repair points are nonreal. Choose the bridges and enhancements in conjugate pairs. Conjugation permutes their independent generators and keeps \(t\) fixed, so the chart and central map descend as \'etale log structures. Labels of markings are indexed by a finite \'etale scheme so conjugation can exchange paired markings. The marked deformation stack is the corresponding descent of the ordered one, and smoothness is checked after complexification. There is a smooth scheme chart through the descended real point with the same residue field by \cite[Th\'eor\`eme~6.3]{LMB}. Its projection to the affine line is smooth at that point. The real implicit function theorem gives smoothings for arbitrarily small nonzero real \(t\). Over \(\C\) the same argument uses the complex chart.

The smoothed map has degree one over \(B\), so its smooth connected domain is identified with \(B\). Near \(K\), the limiting map on the handle is a section and no nodes occur there. Local analytic inverses of the base projection in the smoothing family give convergence on a neighborhood of \(K\) to the handle, hence \(C^m(K)\) convergence by Cauchy estimates. Incidence at the marked final points preserves the original jets. Over the unprotected values the model is unchanged, and the nondegenerate log map has only its specified marked contacts with the boundary. Over each protected value the orders downstairs were fixed by the original jets. This gives the global Campana condition.

The number of repairs is at most \(\sum_i\mathcal D_i\cdot s(B)\), which is \(O(\deg(s^*H))\) because \(H^{\otimes c_i}(-\mathcal D_i)\) is ample for fixed \(c_i\) sufficiently large. Degree downstairs in the smoothing is the sum of the handle degree and the degrees of the teeth, each at most \(M_H\). The contracted bridges have degree zero. The smoothed degree is therefore \(O(\varepsilon^{-A})\). Choose the smoothing close enough on \(K\) for the desired error. Compact exhaustion now proves the convergence assertion.
\end{proof}

Chen, Lehmann and Tanimoto proved finite Campana jet interpolation over algebraically closed fields of characteristic zero under strong Campana uniruledness and ordinary weak approximation at the chosen places \cite[Theorem~1.6]{CLT}. Their toric case includes arbitrary fibers supporting the jets \cite[Corollary~1.10]{CLT}. Dandapat obtained all-places Campana weak approximation for the del Pezzo orbifolds in \cite[Theorem~E]{DandapatCampana}. Chen and Tseng proved the smooth weak log Fano case of the Campana rational connectedness conjecture formulated in \cite[Conjecture~1.2]{LTcampana}, with the ample case stated in \cite[Conjecture~5.8]{CLT}. Their proof supplies the exact-contact free curves used here \cite[Theorem~1.1 and \S5.4]{ChenTseng}. Theorem~\ref{thm:campana} adds compact-set approximation, conjugation-equivariant interpolation over \(\R\) and a degree estimate. The contact repair uses the jet blowups and the log gluing by contracted bridges of \cite[Sections~2.6--2.7 and~6.2--6.3]{CLT}. Here the repairs are made outside the approximation compact set, and the simultaneous smoothing is chosen equivariantly for the pairs of conjugate bridges. The degree constants depend on the fixed holomorphic section, compact set and interpolation data.

\subsection{Campana sections and finite covers}\label{sec:campana-finite-etale}

For the data of Theorem~\ref{thm:campana}, put \(K_*=K\) over \(\C\).
Over \(\R\), put \(K_*=K\cup\overline K\cup B(\R)\), where
\(\overline K\) is the conjugate of \(K\). The following application
requires \(B(\C)\setminus K_*\) to contain a nonempty open set.

\begin{corollary}\label{cor:campana-finite-etale}
Let \(B^\circ\subset B\) be a nonempty open subcurve and let
\(U\subset p^{-1}(B^\circ)\setminus\bigcup_i\mathcal D_i\)
be a dense open subscheme. Suppose \(U\to B^\circ\) is smooth with
geometrically connected fibers and \(U(\C)\to B^\circ(\C)\) is a
locally trivial topological fiber bundle. For some \(C,A>0\) and
every \(0<\varepsilon<1\), there is a \(k\)-algebraic Campana section
\(s\) with the prescribed jets and satisfying \eqref{eq:effective-tight}
for these constants. The open \(V=s^{-1}(U)\) is nonempty and
\(\pi_1^{\mathrm{et}}(V)\to\pi_1^{\mathrm{et}}(U)\) is surjective.
For every connected component \(B'\) of a finite \'{e}tale base change
of \(B^\circ\), the map \(V_{B'}\to U_{B'}\) also induces a surjection.
Thus each connected finite \'{e}tale cover of \(U_{B'}\) remains
connected over \(V_{B'}\). If \(B'\to B^\circ\) has degree \(d\), the
pulled-back section has \(H\)-degree over \(k\) at most
\(d C\varepsilon^{-A}\), and degree normalized by \(d\) at most
\(C\varepsilon^{-A}\).
\end{corollary}

\begin{proof}
Choose a small coordinate disk \(D\) in \(B^\circ(\C)\) outside
\(K_*\) and the prescribed jet centers. Over \(\R\), take \(D\)
disjoint from its conjugate and avoid the conjugates of the centers
as well. We may take \(D\) inside a single topological trivialization.
The fundamental group of the fiber over a point of \(D\) is finitely
generated, since a smooth complex algebraic variety has finite CW
type. Choose finitely many based loops generating it. Embed a finite
connected graph \(\Gamma\) in \(D\) with that many independent loops,
or take a point if the group is trivial. The trivialization defines
a continuous section \(v\) on \(\Gamma\) by assigning the chosen
fiber loops to the generators of \(\pi_1(\Gamma)\).

Delete one point from each relatively compact complementary component
of \(D\setminus\Gamma\) and denote the resulting open Riemann surface
by \(D'\). The scalar Mergelyan theorem gives the Mergelyan property
for \(\Gamma\) in \(D'\). The local
Mergelyan theorem for manifold-valued maps approximates \(v\) by
holomorphic maps to \(U(\C)\) on neighborhoods of \(\Gamma\)
\cite[Theorem~1.4]{ForstnericLocalMergelyan}. To make such a map
\(\widetilde v\) a section, fix a Stein neighborhood of the original
graph \(\{(z,v(z)):z\in\Gamma\}\) in \(D'\times U(\C)\), small
enough that the second projection lands in \(p^{-1}(D)(\C)\)
\cite[Theorem~2.2 and the proof of Theorem~1.4]{ForstnericLocalMergelyan}.
On this neighborhood, Cartan's Theorem~B lifts the constant
vector field \(\partial_z\) via
\[
d(p\circ\operatorname{pr}_2):
\ker d\operatorname{pr}_1\longrightarrow
(p\circ\operatorname{pr}_2)^*T_D.
\]
Use its holomorphic flow for time \(z-p(\widetilde v(z))\) at
\((z,\widetilde v(z))\), using the coordinate on \(D\) to take
the difference. For a sufficiently close approximant the flow
is defined near the compact graph and yields a holomorphic
section \(v'\) close to \(v\) on \(\Gamma\). Over \(\R\), prescribe
the conjugate section near the conjugate graph.

Apply Theorem~\ref{thm:campana} on the union of these neighborhoods
and a neighborhood of the original protected data disjoint from
them. The image of each graph lies in \(U(\C)\). Uniformly close
maps of the graph to \(U(\C)\) are homotopic there, with one
path relating their base points. Thus an approximation tolerance
\(\delta\varepsilon\) on \(K\) and the auxiliary graphs, with fixed
\(0<\delta\leq 1\), ensures that the resulting section has
\(\Gamma\subset s^{-1}(U)(\C)\) and carries all the chosen
fiber generators into \(\pi_1(U(\C))\), up to the same change of
base point. The degree bound changes only by a fixed factor.

The homotopy sequence of the bundle identifies the image of the
fiber group with the kernel of
\(\pi_1(U(\C))\to\pi_1(B^\circ(\C))\).
The inclusion \(s^{-1}(U)(\C)\subset B^\circ(\C)\) removes only
finitely many points and surjects on fundamental groups. The
section map therefore surjects on \(\pi_1(U(\C))\).
The Riemann existence theorem gives geometric \'{e}tale
surjectivity, and over \(\R\) the fundamental group exact
sequences with quotient \(\operatorname{Gal}(\C/\R)\) give
arithmetic surjectivity. The maps from the fundamental groups of
\(s^{-1}(U)\) and \(U\) to that of \(B^\circ\) are surjective.
Indeed \(s^{-1}(U)\) is dense open in the smooth base curve, and
its inclusion factors through \(U\) by \(s\). For a connected
finite \'{e}tale base change \(B'\), both pullbacks are thus
connected. Their fundamental groups are the inverse images of the
open subgroup corresponding to \(B'\) in
\(\pi_1^{\mathrm{et}}(B^\circ)\). Restricting the original
epimorphism to these inverse images gives the new epimorphism.
This applies to every connected component of a base change. Let
\(g:B'\to B\) be the finite map of smooth projective completions.
The pullback section \(s'\) and the pulled-back polarization
\(H'\) satisfy \((s')^*H'=g^*s^*H\). The projection formula gives
\(\deg_k((s')^*H')=(\deg g)\deg_k(s^*H)\), as claimed.
\end{proof}

\begin{corollary}\label{cor:campana-nonthin}
Under the hypotheses of Theorem~\ref{thm:campana} and the same
condition on \(K_*\), the Campana sections with the prescribed
jets in any \(C^m(K)\)-approximation neighborhood of \(u\) form a
nonthin subset of the generic fiber over \(k(B)\). Given finitely
many finite extensions \(L/k(B)\), the section may be chosen
outside any prescribed thin subset after base change to each \(L\).
\end{corollary}

\begin{proof}
Write \(F=k(B)\), and include \(L=F\) to test nonthinness over
\(F\). For each of the finitely many extensions \(L/F\),
write the prescribed thin subset using proper closed subsets and
images of integral dominant generically finite covers of degree
greater than one. On the normalization \(B_L\) of \(B\) in \(L\),
spread out the covers and choose a dense open \(V_L\) of the
base-changed model that avoids the closures of the proper closed
subsets and over which the covers are finite \'{e}tale of their
original degrees. For a cover given by a rational map, work on
its domain of definition. If a dense source open was used to
spread out the cover, also remove the closure of the image of
its complement in the original domain. This image has smaller
dimension than the generic target.

Shrink \(B^\circ\) so that \(B_L^\circ\to B^\circ\) is finite
\'{e}tale for every \(L\). For the finite projection
\(\pi_L:\cX_{B_L^\circ}\to\cX_{B^\circ}\), remove
\(\pi_L(\cX_{B_L^\circ}\setminus V_L)\) from
\(\cX_{B^\circ}\setminus\bigcup_i\mathcal D_i\), for every \(L\).
The images are closed and proper on the generic fiber by
finiteness and the dimension bound on each complement. The
result is a dense open \(U\), whose entire pullback to
\(B_L^\circ\) lies in \(V_L\). Each cover restricted to this
pullback is integral, since its source is a nonempty open
in the integral source of the original cover. Since \(U\) has
a smooth geometrically integral generic fiber, its fibers are
smooth and geometrically connected after shrinking the base.
For local topological triviality, compactify and resolve with
simple normal crossings boundary outside \(U\), then shrink so
that the compactification and all boundary strata are smooth over
the base. Smooth horizontal lifts of base tangent fields, tangent
to the boundary strata, give the local trivializations by
properness.

Apply Corollary~\ref{cor:campana-finite-etale} to this \(U\)
and its base changes to the \(B_L^\circ\). The chosen generic
point avoids each proper closed subset. Every restricted cover
pulls back to a connected finite \'{e}tale cover of the nonempty
open \(s_L^{-1}(U_{B_L^\circ})\) in the smooth integral curve
\(B_L^\circ\). Its source is normal and connected, hence
integral. Its generic fiber is therefore the spectrum of a
field of degree greater than one over \(L\) and has no
\(L\)-rational point. This excludes every prescribed thin set.
\end{proof}

The degree constants in Corollary~\ref{cor:campana-finite-etale}
depend on the fixed open model through the auxiliary graph data.
The same constants bound the normalized degree on every finite
\'{e}tale base change.

\Needspace{11\baselineskip}
\section{Approximation over function fields of real surfaces}
\label{sec:real-surface-approx}

Let \(S/\R\) be a smooth projective geometrically integral surface with
\(S(\R)=\varnothing\), and let \(U\subset S\) have finite complement.
Consider a flat projective morphism \(f\colon X\to U\) with smooth total
space, geometrically integral fibers, and a relatively very ample line
bundle \(L\). We assume that the smooth geometric generic fiber with
this polarization satisfies Hypothesis~7.8 of \cite{DHS} and contains a very
twisting scroll whose rulings have \(L\)-degree one. These geometric
conditions hold over a nonempty open \(V\subset U\). Let \(\mathfrak D\)
be the reduced divisorial discriminant. There is a finite subset
\(T\subset S\), containing \(S\setminus U\), such that \(\mathfrak D\) is
smooth away from \(T\) and
\begin{equation}\label{eq:real-surface-critical-image}
df(T_xX)=T_y\mathfrak D
\quad\text{if }y=f(x)\notin T \text{ and } f \text{ is not smooth at }x.
\end{equation}
To obtain \(T\), stratify the reduced critical locus over each smooth
discriminant component. A stratum dominating the component is smooth
over a dense open of it by generic smoothness in characteristic zero.
The remaining strata have finite images. After removing those images,
the tangent image contains the tangent line of the component and is
one-dimensional at every critical point.

\begin{theorem}\label{thm:real-divisorial-approx}
Under the above assumptions, \(X_{\R(S)}\) has an \(\R(S)\)-point.
Assume in addition that the relative dimension of \(f\) is at least
three. Let \(D\) be a smooth geometrically integral real divisor of
\(S\), contained in \(S\setminus T\), meeting \(V\), and transverse
to \(\mathfrak D\).
Put \(K=\R(S)\) and
\(K_D=\operatorname{Frac}\widehat{\cO}_{S,\eta_D}\).
Then \(X_K(K)\) is dense in \(X_K(K_D)\) for the valuation topology.
\end{theorem}

The complex surface existence theorem is \cite[Corollary~13.2]{DHS}.
We use their canonical components of stable sections and ruled-surface
deformations. The real step uses Theorem~\ref{thm:main} and the
fixed-jet tangent vanishing in Proposition~\ref{prop:fixed-order}
both to prepare real families of sections and to interpolate in the
fixed Abel fiber. Divisorial interpolation will require a singular
member of the pencil.

\subsection{Real paths at a fixed Abel class}
\label{subsec:real-surface-paths}

Consider a real pencil of curves on a modification of \(S\). Denote
its base by \(B\), a Zariski open of \(\Pone_\R\) containing its real
locus, and its fibers by \(C_t\). In this subsection all the parameters
under consideration have smooth geometrically integral \(C_t\),
with \(X_t=X\times_U C_t\) smooth and \(C_t\) meeting \(V\).
Write \(g=g(C_t)\). We give the connecting argument on intervals
with these properties. We use the peaceful locus over the open of
\(C_t\) lying over \(V\), so the hypotheses on line evaluation in
\cite[Lemma~7.11]{DHS} apply there.

At a real parameter \(t\), the curve \(C_t\) is point-free. Theorem~\ref{thm:main}
and the tangent vanishing in Proposition~\ref{prop:fixed-order} give
a real section \(s\) of \(X_t\to C_t\), meeting the peaceful locus of
\cite[Definition~7.4]{DHS}, with \(H^1(C_t,s^*T_{X_t/C_t}(-P))=0\)
for a fixed effective divisor \(P\) of arbitrarily large degree. By
Riemann--Roch this implies the same vanishing with \(P\) replaced by
every effective divisor of degree at most \(2g+2\). The section scheme
is smooth at \(s\) over \(B\), so there are smooth real local branches
with the same properties. One can also start with a specified local
branch having these properties.

On an overlap, use the product of the two smooth relative section
charts as a parameter space for a pair of endpoints. Pull back the
adapted-chain open with smooth endpoint evaluation of
\cite[Proposition~7.12]{DHS} to the universal source curve over
this parameter space above \(V\). Compactify by taking its closure in the
projective incidence scheme of ordered \(L\)-lines and their marked
meeting points. The pulled-back open is smooth over the smooth
universal curve, so take a projective resolution of the dominating
closure that is an isomorphism on this open. By generic smoothness
the resolved slices over a dense open of the endpoint parameters
are smooth. Over a dense open of the universal curve the projection
from the resolution has smooth rationally connected fibers, since a
general fiber is a smooth projective model of the evaluation fiber
in Proposition~7.12. The image of this open in the endpoint
parameter space is open because the universal curve is smooth
over that space. After shrinking the dense parameter open, each
slice therefore has rationally connected geometric generic fiber
over its source curve. A real pair can be chosen in this parameter
open arbitrarily near the initial branches, as the smooth real
endpoint charts are Zariski dense and their value maps at general
source points are submersions.

Choose a nearby real pair of sections in this parameter open,
over a real parameter \(t\) in the overlap. Work for the moment on
its fixed source curve \(C_t\), allowing both sections to vary in
smooth local charts inside the parameter open. Choose an effective reduced real divisor \(P_0\)
of degree \(2g+2\) in the open where both endpoints are peaceful
and the resolved family contains the open of adapted chains.
Its points form conjugate pairs. Let \(\mathcal T_t\) denote
the product of the two section charts over this fixed curve.
Let \(\mathcal E\) be the Weil restriction from \(P_0\) to
\(\R\) of the space of endpoint pairs, and let \(\mathcal G\)
be the Weil restriction of the adapted-chain open at \(P_0\).
The fiber product \(\mathcal T_t\times_{\mathcal E}\mathcal G\)
records the varying sections and choices of chains at these
points. It has a real point over the chosen pair, since choices
at a point with residue field \(\C\) determine the choices at
its conjugate. These chains lift to the resolved family on the
preserved adapted-chain open.

The evaluation \(\mathcal T_t\to\mathcal E\) is smooth near
the chosen pair by the tangent vanishing after subtracting
\(P_0\). Thus the projection of the fiber product to
\(\mathcal G\) is smooth there. For each intermediate vertex,
the simultaneous vertex evaluation from \(\mathcal G\)
is smooth by \cite[Lemma~7.11(3)]{DHS} and smoothness
under Weil restriction along a finite \'{e}tale morphism.
For each degree \cite[Proposition~4.15]{DHS} excludes a proper
closed subset of the space of values on this fixed divisor.
Its inverse image under the composite vertex evaluation is
nowhere dense in a small real smooth neighborhood in the
fiber product. Indeed, the values on the conjugate pairs form
a real form of the product over all complex roots, and a
nonempty real open in this smooth value space is Zariski
dense. Choose one point outside the countable union of the
excluded sets, keeping the handle pair as close as desired
to its initial value. Now fix this pair. The chosen chains
can serve as interpolated values in Theorem~\ref{thm:main}
on the resolved chain family over \(C_t\).
Every intermediate section has weak \(P_0\)-freeness and
hence is free by \cite[Lemma~4.12]{DHS}.

The projective incidence supplies a ruled surface for each
adjacent pair of sections, also at the boundary of the chain
open. By Proposition~\ref{prop:fixed-order} the lifted chain
section can be chosen with vanishing \(H^1\) of its pulled-back
relative tangent bundle. The section scheme of the resolved
chain family is therefore smooth over the endpoint parameter
space at this section. Here the resolved slices and the source
curves are smooth over that parameter space, and on the slice
the projection to \(C_t\) is smooth along the section.
The chain, its ruled surfaces, and its free handles extend
on a smaller real interval. At each point of this interval
the two resulting endpoint handles can be joined to the
original handles by short paths in their respective smooth
section charts. Shrink the charts so these paths consist of
free handles meeting the peaceful locus. Their \(L\)-degrees
are constant. A finite covering of the compact parameter
intervals thus uses finitely many handle and ruled-surface
families of bounded degrees.

Choose a degree \(e\) large enough for the bounds in
\cite[Lemma~9.5]{DHS} for the finitely many ruled-surface
families used on the overlaps. Also impose the bounds of
\cite[Theorem~13.1(3)]{DHS} on the geometric generic fiber
for the section components of the degrees in use. There are
finitely many such components: sections of a fixed \(L\)-degree
have a fixed Hilbert polynomial after twisting \(L\) by a
sufficiently ample line bundle from the source curve.
By \cite[Lemma~8.4 and Corollary~9.8]{DHS}, one bound
for each free section component suffices. Let \(Z_e\) be the
closure of the canonical component over the pencil. It descends
to the generic pencil field by asymptotic uniqueness, as in
\cite[Corollary~13.2]{DHS}. Choose a line bundle \(A\) on the
pencil surface, of sufficiently large degree \(e\) on the fibers,
whose geometric generic class is in the open of the Abel map
where the fiber on \(Z_e\) is geometrically integral,
rationally connected, and meets the genuine-section locus.
The constructions of \(A\) below allow arbitrarily large \(e\).
Let \(W_A\) be the closure of that generic fixed-Abel fiber
in the projective coarse space of stable sections over the pencil.
Take a projective resolution of \(W_A\) that is an isomorphism
over its regular locus.

We use a local observation to identify its regular charts. Let
\(M\) be an open chart of stable sections on which the Abel class
defines a morphism to the relative Picard space of the pencil.
Suppose the pencil near the parameter has geometrically reduced
and connected fibers. Its structure sheaf has universal direct
image equal to the base structure sheaf by
\cite[Lemma~36.32.6]{StacksProject}. The relative Picard diagonal
is an immersion \cite[Lemma~99.11.9]{StacksProject}, so the
equalizer \(F_A\) of the Abel morphism and the class of \(A\) is
closed in \(M\) after shrinking around the point under study.
Suppose \(F_A\) is smooth over the pencil there. Shrink further
so \(F_A\) is integral, using regularity of the base and of this
smooth chart. The intersection \(W_A\cap M\) is contained in
\(F_A\) scheme-theoretically since its defining generic fiber
lies in \(F_A\). If the geometric generic fiber of \(F_A\)
meets an open of the ambient stable-section space where \(Z_e\)
is the unique component, \(W_A\) contains a nonempty open of
that geometric generic fiber. A generic point of any geometric
generic component meeting this open maps to the generic point
of \(F_A\) by flatness of the field extension. Thus the closed
subscheme \(W_A\cap M\) contains the generic point of the
integral \(F_A\), so \(W_A\cap M=F_A\) in this chart.

For each chosen free handle \(s\) at a real parameter \(t\),
the real line bundle \(A_t\otimes (s^*L)^{-1}\) has even degree
since \(C_t(\R)=\varnothing\). In sufficiently large degree a
section over \(\R\) has a reduced effective zero divisor \(E\)
supported at peaceful values. Indeed, vanishing at a forbidden
pair of points imposes two independent real linear conditions by
Riemann--Roch. A double zero and its conjugate impose four, and
the location of the pair varies in real dimension two. Thus the
excluded set in the real linear system has real codimension at
least two. Its complement is path connected, by perturbing
paths relative to their endpoints.

Attach free lines at \(E\), in conjugate pairs. This stable section
has Abel class \(A_t\) by \cite[Lemma~6.7]{DHS}. It has no
automorphisms: the handle projects isomorphically to \(C_t\),
the attachment points are distinct, and each tooth has \(L\)-degree
one and is embedded. Freeness on the handle and the teeth gives
vanishing of \(H^1\) of the pulled-back relative tangent bundle over
the curve \(C_t\), by normalization and evaluation at the nodes.
The lifting argument for the free combs of
\cite[Lemmas~8.4--8.5]{DHS} then gives smoothness of the
stable-section space over the pencil at that point. Equisingular
deformations moving \(E\) lift
by the smooth evaluation of free lines. Their Abel differential
is surjective by \(H^1(C_t,\cO(E))=0\), using the tangent
sequence of the divisor \(E\). Hence the comb is a smooth point
of the fixed-Abel space over \(B\). The locus of these combs
with their roots in the peaceful open is itself smooth over \(B\)
at the indicated point: locally it is obtained from the free handle
space, the relative space of reduced effective divisors, and the
smooth free-line evaluations by imposing the same Abel equation.
This smooth comb locus meets the geometric generic fiber. At such
a comb \(Z_e\) is the unique component of the ambient stable-section
space by the choice of \(e\) and \cite[Theorem~13.1(3)]{DHS}.
The local observation identifies \(W_A\) with the smooth fixed-Abel
chart near the comb over \(t\). In particular the comb lifts
uniquely to the regular chart of the resolution. This argument
also gives local smooth real branches through these combs over
parameters where other points of the fiber may be singular.

For a fixed handle, the path-connected open of divisors just
constructed connects the root divisors of any two such combs.
Label the pairs along an interval and choose a lift of each pair
to \(C_t(\C)\) there. The family of lines through its image under
the handle is the pullback of a smooth proper evaluation over
the peaceful locus. Its complex fibers are connected by
\cite[Lemma~7.11]{DHS}. The path in root divisors therefore
lifts to paths of lines, with any prescribed endpoint choices,
by local triviality and by inserting paths in the fibers.
Use conjugate choices over the other branches. This gives a
real path in the smooth fixed-Abel comb locus.

The same locus also lifts the short paths of \emph{varying}
handles used to reach the ruled-surface endpoints. Fix \(t\)
and a free-handle chart of the indicated degree. These paths
require only the fixed-Abel comb chart, so the auxiliary open
on which the chain family was constructed need not contain
them. Over the
chart, the choices of reduced root divisors and of free marked
lines at their peaceful values are smooth. Their map to the
relative Picard space is submersive with the handle fixed,
by the root-divisor \(H^1\) vanishing above. Consequently
imposing the Abel class \(A_t\) still gives a smooth morphism
from this comb locus to the handle chart. Every handle on the
short path has a real lift by the large-degree construction,
and the real fiber of this morphism over that handle is
path connected by the fixed-handle argument. Subdivide the
handle path into intervals with local smooth lifts. At a
subdivision point join the adjacent lifts by a path in the
fiber and keep the handle constant during this join. This
gives a path of combs over a reparametrization of the handle
path, with prescribed lifts at its two endpoints if desired.
All these combs lie in the regular chart of \(W_A\), so the
path also lifts to its resolution.

For successive handles penned by the same ruled surface, choose
root divisors in the common open of \(C_t\) where both handles have
peaceful values. This open is nonempty because each handle meets
the peaceful locus. Each ruling through such a root is free by
\cite[Definition~7.4]{DHS}. The preceding fixed-handle paths move
the combs to ones with their teeth in the ruled surface.
The two real comb divisors are linearly equivalent by
\cite[Lemma~9.3]{DHS}. The complex linear equivalence descends,
since restriction of invertible sheaves from the real ruled
surface to its complexification is injective by Hilbert 90.
For a nodal divisor \(D\) in their linear
system, its Abel class is the ratio of the determinants of the derived
pushforwards of \(L|_D\) and \(\cO_D\) to the source curve.
The restriction sequences on the ruled surface show that this ratio
depends only on the divisor line bundle. Thus the Abel class is
constant throughout the nodal locus. Write \(\Lambda\) for the complete
linear system on the ruled surface. As in \cite[Lemma~9.5]{DHS},
the nonnodal locus has complex codimension at least two for
our choice of \(e\). For each endpoint, the closed cone of lines
from it meeting that locus is consequently a proper subset of
\(\Lambda\). The inverse image of the regular locus of \(W_A\)
on the nodal locus of \(\Lambda\) is nonempty. More explicitly,
the two endpoint combs are smooth in the full fixed-Abel space,
and near either endpoint the nodal-divisor family lies in its
unique local component \(W_A\). The image of the whole nodal
open therefore lies in that component, by irreducibility and
closedness in the fixed-Abel space. Choose a real common point
in this regular inverse image, outside both cones. Such points
exist by Zariski density of the real points of \(\Lambda\).
The two joining real projective lines are entirely nodal,
and map to \(W_A\) with fixed Abel class. Each lifts from its
regular inverse image to a morphism on the whole line by
properness and normality of the line. The two lifts agree at
the chosen common point because it is in the regular locus.
Their real points give the required paths on the resolution.

After the handle and ruled-surface families have been fixed on
overlaps, the real transition parameters can be taken outside
the finitely many values where the projective resolution over
\(B\) is not smooth. For each handle interval, take local
lifts to the smooth comb locus, using the relative smoothness
proved above. On a compact interval a finite number suffices,
and two such lifts of the same handle can be joined at a
transition by the root-divisor path. For different handles,
use the lifted short paths to reach the chosen endpoint combs,
then use the ruled-surface paths. Near each chosen transition,
trivialize the smooth proper map on the
real locus, view the two local branches in a single fiber,
and interpolate along the chosen path in that fiber with
constant collars to splice the branches smoothly. Take the
transition intervals disjoint. This gives a smooth real
section, and any prescribed central collar is left fixed.

For the existence assertion of Theorem~\ref{thm:real-divisorial-approx},
take a sufficiently ample real pencil. Start with a smooth
real member \(C\) avoiding \(T\) and transverse to the
discriminant branches. The other member will be chosen through
the two conjugate pairs specified below. Choose the system
ample enough to separate first jets at two additional points
after these four point conditions. Away from \(C\) the ratio
of the moving equation to the equation of \(C\) then has
independent first jets at two distinct points. For two
smooth strata of dimensions \(d_1,d_2\), criticality at the
two points and equality of their values impose \(d_1+d_2+1\)
conditions, whereas the points vary in dimension \(d_1+d_2\).
This gives distinct critical values on the surface, and
distinct critical values on the discriminant branches also
across different conjugate branches. The fixed member
\(C\) is already transverse, and the generators can be
chosen transverse at their base points, which avoid the
discriminant. The finite set \(T\) can be assigned nonreal
pencil values by the same separation of values at conjugate
points. Since the surface has no real points, a real critical
value would agree with the value at a distinct conjugate
critical point. Thus \eqref{eq:real-surface-critical-image}
gives smooth \(X_t\) over the real parameters. The ample
class can also be taken of larger degree than the divisorial
components of \(S\setminus V\), so each smooth connected
member meets \(V\).

To choose \(A\), choose two conjugate pairs on \(C\) so that their
difference generates a Zariski-dense cyclic subgroup of
\(\operatorname{Pic}^0(C_\C)\). Indeed, the locus
\(\{(p,\bar p):p\in C(\C)\}\) is Zariski dense in
\(C_\C\times C_\C\), even locally on any open disk in
the \(p\)-coordinate, by the identity principle on a totally
real submanifold. Fix the second conjugate pair. The Abel image
of the first pair minus the second cannot be contained in a
proper algebraic subgroup, since the complex image of varying
one point with the other fixed is a translate of the Abel
image of \(C_\C\), which generates the Jacobian.
The proper algebraic subgroups of an abelian variety form a
countable set. Their inverse images among conjugate pairs are
closed with empty interior, so Baire gives the choice, with all
four points distinct. Require the pairs to be transverse pencil
base points. The corresponding difference \(\Delta\) of horizontal
divisors generates a Zariski-dense cyclic subgroup of the geometric
generic Jacobian as well. Otherwise a proper algebraic subgroup
containing all multiples descends to a finite extension of the
generic pencil field. Take a trait over the parameter of \(C\) in
the normalization in that extension. The family of Jacobians is
an abelian scheme there, and the scheme-theoretic closure of the
subgroup is flat over the trait, since its structure sheaf is
torsion free over the discrete valuation ring. Its special fiber
has the dimension of the proper generic subgroup and contains
every multiple of \(\Delta|_C\), a contradiction. Positive
powers of a surface line bundle give arbitrarily large degrees,
and tensoring with the powers of \(\Delta\) gives the required
generic Abel class in each such degree.

The preceding real paths give a smooth real section on the
parameter circle. Compactify its resolution over the parameter
line and resolve away from the chosen smooth open. Use
\cite[Proposition~2.3]{BW} to obtain an equivariant holomorphic
section near the real locus, then apply Theorem~\ref{thm:main}.
Prescribe also a value in the genuine-section open at a nonreal
point and its conjugate. Its generic fiber is rationally
connected by the choice of \(A\). The resulting algebraic
section represents a genuine graph over the generic pencil
field and yields a rational section of \(f\).

\subsection{A nodal member carrying divisorial data}
\label{subsec:real-surface-nodal}

Fix a point in \(X_K(K_D)\) and an integer \(m\geq 1\).
Properness gives a section over
\(\widehat{\cO}_{S,\eta_D}/\mathfrak m_D^m\). Its reduction
extends to \(D\). By \eqref{eq:real-surface-critical-image},
\(X_D\) is smooth. The reduction is a section of \(X_D\to D\),
so its graph is in the relative smooth locus, also for \(f\).
By \cite[Lemma~11.1]{DHS}, the smooth locus of each geometric
fiber over \(D\) is strongly rationally connected.

The prescribed truncated section spreads to an open of \(mD\).
At an omitted closed point \(p\), choose parameters \(z,v\) in
\(\cO_{S,p}\) with \(v=0\) defining \(D\), and choose relative
\'{e}tale coordinates near the image of the reduction at \(p\).
Write each coordinate of the jet modulo \(v^m\) as
\(\sum_{j=0}^{m-1}g_jv^j\), with \(g_0\in\cO_{S,p}\)
and \(g_j\in\cO_{S,p}[1/z]\) for \(j>0\).
Such coefficients can be chosen successively modulo \(v\),
since \((\cO_{S,p}/(v))[1/z]=\kappa(D)\) and the reduction is
regular at \(p\). In the chart near the strict transform after
\(a\) successive blowups, \(v=z^a w\). Choose \(a\) so that
all coefficients \(z^{aj}g_j\) for \(j>0\) are regular.
The resulting coordinates lift the reduced section uniquely
by formal \'{e}taleness, now over the algebraic thickening
defined by \(w^m=0\). The local lifts agree with the given
section on the overlaps. This follows from agreement at the
generic point, since the successive layers of the thickening
are invertible sheaves on the reduced curve and contain no
sections supported at isolated points. Perform the blowups
at all omitted points, in conjugate pairs. Denote the strict
transform by \(D'\), and let \(F\) be the reduced union of
\(D'\) and the exceptional chains. This produces a section
over \(mD'\) without changing the valuation at \(\eta_D\).

Add a sufficiently ample smooth curve \(C'\) meeting \(F\)
transversely at smooth points. Require it to be transverse
to the reduced discriminant of the pulled-back family and
to avoid its finite exceptional set and \(S\setminus U\).
Then \(X_{C'}\) is smooth. We take the intersection of \(C'\)
with every component of \(F\) to be large. Now \(C_0=F+C'\)
is a connected nodal curve in a sufficiently ample system. We
take a pencil with this central member and the following
properties.
\begin{enumerate}
\item Every other real member \(C_t\) is smooth and geometrically
integral, \(X_t\) is smooth, and \(C_t\) meets \(V\).
\item The base points avoid the nodes and the finite bad sets.
There are horizontal exceptional divisors allowing independent
changes of the real multidegrees on \(C_0\), by even integers
on the geometrically integral real components, and by equal
integers on exchanged components.
\item There are horizontal divisor classes \(H,\Delta\) with
positive and zero degree, respectively, supported on base-point
sections meeting \(C'\setminus F\) in the central member.
The generic Jacobian class of \(\Delta\) generates a
Zariski-dense cyclic subgroup.
\end{enumerate}
For the first property, over a point \(p\) blown up on \(D\),
the differential of the first blowdown has image the
direction represented by the point of the exceptional
curve. Only the discriminant direction can give singular
points of the pulled-back total space. Over the general
point of each subsequent exceptional curve that image is
\(T_pD\), transverse to the discriminant by assumption.
There are therefore only finitely many exceptional points
to avoid in noncentral members. The exceptional curves
themselves are discriminant components when the original
fiber is singular. At their base points the pencil is
transverse to them for \(t\ne 0\), because the two pencil
equations form independent local coordinates. For the other
points take \(C_\infty\) smooth and transverse to the
discriminant branches. On the complement of \(F\cup C_\infty\),
the parameter is the product of the equations of \(F\) and
the varying \(C'\), divided by the equation of \(C_\infty\).
Multiplication by the fixed invertible factor preserves the
first-jet freedom of \(C'\), even after imposing the four
chosen base points. The same two-point count gives distinct
critical values for noncentral members, including across
conjugate discriminant branches. Points of \(F\) that are
not base points belong only to the central member. The
points of \(T\) outside \(C_0\) can again be assigned nonreal
values. Thus the real noncentral members are regular.
Ampleness before the pencil blowup gives geometric
connectedness, and large degree ensures that each such
smooth member meets \(V\).

For the last two properties, first take the other smooth
member in the ample system, and choose two conjugate
pairs on it away from \(F\) by the cyclic-density
argument above. Take \(C'\) generally through them. The
additional base points on all components come in the
required conjugate orbits. Their horizontal exceptional
divisors give the multidegrees, while the two selected
pairs give \(H\) and \(\Delta\). In what follows we also
use \(F,C',C_0\) for the strict transforms. Pulling back
the prescribed data through the pencil blowup preserves
the specified thickening of \(D'\) and its generic jet.

Extend the section along each exceptional chain starting at
\(D'\), then extend across \(C'\). For a new component \(G\),
first prescribe its reduced section along the length-\(m\)
divisor at every junction with a previously treated branch.
On an exceptional component the target for this interpolation
is the smooth locus of the original geometric fiber at the
blowup center. Use \cite[Theorem~15]{HTbad} on this constant
strongly rationally connected family, which has a constant
section. On \(C'\) use Theorem~\ref{thm:main} with the
fixed-jet tangent vanishing and prescribe also one value
in the peaceful locus. Since \(X_{C'}\) is smooth, the
differential of a section places its image in the relative
smooth locus. In either case one can arrange that
\[
 H^1(G,h_G^*T_{X/U}(-jG-P))=0\qquad (0\leq j<m),
\]
where \(h_G\) is the reduced section and
\(P=m\sum_x[x]\) is the divisor of prescribed junctions.
For an exceptional component use the very free comb
construction with the jets fixed, in the smooth locus.
For \(C'\) use Proposition~\ref{prop:fixed-order} with a
larger interpolation divisor. In both cases Riemann--Roch
lets one choose that extra vanishing for a divisor \(A\)
with \(\cO_G(A-jG-P)\) effective for all the indicated \(j\).
The inclusion of the corresponding twists then gives the
displayed vanishing.

At a junction write \(u=0\) for \(G\) and \(v=0\) for the
previous branch. In passing from \(jG\) to \((j+1)G\), the
prescribed intersection changes from \((u^j,v^m)\) to
\((u^{j+1},v^m)\). Formal smoothness gives local lifts
matching this intersection. Their differences are sections
of \(h_G^*T_{X/U}\otimes\cO_G(-jG)\) vanishing on \(P\),
because the square-zero layer is \(u^j\) modulo \(u^{j+1}\)
and its restriction at each junction is taken modulo \(v^m\).
The displayed vanishing glues the local lifts at every step.
Finally \((u^m)\cap(v^m)=(u^mv^m)\) gives the fiber-product
description of \(mC_0\) at a node. The sections on the
thickened components thus glue to \(h_m\) over \(mC_0\),
with reduction \(h\). All choices on the exceptional
components are made in conjugate pairs.

\subsection{A fixed Abel deformation retaining the interpolation}
\label{subsec:real-surface-fixed-nodal}

We first build a local branch of bounded degree over the
central parameter. Attach embedded very free curves to
the graph of \(h\) over distinct smooth points of the
components of \(C_0\), in conjugate pairs. The available
tangent lines are dense at a general attachment point
by smoothness of the evaluation of first jets of a very
free curve. Such curves can be chosen embedded in relative
dimension at least three by sufficiently free deformation,
using two-point and tangent evaluation. On \(D'\) and
\(C'\) the attachment points may lie over the smooth
good open, and on an exceptional component use the
strongly rationally connected smooth fiber locus.

We record the normal calculation. On each component of
the handle the normal sheaf of the embedded comb is
obtained from \(h^*T_{X/U}\) by positive elementary
modifications in the tangent directions of the teeth.
At a node of the handle itself the graph is a regular
embedding in a smooth morphism over \(uv=0\), so its
normal sheaf there is \(h^*T_{X/U}\). At the attachment of
a tooth, in the smooth part of the central target, the
node-smoothing map is the fiber quotient of the positive
elementary modification. These identifications follow
from the transverse node equation and are the normal
calculation used in \cite[Section~2]{HT}. Given a fixed
divisor of large multidegree away from the attachments,
the dual of the first cohomology on each modified handle
component consists of sections annihilating the chosen
tangent lines. At a point where a nonzero dual section
does not vanish, a choice of tangent line eliminates
it. Iterate in conjugate pairs, using two disjoint
sets of teeth, each sufficient for the vanishing.
This preserves the vanishing even after one conjugate
pair of teeth is omitted. The resulting modified
normal bundles have vanishing first cohomology even with
the large divisor subtracted. Its degrees can be chosen
to give evaluation surjectivity at the nodes of \(C_0\)
and at any additional specified one or two points, by
Riemann--Roch. The normalization sequence gives the
vanishing on the whole handle. On each tooth the
unmodified normal bundle is an extension of the free
fiber normal bundle by a trivial line, and the node
modification preserves generation and \(H^1=0\).
The normalization sequence then gives the vanishing
on the entire comb, also with the large divisor subtracted.

To see that the vanishing supplies all required
smoothings, fix one tooth node. The handle normal bundle
can attain a nonzero value in its node-smoothing
quotient while matching the handle nodes, by the
evaluation surjectivity just proved. Lift this section to
the whole comb using generation on the teeth. Thus each
tooth smoothing coordinate is nonzero on some tangent
vector in the Hilbert space. A general vector has all
coordinates nonzero. This calculation applies in the
relative Hilbert space over the pencil. The graph at a
handle node remains a regular embedding when the base
equation deforms to \(uv=t\), and all the tooth
attachments are over its relative smooth locus.
For the relative Hilbert lifting, at a handle node the
completed ambient ring over the pencil has equation \(uv=t\)
with the smooth relative coordinates \(z_1,\ldots,z_r\).
The graph at \(t=0\) is cut out by those relative coordinates.
At a tooth node the base curve of the central fiber is smooth
and the comb is cut out in suitable local coordinates by
\(u z_1,z_2,\ldots,z_r\), with \(u\) a parameter on that curve.
Both are regular embeddings in the central fiber. Lifting
the defining regular sequences provides local flat embedded
lifts over a square-zero extension of the parameter scheme.
On intersections, differences of the lifts are sections of
the normal sheaf with coefficients in the square-zero ideal.
The vanishing of \(H^1\) therefore gives smoothness of the
relative Hilbert scheme over the pencil at the comb.
The closed-immersion locus of stable maps is identified
with the corresponding open Hilbert locus by taking the
image. In particular, this chart carries the actual universal
curve, including its restrictions to nonreduced bases. In a smoothing
the domain maps to a smooth \(C_t\) with degree one,
hence is a genuine section. The twisted normal
vanishing gives \((2g+2)\)-freeness by semicontinuity,
taking the original large divisor to have total degree
greater than \(4g+2\) and extending it off the nodes.

To impose an Abel class, add a second group of very
free teeth of one common \(L\)-degree \(a>0\), with
attachment divisor of fixed, sufficiently large
multidegree \(\mathbf b\). On each component such teeth
are available over a dense open. For a common degree
use a multiple of the degrees available at the
geometric generic points and deform the resulting
coverings to embeddings with the root fixed. Spread
by smooth evaluation. Choose an actual line bundle
\(A_*\) on the pencil surface with the multidegree
of the handle, the first teeth, and \(a\mathbf b\)
combined, using the horizontal divisors.

The group of genuine real line bundles of
multidegree zero on \(C_0\) is connected. After
translation by a fixed high-multidegree line bundle,
every element is represented by a real effective
Cartier divisor away from the nodes. This follows
from Riemann--Roch and evaluation at the nodes on the
normalization. The spaces of such divisors are
products of symmetric powers, with conjugate points
paired on real components. They are connected.
Their Abel images cover the group, and evaluation in
large degree shows that the genuine group contains
a neighborhood of its identity in the real Picard
group. Multiplication by \(a\) has invertible
differential and hence is onto the connected group.
Choose the second attachment divisor so its weighted
class supplies the difference to \(A_*\). High
multidegrees allow reduced support at the permitted
points. Varying those support points is Abel
submersive by \(H^1(C_0,\cO(E))=0\) for their divisor
\(E\) and by the invertible differential of
multiplication by \(a\). These variations are
equisingular. They can correct the Abel differential
of any chosen node-smoothing tangent vector without
changing its node coordinates. The extra positive
normal modifications preserve the previous
vanishing.

On a smooth Hilbert chart \(M\) over the pencil, write
\(q:\mathcal C\to S_M\) for the projection of the universal
stable domain to the pulled-back pencil surface and
\(\mu\) for its map to the pulled-back target. Use the Abel
line
\[
  \mathcal A=\det Rq_*\mu^*L\otimes(\det Rq_*\cO_{\mathcal C})^{-1}.
\]
The scheme \(S_M\) is regular because it is smooth over the
regular pencil surface. Properness gives bounded coherent
pushforwards, hence perfect complexes on \(S_M\).
This formula commutes with any change \(M'\to M\).
Both \(\mathcal C\) and \(S_M\) are flat over \(M\).
In fact, for affine rings \(R\to A\to C\) in this diagram,
the tensor product \(C\otimes_A^{\mathbf L}
(A\otimes_R R')\) is \(C\otimes_R R'\) in degree zero,
by the flatness of \(A\) and \(C\) over \(R\).
Thus the square for \(q\) is Tor-independent. Derived base
change and compatibility of determinants with pullback apply,
including over the Artin parameter.

When \(q\) is a genuine section's domain projection,
\(\mathcal A\) is its pullback of \(L\). On a comb, the
gluing exact sequence for \(\mu^*L\) gives the line on
the handle together with one summand for each tooth, twisted
down by its attaching point. A degree-\(b\) tooth has
\(b\) sections of that twist and no \(H^1\). Its pushforward
is a rank-\(b\) bundle supported on the root divisor \(E\)
in the target curve. The effective Cartier divisor resolution
gives determinant \(\cO(bE)\). For \(\cO_{\mathcal C}\) the
corresponding tooth contribution vanishes, since
\(\cO_{\Pone}(-1)\) has no cohomology. These computations hold
for families of combs with finite \'{e}tale root divisors, after
an \'{e}tale splitting and descent. They require no normalization
of the handle, which can itself be nodal. On the locus of
smooth source curves the construction agrees with the Abel
map of \cite[Lemma~6.7]{DHS}.

The relative Picard space of the pencil is smooth near the
indicated class. The obstruction to extending an invertible
sheaf over a square-zero extension lies in \(H^2\), which
vanishes on the curve fibers. Since the chart and the Picard
space are smooth over the pencil, the surjective relative Abel
differential computed from the moving attachments makes the
Abel morphism smooth at the comb. Pulling it back along the
class of \(A_*\) gives a smooth fixed-\(A_*\) chart over the
pencil at that point.

Write \(S_{\rm pen}\) for the pencil surface. Put
\(B_j=\Spec(\cO_{B,0}/\mathfrak m_0^j)\) and
\(C_{B_j}=S_{\rm pen}\times_B B_j\), so \(C_{B_m}=mC_0\).
Keep the handle over
\(B_m\) equal to \(h_m\), and lift the first group of root
divisors in the smooth locus of the pencil. Smooth pointed
evaluation lifts the first group of teeth with their roots
on the prescribed handle. The first group and the handle
determine a line bundle on \(C_{B_m}\) by the comb formula above.
To reach \(A_*\), lift the root line bundle of the second
group to have \(a\)-th power equal to the required difference.
Over \(B_{j+1}\to B_j\), the ideal on the curve is square-zero.
The sequence of units gives existence of a lift by \(H^2=0\),
and the kernel on line-bundle classes is the additive \(H^1\)
of this ideal. Indeed, global scalar units lift, since the
proper flat curve has geometrically connected reduced fibers
and \(H^0(\cO_{C_{B_j}})=H^0(\cO_{B_j})\) by base change.
Multiplication by \(a\) on the kernel is invertible in
characteristic zero. Correct the chosen lift at each step
so that its \(a\)-th power has the required class.
The vanishing of \(H^1\) for the central root line bundle
lifts its section through the same successive steps.
Its simple zeros at smooth points give a relative effective
Cartier divisor, finite \'{e}tale over \(B_m\).
Lift the second group of teeth by smooth pointed evaluation
along that divisor. All lifts are made over the residue fields
of the roots and descended in conjugate pairs.
Gluing along the marked sections gives a flat family of
combs. Its map is a closed immersion over \(B_m\), since the
central map is one and the base is an Artin thickening.
The Abel formula gives \(A_*\) on all of \(C_{B_m}\), yielding a
\(B_m\)-point of the fixed-\(A_*\) Hilbert chart with the
specified handle. To choose an algebraic branch with this jet, take an
\'{e}tale chart of the smooth fixed-\(A_*\) locus over
\(\mathbb A^r_B\) at the central real point. Lift its coordinates
on the truncated base to polynomials in a local parameter, and
pull back the chart along the resulting section of
\(\mathbb A^r_B\). The pointed \'{e}tale base change over \(B\)
recovers the given jet by formal \'{e}taleness. Choose the lift to
meet the open of simultaneous smoothings on its branch. Indeed,
at a sufficiently small fixed nonzero real parameter, bounded
variations of the coefficients of \(t^m\) give an open ball of
coordinates in the real local chart, and the dense smoothing open
meets that ball. The corresponding point belongs to the pointed
branch, since the \'{e}tale chart is a real local diffeomorphism
and these bounded perturbations stay in its neighborhood on the
interval to that parameter. By the twisted normal vanishing, the
geometric generic handle is a free genuine section.
Shrink the real interval so that this holds at its noncentral
points. The branch has fixed degree \(e_0\) and is chosen before
the large degree \(e\) used for the global Abel paths.

\subsection{Completing the approximation}

Use the specified branch near the central parameter,
and construct the ordinary branches and their finite
transition intervals of \S\ref{subsec:real-surface-paths} away from that
parameter. Fix the bridge bounds and apply
\cite[Theorem~13.1(3)]{DHS} to the free handle defined by the
chosen low-degree branch over the geometric generic point.
Choose \(l\) large and then \(k\) so that
\[
 A=A_*\otimes H^{\otimes l}\otimes\Delta^{\otimes k}
\]
has degree above all bounds, including the bound for that
geometric generic handle, and belongs to the generic good Abel
open. Both requirements can be met by the cyclic-density property.

Near the central member the extra free lines can be
added without changing the specified branch. Write
\(Q=H^{\otimes l}\otimes\Delta^{\otimes k}\).
Its presentation by horizontal divisors trivializes
it on \(F\). Since \(F+C'\) is a fiber, locally
over the parameter
\[
 Q(-F)|_F=\cO_F(C'),\qquad
 Q(-F)|_{C'}=Q|_{C'}(-F\cap C').
\]
Our choice of \(C'\) gives vanishing of \(H^1\)
on \(F\) in the first formula, and large \(l\)
gives vanishing and node evaluation on \(C'\)
in the second, independently of \(k\).
Thus \(H^1(C_0,Q(-F)|_{C_0})=0\).
The restriction sequence for \(F\) now gives a
local section of \(Q\) equal to \(1\) on \(F\).
Vary it by sections of \(Q(-F)\), using
cohomology and base change. Sections on \(C'\)
vanishing at its junctions with \(F\) extend
by zero to \(F\) as sections of \(Q(-F)|_{C_0}\).
They allow simple zeros at admissible points of \(C'\), away
from the nodes and the old teeth. The zero divisor is finite
\'{e}tale over the pointed parameter after shrinking. Near a
central zero the projection of the low-degree stable domain is
\'{e}tale at its unique preimage, since source and target are
smooth over the parameter there and the fiberwise differential is
an isomorphism. Properness then gives a finite \'{e}tale map of
degree one over a neighborhood of the zero. Pull the zeros back
to the algebraic low-degree branch and attach free \(L\)-lines.
Smooth evaluation on the space of marked free lines lifts these
attachments over a further pointed \'{e}tale neighborhood, with
conjugate choices over conjugate zeros. The projection is
unchanged on the old handle. In particular the extra attachments
are disjoint from the prescribed Artin section over the generic
point of \(D'\), and the original jet is retained. The Abel
contribution of the new teeth is the zero divisor of \(Q\),
so the enlarged family has fixed Abel class \(A\).

At the enlarged central comb, the normal vanishing above still
holds by the positive modifications of the handle normal bundles
and evaluation surjectivity on the new teeth. The equisingular
directions moving the old attachment divisor lift through the
smooth evaluations of the new teeth. They preserve the node
coordinates and still surject onto the Abel tangent space.
Thus the fixed-\(A\) chart at the enlarged point is smooth over
the pencil. At the geometric generic point of the pointed
algebraic branch, this enlargement is an ordinary comb whose
handle is the chosen free section and whose teeth are precisely
the new free lines. It belongs to \(Z_e\) by the degree bound
in \cite[Theorem~13.1(3)]{DHS}, and to its fiber over \(A\) by the
Abel calculation. In the geometric generic fixed-\(A\) chart
at this comb, \(Z_e\) is the unique component of the ambient
degree-\(e\) stable-section space. This pointed branch supplies
the geometric generic point required by the local observation
above. The central pencil fiber is geometrically reduced and
connected, and the other nearby fibers are smooth and
geometrically integral. Applying the observation in the open
Hilbert chart of embedded maps identifies \(W_A\) there with
the smooth fixed-\(A\) chart as schemes, so the prescribed
\(B_m\)-point lies in \(W_A\).
Choose the projective resolution to preserve this smooth central
chart. At sufficiently close noncentral parameters the enlarged
branch is an ordinary comb on a free handle with its new roots
in the peaceful open, by choosing the zeros on the corresponding
open of \(C'\). The paths in \S\ref{subsec:real-surface-paths}
then join it to the branches on the noncentral intervals, leaving
its prescribed jet at the central parameter fixed. Use the algebraization at the end of that
subsection, prescribing a value in the genuine-section open at a
nonreal pair. For the holomorphic preparation with the specified
parameter jet, use the argument in the proof of
Corollary~\ref{cor:closed-point-interpolation}.

Matching the \(B_m\)-jet takes place in the regular
Hilbert chart, where the resolution is an isomorphism.
The two \(B_m\)-points therefore determine the same embedded
family of stable sections. Over the generic point of \(D'\)
this family is exactly the graph of \(h_m\), because all
old and new attachments are over closed points of \(C_0\).
The projection of the final stable family is an isomorphism
near that point: on the central fiber it is an isomorphism
there, and properness gives a finite map after shrinking.
The locus where this finite map is an isomorphism contains
the specified neighborhood of the central fiber by flatness
over the parameter and Nakayama's lemma.
Thus the graph of the resulting genuine section over the
generic pencil field extends over the generic point of \(D'\)
and agrees with \(h_m\) on its length-\(m\) thickening.
The parameter \(t\) is a uniformizer along \(D'\), since
\(D'\) has multiplicity one in \(C_0\). The section has the
specified truncation modulo \(\mathfrak m_D^m\), and letting
\(m\) vary gives the density in
Theorem~\ref{thm:real-divisorial-approx}.

\subsection{Quartic hypersurfaces in seventeen variables}

The following geometric lemma verifies the generic-fiber hypotheses
of Theorem~\ref{thm:real-divisorial-approx} for smooth quartics.
DeLand developed twisting-surface constructions for general complete
intersections in this degree range.
In \cite[Theorem~1.5, arXiv version]{DeLand} the construction applies
to smooth complete intersections whose Fano scheme of lines is smooth.
He discusses the problem of spanning the tangent space of a
line-evaluation fiber by the positive directions of curves in that fiber
\cite[Remark~7.11, arXiv version]{DeLand}. Here first jets of degree-two
curves in that fiber give the spanning directions needed for every
smooth quartic in \(\mathbf P^{16}\).

\begin{lemma}\label{lem:every-smooth-quartic-scroll}
Over an algebraically closed field of characteristic zero,
every smooth quartic \(Y\subset\mathbf P^{16}\) satisfies
Hypothesis~7.8 of \cite{DHS} and has a very twisting
scroll ruled by free lines.
\end{lemma}

\begin{proof}
All lines through a general point are free, by the
generic freeness of a dominating family of rational
curves, as in \cite[Lemma~7.9]{DHS}. The evaluation
fiber is a smooth complete intersection of type
\((1,2,3,4)\) in \(\mathbf P^{15}\). For general
endpoints \(x,y\) the connecting line is not contained
in \(Y\). The middle point of a two-line chain satisfies
equations of degrees \((1,2,3,1,2,3,4)\) in
\(\mathbf P^{16}\), and cannot be either endpoint.
Both constituent lines are free. The smoothness of
the free chain space and generic smoothness of the
endpoint evaluation give a smooth fiber, of the
expected dimension. The defining equations give a
nonempty fiber by the dimension theorem, and the
sum of their degrees is \(16\). Thus this fiber and
the one-point line fiber are smooth connected Fano
complete intersections. This is the incidence
calculation of \cite[\S17]{DHS}, and it gives
Hypothesis~7.8.

Write \(W\) for the smooth space of pointed free
lines, \(e\colon W\to Y\) for evaluation, and
\(\pi\) for the forgetful map. The fiber \(V\) of \(e\)
over a general point is a smooth complete intersection
of type \((2,3,4)\) in \(\mathbf P^{14}\), of dimension
eleven and index six. Since \(2+3+4=(11+1)/2+3\),
Bonavero and H\"oring's Corollary~1.2
\cite{BonaveroHoeringConics} gives a smooth conic \(Q\subset V\)
that is a quasi-line. Its tangent pullback is
\(\cO(2)\oplus\cO(1)^{10}\), by the tangent-normal sequence.
The two-point evaluation on its deformations is smooth there,
since \(H^1(Q,T_V|_Q(-2))=0\). The same vanishing gives
smooth evaluation of a first jet at a marked point.
Consequently, for a general point of \(V\), the tangent
directions of these conics contain a dense open.

Along a fixed free line in \(Y\), moving the mark
gives a curve \(h\) in \(W\). The line's normal
bundle is \(\cO(1)^{11}\oplus\cO^3\), so
\(h^*T_e=\cO(-1)^{11}\) and \(h^*T_\pi=\cO(2)\).
The first identification follows by taking sections
of the normal bundle that vanish at the moving
point, and the second from the marked-line fiber
of \(\pi\). For a very free conic \(v\) in the
evaluation fiber \(V\), one has
\(v^*T_e=\cO(2)\oplus\cO(1)^{10}\) and
\(v^*T_\pi=\cO(-2)\), using the tautological
line of directions through a fixed point for the
last equality. Both pullbacks of \(T_W\) are
globally generated. Fix a general point \(q\) on
the conic and attach \(h\) at a different point
of \(v\), choosing the corresponding fixed line.
The pair smooths with \(q\) fixed since the
normalization sequence gives \(H^1\) vanishing
for \(T_W\) with that mark subtracted. Let \(E\)
denote the pullback of \(T_e\). At another mark
\(p\) on \(v\), disjoint from \(q\) and the node,
\(H^1(E(-p))=0\) on the pair. Its single
section is along the conic tangent, vanishing
at the node and \(p\) and nonzero at \(q\).
After smoothing, \(E\) has splitting
\(\cO(1)\oplus\cO^{10}\), and the positive
direction at \(q\) specializes to that tangent.
The degree of the forgetful relative tangent
line on the smoothed curve is zero.

Choose enough such curves with common marked
image \(q\) so their positive directions span
\(T_{e,q}\). Attach them to a constant
\(\Pone\) at \(q\). A section of \(E^\vee\)
on this tree would restrict on the constant
handle to a covector annihilating all the
positive directions. It must vanish on the
handle and then on each tooth. The absolute
tangent bundle has \(H^1=0\) on the tree,
also after subtracting a mark on the constant
handle, by normalization. Hence we can smooth
the tree to a free map \(g\colon\Pone\to W\).
The resulting \(g^*T_e\) is ample by semicontinuity
of \(H^0(E^\vee)\), and \(\deg g^*T_\pi=0\).
The composition with \(e\) is free, since \(e\)
is smooth along the image. Apply Lemma~12.6 of
\cite{DHS} with \(m=2\) and then Definition~12.7
to obtain the very twisting scroll. Its rulings are lines in
\(Y\), since \(g\) parametrizes pointed lines.
\end{proof}

\begin{corollary}\label{cor:real-surface-quartics}
Let \(S/\R\) be a smooth projective geometrically integral
surface with \(S(\R)=\varnothing\). A quartic in seventeen
variables over \(\R(S)\) with smooth geometric projective
hypersurface has a nonzero zero if it has a flat quartic
hypersurface model whose total space is smooth off the
inverse image of a finite set of \(S\).
\end{corollary}

\begin{proof}
Over a geometric point where the total space is smooth,
the fiber is integral. Otherwise its singular locus has
a component of dimension at least two, since the fiber
is a reducible or nonreduced degree-four hypersurface
in \(\mathbf P^{16}\). The two derivatives in the base
directions are quartic sections on this singular locus.
They have a common zero there by projective dimension,
which would be a singular point of the total space.
The hyperplane bundle gives the polarization. Now use
Lemma~\ref{lem:every-smooth-quartic-scroll} and
Theorem~\ref{thm:real-divisorial-approx}.
\end{proof}

The model hypothesis in Corollary~\ref{cor:real-surface-quartics}
restricts the allowed divisorial degenerations. The arbitrary
coefficient case of Lang's \(C_2\) question also includes
geometrically singular generic hypersurfaces.

\section{Further applications}\label{sec:further-apps}

\subsection{Degrees of curves on real quartic threefolds}

Benoist and Wittenberg asked whether the minimal degree and the minimal geometric genus of a geometrically integral curve of even geometric genus could be bounded uniformly on smooth real quartic threefolds without real points \cite[Example~7.3]{BWhodgeII}. By adjunction such a threefold is Fano, so Corollary~\ref{cor:conic} answers the genus question with the bound zero. We study the degree question below.

The deformation and local boundedness of anisotropic-conic maps were already studied by Allums \cite[Theorem~2.7.2 and its proof]{Allums}. Here smoothness of the relative space of maps at an unobstructed map from the conic gives the compact-subset assertion. The explicit smooth-branch double-quadric family in Proposition~\ref{prop:quartic-double-quadric} reaches a different part of the boundary of the smooth quartic locus than the ordinary-double-point families of \cite[\S2.7]{Allums}.

\begin{proposition}\label{prop:quartic-compact}
Let \(K\) be a Euclidean compact subset of the parameter space of smooth quartic threefolds in \(\mathbf P^4_\R\) without real points. There is an integer \(D_K\) such that each quartic in \(K\) contains a geometrically integral curve of geometric genus zero and degree at most \(D_K\).
\end{proposition}

\begin{proof}
For a fixed quartic \(Y\) in this locus, apply the construction on the anisotropic conic \(Q\) to \(Y\times Q\to Q\), with no interpolation. One may prescribe a germ on a small closed disk and its conjugate, using a constant map into \(Y_\C\) on one disk and the conjugate constant on the other. Take the packing parameter small enough to apply the linearization in the proof of Proposition~\ref{prop:fixed-order}. The resulting real morphism \(f:Q\to Y\) then satisfies \(H^1(Q_\C,f_\C^*T_{Y_\C})=0\). It is nonconstant since \(Y(\R)=\varnothing\).

For the universal family over the smooth quartic parameter scheme, this vanishing makes the relative scheme of morphisms from \(Q\) smooth at \([f]\) over the parameter scheme. Its real points map onto a Euclidean neighborhood of \(Y\), with the degree \(d_f=\deg(f^*\cO_Y(1))\) locally constant. Restrict to quartics without real points. The reduced image of each nearby morphism is geometrically integral of degree at most \(d_f\), and its normalization is geometrically rational by Riemann--Hurwitz. A finite subcover of \(K\) and the maximum of the associated degrees give \(D_K\).
\end{proof}

The global bound for even geometric genus is equivalent to the existence of a geometrically integral curve of even geometric genus on every smooth quartic threefold without rational points over every real closed field. To see this, let \(R\) be real closed and \(C\) a geometrically integral curve on such a threefold. In the normalization sequence for \(C\), the quotient sheaf has finite support and its Euler characteristic is even, since each supporting closed point has residue field \(R(\sqrt{-1})\). Thus the parity of \(\chi(\cO_C)\) equals that of \(1-g(\widetilde C)\), where \(\widetilde C\) is the normalization. For fixed degree \(d\), the arithmetic genus has only finitely many possibilities. Over the algebraic closure a general plane projection gives a finite birational morphism of \(C\) onto a plane curve of degree \(d\), so the arithmetic genus is between zero and \((d-1)(d-2)/2\). For a degree bound \(D\), the existence condition is therefore first-order over real closed fields, using the geometrically integral loci of finitely many relative Hilbert schemes with odd Euler characteristic. A bound over \(\R\) transfers to every real closed field. Conversely, if no bound over \(\R\) exists, the theory of real closed fields together with the existence of a smooth quartic without rational points and with no such curve of degree at most \(D\), for each integer \(D\), is finitely satisfiable. Compactness gives a real closed field with a quartic that has no such curve at all. This is the bounded-degree transfer principle in this setting \cite[Proposition~9.1 and \S9.3]{BWhodgeII}.

\Needspace{11\baselineskip}
\begin{proposition}\label{prop:quartic-double-quadric}
Let \(q\) be a positive definite real quadratic form in five variables and let \(g\) be a real quartic form. Suppose \(\{q=g=0\}\) is a smooth complete intersection over \(\C\). There exist Euclidean neighborhoods of \(q\) and \(g\), a number \(a>0\), and an integer \(D\) such that, for any forms \(q',g'\) in these neighborhoods and every real \(s\) with \(0<|s|<a\), the quartic
\[
  \{(q')^2+s^2g'=0\}\subset\mathbf P^4_\R
\]
is smooth, has no real points, and contains a geometrically integral curve of geometric genus zero and degree at most \(D\).
\end{proposition}

\begin{proof}
With \(y\) of weight two and the five \(x_i\) of weight one, consider the projective family in \(\mathbf P(1,1,1,1,1,2)\) defined by
\[
  q'(x)=s y,\qquad y^2+g'(x)=0 .
\]
It avoids the point where all \(x_i\) vanish. The charts with \(x_i\ne 0\) are smooth, and projection to the \(x_i\) defines a morphism to \(\mathbf P^4\). At \((q',g',s)=(q,g,0)\) the fiber \(Y_0\) is a double cover of the smooth quadric \(\{q=0\}\), branched along the smooth divisor \(\{q=g=0\}\). The cover is smooth and geometrically integral. Its canonical bundle is the pullback of \(\cO_{\{q=0\}}(-1)\), since the branch divisor is twice the class of \(\cO_{\{q=0\}}(2)\). Thus \(Y_0\) is Fano and has no real points. The relative Jacobian criterion on the weight-one charts and properness give a smooth projective family over an open neighborhood of this parameter point.

For \(Y_0\), the unobstructed conic map and the relative Hom argument
in the proof of Proposition~\ref{prop:quartic-compact} give conic maps
to all sufficiently nearby real fibers. Their degree \(D\) for
the pullback of \(\cO_{\mathbf P^4}(1)\) by projection is locally constant.
For \(s\ne0\), elimination of \(y\) identifies the fiber with the
displayed quartic. On the real unit sphere the forms \(q'\) are
uniformly positive and the \(g'\) uniformly bounded after shrinking
the neighborhoods. The quartics thus have no real points for all
sufficiently small real \(s\). The images of the conic maps have
geometric genus zero and degree at most \(D\), as in
Proposition~\ref{prop:quartic-compact}.
\end{proof}

The following result extends the degree bound to isolated
singular limits. Its use of disk replacement on the very free
locus requires keeping the crossings with its divisorial complement
in protected disks.

\begin{proposition}\label{prop:pointfree-isolated-klt}
Let \(Y\) be a geometrically normal real Fano hypersurface
of dimension at least three. Suppose \(Y_\C\) is klt and
has only isolated singularities.
In a Euclidean neighborhood of its equation, the
smooth point-free hypersurfaces of the same degree contain
geometrically integral curves of even geometric genus and
uniformly bounded degree. If \(Y_\C\) is \(\mathbf Q\)-factorial,
one can take maps from the anisotropic conic and hence
curves of geometric genus zero.
\end{proposition}

\begin{proof}
Write \(U=Y_{\rm reg}\), and let \(V_\C\subset U_\C\) be
the locus swept by the images of very free morphisms from
\(\Pone_\C\) to \(U_\C\). The smooth locus \(U_\C\) is rationally
connected by \cite[Corollary~1.2]{XuFreeSingular}. A component of
the space of morphisms \(\Pone_\C\to U_\C\) therefore has
dominant evaluation at \(0\) and \(\infty\). By generic
smoothness, on the smooth locus of the reduced component its
differential is surjective at a general morphism \(g\). The
differential is the restriction map from \(H^0(\Pone_\C,g^*T_U)\)
to the two fibers, restricted to the component tangent space.
A splitting summand of degree at most zero cannot surject
onto its two fibers, so all summands of \(g^*T_U\) have
positive degree. This gives a very free morphism. Its
evaluation, and that of every very free morphism, is smooth.
Hence \(V_\C\) is nonempty, open, and invariant under conjugation. The image of every such map lies entirely in \(V_\C\).
Let \(D_1,\ldots,D_r\) be the divisorial components of
\(U_\C\setminus V_\C\). Include their singular and pairwise
intersection loci in a closed subset \(Z\) of \(U_\C\) of
codimension at least two, together with the other components
of \(U_\C\setminus V_\C\). Take this subset conjugation
invariant.

The manifold \(U(\C)\) is simply connected.
Compare a smooth hypersurface of the same degree with
\(Y_\C\) by smoothing and removing pairwise disjoint Milnor
balls at the isolated singular points. One may choose a smoothing
whose added coefficient is nonzero at all these points.
Locally the smoothing is then the Milnor fibration of the
defining function divided by this unit.
Outside the balls the smooth and singular fibers have
diffeomorphic exteriors by the proper submersion theorem
with boundary. Since the complex dimension is at least
three, the links and Milnor fibers are simply connected
\cite[Theorems~5.2 and 6.5]{MilnorSing}.
By van Kampen the exterior has the fundamental group of
the smooth hypersurface, which is trivial by Lefschetz.
The punctured conical neighborhoods at the singularities
give the same fundamental group for \(U(\C)\).

\emph{A smooth source with positive crossings.}
There is a smooth equivariant map from the sphere of
the anisotropic conic into \(U(\C)\), with two distinct
conjugate values in \(V\) on a pair of disks.
Indeed, a path from a point to its conjugate and its
conjugate path form the boundary loop of one hemisphere.
Fill it using simple connectivity, extend to the other
hemisphere by conjugation, and smooth equivariantly,
preserving the constant germs on the disks.
Make the map transverse to the \(D_i\setminus Z\) and
disjoint from \(Z\) by a small relative equivariant
perturbation. Transversality applies on the quotient
of the free source, using conjugate choices of local
perturbations. The crossings are finite. Conjugate
crossings have the same sign since conjugation reverses
the orientation of both the source and the normal
complex line of the divisor.

For each conjugate pair of negative crossings, take a
fresh pair of conjugate smooth complex complete
intersection curves in \(Y_\C\). Choose general
hyperplanes so these curves lie in \(U_\C\), avoid \(Z\),
and meet the \(D_i\setminus Z\) transversely.
Each hyperplane intersection with the closure of
the indicated divisor has positive degree, so the
curve has a positive crossing on that same divisor.
Choose one such crossing and the conjugate one on the
mate. Connect the negative crossing and the chosen
positive crossing by a path in the smooth irreducible
divisor away from \(Z\) and the other divisors.
In small normal coordinates around the endpoints,
flatten the tangential part of each map and
remove a disk whose boundary maps to a normal circle.
Join these boundaries by an annulus mapping to the
normal circle bundle along the path, inside \(V(\C)\).
The opposite crossing signs make the boundary
identifications orientation reversing with the
boundary orientations of the two source complements:
the negative disk complement maps with degree \(+1\)
to the normal circle, and the positive complement
maps with degree \(-1\). Trivialize the normal line
along the path and keep the path constant in small
end collars to smooth the attached map.
Attach the conjugate tube at the conjugate crossing.
This operation removes both negative crossings and
introduces only the remaining positive crossings on
the added curves.

After all such operations the connected source is a
closed oriented surface of even genus. Each conjugate
pair of complete intersection curves adds twice their
genus, and the involution still acts freely and reverses
orientation. At each positive crossing, its normal
coordinate gives an oriented coordinate on the source
disk. Flatten the tangential part using that coordinate.
Choose a conformal structure for which those coordinates
are holomorphic and the involution antiholomorphic.
For instance, average a metric with its conjugate
while keeping the prescribed conformal classes on
the conjugate disks. The resulting real Riemann surface
comes from a smooth projective real curve \(C\).
The map is now holomorphic near all crossings and the
two originally protected disks, and its remaining
image lies in \(V\).

\emph{Replacement relative to the crossings.}
Apply the construction of Lemma~\ref{lem:positive} and
Proposition~\ref{prop:analytic} to the product with \(U\).
In these arguments projectivity is used for the source and the
eventual algebraization in \(Y\). Along the compact graph, the
relative charts are in the smooth open \(U\); the two preparation
disks and every replacement center only need very free maps
through their images.
The very free spheres inserted in the bundle preparation
are inserted in a conjugate pair of disks whose images
lie in \(V\), away from the protected disks.
The only centers needed in the replacement
are on the compact complement of smaller crossing disks
where the original map is already holomorphic.
Their images lie in a compact subset of \(V\).
For the local first-jet families at those centers, precompose
one of these very free maps by a ramified self-map of
\(\Pone\) to get zero derivative at the marked point.
Length-two jet evaluation is smooth by
\(H^1(\Pone,g^*T_V(-2))=0\). This provides holomorphic
families for small derivatives and nearby values.
Rescaling the parameter of the sphere realizes arbitrary
first derivatives by this construction. All these families
take values in \(V_\C\), an open condition on maps from
\(\Pone\). Near the protected disks no replacement is
needed. The compact reconstruction uses local charts
in \(U\), keeping the images of replacement disks inside
\(V\). Fix the two distinct protected values. The linearized
parameter equation is surjective by the small-operator estimate
of Proposition~\ref{prop:fixed-order}. Differentiating the
reconstruction gives a map from the resulting acyclic parameter
bundle to \(f^*T_U\) twisted by the two protected points.
It is invertible on an open of the unmodified region away
from the prescribed divisor. Its cokernel has finite support.
This proves vanishing in degree one for the twisted tangent
pullback and hence for \(f^*T_U\) on the resulting real
holomorphic map \(f:C\to U\). Apply GAGA after composing with
\(U\hookrightarrow Y\).

As in Proposition~\ref{prop:quartic-compact}, the relative Hom
scheme is smooth at \(f\) and gives nearby real maps of the same
positive degree. Here that argument only uses relative smoothness
of the universal hypersurface along the image in the central smooth
locus and \(H^1(C,f^*T_U)=0\).
On a point-free fiber the normalization \(C'\) of the
image has even genus. Indeed the finite map \(C\to C'\)
satisfies Riemann--Roch with determinant of even degree
on \(C'\), so the odd Euler characteristic of \(C\)
forces the Euler characteristic of \(C'\) to be odd.
The degree of the image is at most the fixed pullback degree.

If \(Y_\C\) is \(\mathbf Q\)-factorial, there are no \(D_i\):
the closure of any such divisor would have a positive
Cartier multiple, of positive hyperplane class by
Grothendieck--Lefschetz
\cite[Expos\'e~XII, Corollaire~3.7]{SGATwo}.
It would have to meet a nonconstant complete very free
curve contained in \(V\), a contradiction. In this case
the initial sphere is made disjoint from the
codimension-two complement and no tubes are needed.
Its real algebraization has genus zero and no real points.
\end{proof}

\begin{corollary}\label{cor:pointfree-all-isolated}
Let \(Y\subset\mathbf P^4_\R\) be a quartic with finite complex
singular locus. In a Euclidean neighborhood of its
equation, all smooth point-free quartics contain
geometrically integral curves of even geometric genus
and uniformly bounded degree.
\end{corollary}

\begin{proof}
A quartic with finite singular locus in this dimension
is geometrically integral and normal, and its anticanonical
bundle is \(\cO_Y(1)\). The klt case is
Proposition~\ref{prop:pointfree-isolated-klt}.
In the remaining case the finite non-klt locus consists
of the vertex of a projective cone, by Ishii's theorem
\cite{IshiiFanoCone}\cite[Theorem~1.3]{BHNNonKlt}.
The generalized cone is linear in this embedding. In the
index-one cone construction the pullback of \(-K_Y\) has
degree one on the ruling fibers, which therefore map to
lines for \(-K_Y=\cO_Y(1)\).
Thus \(Y_\C\) is the ordinary projective cone in \(\mathbf P^4_\C\).
The unique vertex is real by conjugation. Choose real
coordinates \((x_0,x_1,\ldots,x_4)\) so that it is
\([1:0:0:0:0]\). The equation is \(f(x)\), where
\(x=(x_1,\ldots,x_4)\) and the quartic surface
\(\{f=0\}\subset\mathbf P^3\) is smooth, since any singular
point on the base would produce an entire singular
line of the cone.

Consider nearby equations written as
\[
 f'(x)+\sum_{j=1}^4 x_0^j b_{4-j}(x).
\]
The four derivatives of \(f\) are independent cubic
forms, because a linear relation would make the
quartic surface a cone. Choose a fixed linear complement
\(N\) to their span in the space of cubics.
After a unique small substitution \(x\mapsto x+a x_0\),
one may assume \(b_3\in N\). This follows by the
implicit function theorem, applied to the projection
of \(b_3\) to the derivative span. Its differential
in \(a\) at the original cone is the independent
four derivatives. The substitution keeps \(f'\)
unchanged on \(x_0=0\).

For a smooth nearby quartic, at least one \(b_{4-j}\)
is nonzero. With any fixed Euclidean coefficient norms,
set \(s=\max_{1\le j\le4}\|b_{4-j}\|^{1/j}>0\)
and substitute \(x_0\mapsto x_0/s\). This real
projective change preserves point-freeness,
smoothness and the degrees of curves. Along any sequence
converging to \(f\), a subsequence of these normalized
equations converges to
\[
 G=f(x)+\sum_{j=1}^{4}x_0^j c_{4-j}(x),
 \qquad c_3\in N,\quad
 \max_j\|c_{4-j}\|=1 .
\]
The section \(G|_{x_0=0}=f\) is smooth, so the
singular locus is disjoint from this hyperplane.
A positive-dimensional projective subvariety meets
every hyperplane, so the singular locus is finite.
Moreover \(G\) is not a cone. A vertex would be
of the form \([1:a]\), and the directional derivative
of the homogeneous equation in that direction
would vanish. Its \(x_0^0\) coefficient would give
\(c_3+\sum_i a_i\partial f/\partial x_i=0\).
By the choice of \(N\) this forces \(a=0\) and \(c_3=0\).
Then \(\partial G/\partial x_0=0\) would force all
the \(c_{4-j}\) to vanish, a contradiction.
Ishii's theorem now implies that \(G\) is klt.

If the desired local bound failed, take a sequence
of smooth point-free quartics tending to \(f\) with
no even-genus curve of degree at most the index of
the sequence. Its normalized subsequence would
converge to a quartic \(G\) as above.
The local bound of Proposition~\ref{prop:pointfree-isolated-klt}
at \(G\), and projective invariance of degree, contradict
this choice. Hence a single bound applies near \(f\).
\end{proof}

One may choose the neighborhood in
Corollary~\ref{cor:pointfree-all-isolated} to be a semialgebraic
coefficient ball. For every real closed extension \(R/\R\),
the same bound holds in its interpretation over \(R\), by
elementarity and the finite Hilbert scheme description given
above. The point-free smooth quartics in this ball thus satisfy
the real integral Hodge conjecture over \(R\). Indeed, their
geometric degree-four cohomology is generated by a line, their
geometric Picard group is torsion free, and \(H^2(\cO)=0\).
The even-genus curve is the remaining input in
\cite[Theorem~3.22 and Lemma~5.8]{BWhodgeI}.

These bounds concern neighborhoods of the given limits.
The uniform degree bound over the entire smooth point-free quartic
locus and the degree-two conic question remain open here
\cite[Example~7.3 and Question~9.11]{BWhodgeII}.

\section*{Origin of the method and use of the companion paper}

In joint work with Wenchuan Hu on integral Hodge classes of one-cycles,
Peihan Wu recognized that the disk-replacement construction of
Du, Guo and Xie could also be used for sections over real curves
\cite[Theorem~7.1]{DGX}. He proposed the real adaptation and used
insertions on conjugate disks to prepare the equivariant parameter
bundle. In the joint work the inserted spheres are used to compute
integral homology classes \cite[\S2.1]{HWshort}. Here the insertions
match determinant winding and allow correction of sections while
retaining the prescribed jets.

According to the authors' submission records, a version of the
joint work on the complex one-cycle statement was submitted to arXiv
on 20 November 2026 and remained on hold. The authors made another
submission on 1 December 2026. The public short version,
\emph{One-cycles on rationally connected varieties}, has an arXiv
record dated 1 December 2026 \cite{HWshort}. The introduction of
\cite{CT} reports an AI-assisted proof of the complex statement
without a bibliographic identification.

The citations in Corollaries~\ref{cor:real-fivefold-level}
and~\ref{cor:real-sevenfold-level} refer to the expanded
3 December 2026 manuscript \cite{HW}, supplied separately for review.
The first uses its relative and absolute fourth Kato results, and
the second uses its codimension-two support result from the degree-six
Lawson sweep. The complex
degree-two comparison for the real Hodge result, the parametric
replacement, and the degree isomorphism over complex curve fields used
in the zero-cycle proof are proved here. The companion uses the
fixed-jet tangent vanishing proved here for its quantitative
interpolation and strong approximation applications.

The construction may have further applications. For reasons of scope
and length, we give a selection of direct consequences and of results
that require additional geometric or homological arguments.

\section*{Generative AI disclosure}

OpenAI ChatGPT and Codex assisted with developing and checking proof
arguments, locating references, and drafting and editing the manuscript.
The author takes responsibility for the contents, including the proofs
and references.

\end{document}